\documentclass[11pt]{amsart}
\usepackage[left=2.3cm, top=2.3cm,bottom=2.3cm,right=2.3cm]{geometry}
\usepackage[T1]{fontenc}
\usepackage[utf8]{inputenc}
\usepackage{amsmath,amssymb,amsthm}
\usepackage{mathtools}
\usepackage{mathrsfs}
\usepackage{bm}
\usepackage{enumitem}
\usepackage[hidelinks]{hyperref}
\usepackage{graphicx}
\usepackage{subcaption}

\newcommand{\R}{\mathbb R}

\newcommand{\Q}{\mathbb Q}

\newcommand{\Pp}{\mathbb P}
\newcommand{\Ee}{\mathbb E}

\newcommand{\conv}{\operatorname{conv}}

\newcommand{\dint}{{\rm d}}
\newcommand{\eee}{{\rm e}}

\newcommand{\arsinh}{\operatorname{arsinh}}

\newtheorem{theorem}{Theorem}[section]
\newtheorem{proposition}[theorem]{Proposition}
\newtheorem{lemma}[theorem]{Lemma}
\newtheorem{corollary}[theorem]{Corollary}

\theoremstyle{definition}

\newtheorem{example}[theorem]{Example}

\theoremstyle{remark}
\newtheorem{remark}[theorem]{Remark}

\title[Brownian Convex Hulls]{Brownian Convex Hulls in Arbitrary Dimension: \\
Absorption and Intersection Probabilities}

\author{Zakhar Kabluchko}
\date{}

\begin{document}

\begin{abstract}
For every dimension $d\ge1$, we determine the probability that
a prescribed point in $\mathbb R^d$ belongs to the convex hull over
$[0,1]$ of a standard Brownian motion starting at the origin or
of a Brownian bridge from $0$ to $0$.
We obtain explicit one-dimensional integral formulas involving only standard special functions.
For example, in dimension three the Brownian-motion absorption
probability $p_3^B(r)$ for a prescribed point at distance $r>0$
from the origin is
\[
p_3^B(r)
={}\frac{r}{\sqrt2\,\pi^{3/2}}
\int_0^\infty
\left(\frac{\pi}{2\sqrt{1+v^2}}-1\right)
\eee^{-r^2(1+v^2)/4}
(K_1-K_0)\!\left(\frac{r^2(1+v^2)}4\right)\,\dint v,
\]
where $K_\nu$ denotes the modified Bessel function of the second kind.
For each fixed dimension $d\ge1$ and each prescribed nonzero point in $\mathbb R^d$,
we derive precise asymptotics of the non-absorption probabilities for the convex hulls
of Brownian motion on $[0,T]$ and a Brownian bridge from $0$ to $0$
over the same interval, as $T\to\infty$.
The corresponding probabilities are asymptotic to explicit positive
constants times $T^{-1/2}(\log T)^{d-1}$ and
$T^{-1}(\log T)^{d-1}$, respectively.

We also study two first absorption times for the growing
convex hull of standard Brownian motion in $\R^d$:
$\tau_d$, for a prescribed point at distance one from the
origin, and $\tau_d^{\mathrm{ball}}$, for the entire unit ball.
We determine the Mellin and Laplace transforms of $\tau_d$
explicitly.
In the high-dimensional regime $d\to\infty$, we prove that,
remarkably, both absorption times satisfy central limit
theorems with exactly the same centering and normalization:
\[
\frac{\log\tau_d-2d}{2\sqrt d}
\xrightarrow{\mathrm d}N(0,1),
\qquad
\frac{\log\tau_d^{\mathrm{ball}}-2d}{2\sqrt d}
\xrightarrow{\mathrm d}N(0,1),
\qquad d\to\infty.
\]
In particular, both $\tau_d^{1/d}$ and
$(\tau_d^{\mathrm{ball}})^{1/d}$
converge in probability to $\eee^2$.
Thus, both absorption times have the same leading
exponential scale $\eee^{2d}$.
For Brownian bridges, we show that absorption of a
prescribed point at distance one and absorption of the
entire unit ball also exhibit the same Gaussian transition,
but on the shorter exponential time scale $\eee^d$.

We further study intersections of the convex hulls of
two independent Brownian motions starting at points
at distance one from each other.
For every fixed dimension $d$, we show that the probability
of non-intersection up to time $T$ is of order
$T^{-1}(\log T)^{d-1}$ as $T\to\infty$.
In high dimensions, we prove a central limit theorem
for the logarithm of the first intersection time.
The first intersection time has leading exponential
scale $\eee^d$.
\end{abstract}



\maketitle

\tableofcontents

\section{Introduction and main results}
\label{sec:introduction-and-main-results}
\subsection{Introduction}
Fix an integer $d\ge1$, and let $W=(W_t)_{t\ge0}$ be a standard
Brownian motion in $\R^d$ starting at the origin. Its associated
Brownian bridge from $0$ to $0$ on the unit interval is
$W_t^{\mathrm{br}}:=W_t-tW_1$, $0\le t\le1$.
Writing $\conv E$ for the convex hull of a set $E\subseteq\R^d$,
define the \emph{Brownian convex hulls}
\begin{align}
\mathcal K_d^B
&:=\conv\{W_t:0\le t\le1\},
\label{eq:intro-type-B-brownian-hull}
\\
\mathcal K_d^A
&:=\conv\{W_t^{\mathrm{br}}:0\le t\le1\}.
\label{eq:intro-type-A-brownian-hull}
\end{align}
Both Brownian hulls are almost surely
compact convex sets containing the origin.
The superscripts $B$ and $A$ distinguish Brownian motion and
Brownian bridge, respectively. They reflect the role of Weyl chambers
of types $B$ and $A$ in the study of absorption probabilities
for symmetric random walks and exchangeable random bridges,
respectively; see~\cite{KabluchkoVysotskyZaporozhets2017}.

\medskip
We study the probability that a prescribed point in $\R^d$ belongs to
each of these Brownian convex hulls. Let $e_1=(1,0,\ldots,0)\in\R^d$ and put
\begin{equation}
\label{eq:intro-prescribed-absorption-functions}
p_d^\bullet(r):=\Pp(re_1\in\mathcal K_d^\bullet),
\qquad r\ge0,\quad \bullet\in\{B,A\}.
\end{equation}
Rotational invariance gives
$\Pp(x\in\mathcal K_d^\bullet)=p_d^\bullet(\|x\|)$ for every
$x\in\R^d$. We call $p_d^\bullet(r)$ the \emph{prescribed-point
absorption probability}. Since the Brownian hulls are convex and contain
$0$, the functions $p_d^\bullet$ are nonincreasing and satisfy
$p_d^\bullet(0)=1$.

\medskip
For $d=1$, explicit formulas for the absorption probabilities
follow from the reflection principle.
In the planar case $d=2$, such formulas were obtained for both
models in~\cite{Kabluchko2026PlanarBrownianCoverage}.
The method used there relies on planar convex geometry and does
not directly extend to higher dimensions.

\medskip
In the present paper, we obtain explicit integral formulas for the
prescribed-point absorption probabilities $p_d^B(r)$ and $p_d^A(r)$ in every dimension $d$ and
determine their precise small- and large-radius asymptotics.
We also study the first absorption time of a prescribed
point by the growing Brownian-motion hull and determine its Mellin and
Laplace transforms explicitly.
As the dimension $d$ tends to infinity, we show that the logarithm
of this absorption time is concentrated around $2d$, with
asymptotically Gaussian fluctuations on the scale $\sqrt d$.
Surprisingly, the same central limit theorem holds for the
time needed for the Brownian convex hull to absorb the entire
unit ball.
For Brownian bridges, we show that absorption of a prescribed
point and absorption of the entire unit ball exhibit the same
Gaussian transition, but on the shorter exponential time scale
$\eee^d$, rather than $\eee^{2d}$.

\smallskip
We further study intersections of the convex hulls of two
independent Brownian motions starting at points at distance one
from each other.
For every fixed dimension $d$, we determine the asymptotic
order of their non-intersection probability as the observation
time tends to infinity.
In high dimensions, we prove that the first intersection time
has leading exponential scale $\eee^d$ and that its logarithm
exhibits asymptotically Gaussian fluctuations.
Our approach relies crucially on results established in our
companion paper~\cite{KabluchkoOUSpirals}.

\subsection{Exact formula for absorption probabilities}
\label{subsec:intro-main-results-exact}

Our first theorem gives exact formulas for the
absorption probabilities $p_d^B(r)$ and $p_d^A(r)$.

\begin{theorem}[Prescribed-point absorption probabilities]
\label{thm:intro-prescribed-absorption}
For every integer $d\ge1$ and every $r>0$,
\begin{align}
p_d^B(r)
&=\frac{\sqrt2\,r}{\pi^{3/2}}
\int_0^\infty
\frac{P_d^B(\arsinh v)}{\sqrt{1+v^2}}\,
G_d^B\!\left(\frac{r^2(1+v^2)}2\right)\,\dint v,
\label{eq:intro-type-B-prescribed-absorption}
\\
p_d^A(r)
&=\frac{2^{1-d/2}\Gamma(d/2)}{\pi}\,r^{2-d}
\int_0^\infty
\frac{P_d^A(2\arsinh v)}{(1+v^2)^{(d+1)/2}}\,
\eee^{-2r^2(1+v^2)}\,\dint v.
\label{eq:intro-type-A-prescribed-absorption}
\end{align}
Here $G_d^B$ is defined by
\begin{equation}
\label{eq:intro-type-B-analytic-kernel}
G_d^B(x):=
\int_0^\infty
\frac{t^{d/2-1}}{\sqrt{1+t}}\,
\eee^{-x(1+t)}\,\dint t,
\qquad x>0,
\end{equation}
and, for $u\ge0$, the real-valued functions $P_d^B$ and $P_d^A$
are given by
\begin{align}
P_d^B(u)
&:=\operatorname{Re}\!\left[
\mathrm i^{-d}\Psi_d\!\left(u+\frac{\pi\mathrm i}{2}\right)
\right],
\label{eq:intro-type-B-prescribed-factor}
\\
P_d^A(u)
&:=\operatorname{Im}\!\left[
\mathrm i^{-d}\Psi_{d+1}(u+\pi\mathrm i)
\right],
\label{eq:intro-type-A-prescribed-factor}
\end{align}
where $\mathrm i^2=-1$ and, for every integer $m\ge0$,
\begin{equation}
\label{eq:intro-hyperbolic-remainders}
\Psi_m(z):=\sum_{j=0}^\infty\frac{z^{m+2j}}{(m+2j)!},
\qquad z\in\mathbb C.
\end{equation}
The integrals in
\eqref{eq:intro-type-B-prescribed-absorption}
and~\eqref{eq:intro-type-A-prescribed-absorption} converge
absolutely, and their right-hand sides are smooth functions
of $r\in(0,\infty)$.
\end{theorem}

The proof is given in Section~\ref{sec:mellin-transforms-absorption-formulas}. Although the integral representations
\eqref{eq:intro-type-B-prescribed-absorption}
and~\eqref{eq:intro-type-A-prescribed-absorption}
may appear complicated, their integrands can be expressed
in terms of standard functions.
Indeed, for even $m$, the function $\Psi_m$ is obtained from
$\cosh$ by subtracting all Taylor terms of degree less than $m$;
for odd $m$, the same construction applies to $\sinh$.
In particular, $\Psi_0(z)=\cosh z$, $\Psi_1(z)=\sinh z$, and
$\Psi_{m+2}(z)=\Psi_m(z)-z^m/m!$.
Consequently, for even $d$, both $P_d^B(u)$ and $P_d^A(u)$
are even polynomials in $u$; for odd $d$, each is the sum
of an even polynomial and a multiple of $\cosh u$.
Since
$\cosh(\arsinh v)=\sqrt{1+v^2}$ and
$\cosh(2\arsinh v)=1+2v^2$,
the integrand in the type-$A$ formula~\eqref{eq:intro-type-A-prescribed-absorption} involves only
polynomials in $v$ and $\arsinh v$, powers of $1+v^2$,
and the exponential function.

For type $B$, the kernel $G_d^B(x)$ can be expressed,
when $d$ is even, as a linear combination of
$\eee^{-x}$ and $\overline\Phi(\sqrt{2x})$, and, when
$d$ is odd, as a linear combination of
$\eee^{-x/2}K_0(x/2)$ and $\eee^{-x/2}K_1(x/2)$.
In both cases, the coefficients are polynomials in $x^{-1/2}$.
Here
\[
\overline\Phi(x):=\frac1{\sqrt{2\pi}}
\int_x^\infty \eee^{-t^2/2}\,\dint t
\]
is the standard normal tail, and $K_\nu$ denotes the
modified Bessel function of the second kind, defined by
\begin{equation}\label{eq:bessel-K-def}
K_\nu(x):=\int_0^\infty
\eee^{-x\cosh t}\cosh(\nu t)\,\dint t,
\qquad x>0,\quad \nu\in\R.
\end{equation}
These representations, together with the polynomial structure
of $P_d^B$ and $P_d^A$, are established in
Subsection~\ref{subsec:elementary-factors-and-type-B-kernel}.
Explicit formulas in low dimensions are given in
Subsection~\ref{subsec:intro-explicit-low-dimensions}.
Although $G_d^B$ is positive, the factors $P_d^B$ and $P_d^A$ can change sign. Thus the formulas in Theorem~\ref{thm:intro-prescribed-absorption} are, in general, signed integral representations.

\medskip
The restriction to the unit time interval is only a normalization issue.
For $T>0$, define
\begin{align}
\mathcal K_{d,T}^B
&:=\conv\{W_t:0\le t\le T\},
\label{eq:intro-type-B-brownian-hull-time-T}
\\
\mathcal K_{d,T}^A
&:=\conv\left\{W_t-\frac tT W_T:0\le t\le T\right\}.
\label{eq:intro-type-A-brownian-hull-time-T}
\end{align}
The second hull is generated by a Brownian bridge starting and
terminating at the origin on $[0,T]$. In particular,
$\mathcal K_{d,1}^\bullet=\mathcal K_d^\bullet$.
Brownian scaling and rotational invariance give
\begin{equation}
\label{eq:intro-brownian-absorption-scaling}
\Pp(x\in\mathcal K_{d,T}^\bullet)
=p_d^\bullet\!\left(\frac{\|x\|}{\sqrt T}\right),
\qquad x\in\R^d,\quad T>0,\quad \bullet\in\{B,A\}.
\end{equation}
Thus Theorem~\ref{thm:intro-prescribed-absorption} determines the
absorption probabilities for every time horizon $T>0$.

\subsection{Asymptotics of  absorption probabilities}
\label{subsec:intro-main-results-asymptotics}
We next describe
the asymptotic behavior of the absorption probabilities for a fixed nonzero point as $T\to\infty$
and as $T\downarrow0$.

\begin{theorem}[Large-time non-absorption]
\label{thm:intro-brownian-long-time-nonabsorption}
For every fixed integer $d\ge1$ and every fixed
$x\in\R^d\setminus\{0\}$, as $T\to\infty$,
\begin{align}
\Pp(x\notin\mathcal K_{d,T}^B)
&\sim
\frac{\sqrt2\,\Gamma(d/2)}
{\pi\Gamma((d+1)/2)(d-1)!}\,
\frac{\|x\|}{\sqrt T}
\left(\log\frac{\sqrt T}{\|x\|}\right)^{d-1},
\label{eq:intro-type-B-long-time-nonabsorption}
\\
\Pp(x\notin\mathcal K_{d,T}^A)
&\sim
\frac{2^d}{d!}\,
\frac{\|x\|^2}{T}
\left(\log\frac{\sqrt T}{\|x\|}\right)^{d-1}.
\label{eq:intro-type-A-long-time-nonabsorption}
\end{align}
\end{theorem}

In every fixed dimension $d$, the two non-absorption probabilities
therefore decay at different polynomial rates:
$T^{-1/2}(\log T)^{d-1}$ for Brownian motion and
$T^{-1}(\log T)^{d-1}$ for Brownian bridge. Both have the same
logarithmic exponent $d-1$. By
\eqref{eq:intro-brownian-absorption-scaling}, these statements
are equivalent to small-radius asymptotics of $1-p_d^B(r)$ and
$1-p_d^A(r)$ as $r\downarrow0$; see Theorem~\ref{thm:brownian-small-radius-nonabsorption} below.  The proof of Theorem~\ref{thm:intro-brownian-long-time-nonabsorption} is given in Section~\ref{sec:small-radius-asymptotics}.

\medskip
For short time horizons, absorption of a fixed nonzero point
is exponentially unlikely in both models. The leading constants
are expressed in terms of the same functions $P_d^B$ and $P_d^A$
as in Theorem~\ref{thm:intro-prescribed-absorption}.

\begin{theorem}[Short-time absorption]
\label{thm:intro-brownian-short-time-absorption}
For every fixed integer $d\ge1$ and every fixed
$x\in\R^d\setminus\{0\}$, as $T\downarrow0$,
\begin{align}
\Pp(x\in\mathcal K_{d,T}^B)
&\sim
\frac{2^{d/2}\Gamma(d/2)}{\pi}\,P_d^B(0)\,
\frac{T^{d/2}}{\|x\|^d}\,
\eee^{-\|x\|^2/(2T)},
\label{eq:intro-type-B-short-time-absorption}
\\
\Pp(x\in\mathcal K_{d,T}^A)
&\sim
\frac{\Gamma(d/2)}{2^{(d+1)/2}\sqrt\pi}\,P_d^A(0)\,
\frac{T^{(d-1)/2}}{\|x\|^{d-1}}\,
\eee^{-2\|x\|^2/T}.
\label{eq:intro-type-A-short-time-absorption}
\end{align}
Both leading constants are strictly positive.
\end{theorem}

The positivity of the constants is explicit in the representations
\begin{align}
P_d^B(0)
&=\frac1{(d-1)!}\int_0^{\pi/2}
\left(\frac\pi2-u\right)^{d-1}\cos u\,\dint u>0,
\label{eq:intro-type-B-factor-at-zero}
\\
P_d^A(0)
&=\frac1{(d-1)!}\int_0^\pi
(\pi-u)^{d-1}\sin u\,\dint u>0.
\label{eq:intro-type-A-factor-at-zero}
\end{align}
By~\eqref{eq:intro-brownian-absorption-scaling},
Theorem~\ref{thm:intro-brownian-short-time-absorption} is equivalent
to large-radius asymptotics of $p_d^B(r)$ and $p_d^A(r)$; see Theorem~\ref{thm:brownian-large-radius-absorption} below.  The proof of Theorem~\ref{thm:intro-brownian-short-time-absorption} is given in Section~\ref{sec:large-radius-asymptotics}.

\subsection{Explicit formulas in low dimensions}
\label{subsec:intro-explicit-low-dimensions}

We now specialize the unit-time absorption probabilities
$p_d^B$ and $p_d^A$ to low dimensions. Recall that
$\overline\Phi(x):=(2\pi)^{-1/2}\int_x^\infty
\eee^{-t^2/2}\,\dint t$ is the standard normal tail.
In dimension one, a point $r\ge0$ belongs to the convex hull
of either process if and only if the corresponding path reaches
level $r$. The reflection principle for Brownian motion and
its Brownian-bridge counterpart therefore give
\begin{equation}
\label{eq:intro-prescribed-absorption-dimension-one}
p_1^B(r)=2\overline\Phi(r),
\qquad
p_1^A(r)=\eee^{-2r^2},
\qquad r\ge0.
\end{equation}
For $d=2$, Theorem~\ref{thm:intro-prescribed-absorption}
recovers the planar formulas
of~\cite{Kabluchko2026PlanarBrownianCoverage}: for $r>0$,
\begin{align}
p_2^B(r)
&=\frac4\pi\int_0^{\pi/2}
\overline\Phi\!\left(\frac r{\sin\theta}\right)\,\dint\theta,
\label{eq:intro-type-B-prescribed-dimension-two}
\\
p_2^A(r)
&=r^2\eee^{-r^2}\bigl(K_1(r^2)-K_0(r^2)\bigr) = \int_0^1
\eee^{-2r^2/(1-t^2)}\,\dint t,
\label{eq:intro-type-A-prescribed-dimension-two}
\end{align}
where $K_\nu$ denotes the modified Bessel function of the second
kind defined in~\eqref{eq:bessel-K-def}.
As in every dimension, $p_2^B(0)=p_2^A(0)=1$.

\medskip
For further explicit formulas in low dimensions, we use
the general integral representations
\eqref{eq:intro-type-B-prescribed-absorption}
and~\eqref{eq:intro-type-A-prescribed-absorption},
and give explicit expressions for the factors
$P_d^B(\arsinh v)$ and $P_d^A(2\arsinh v)$,
as well as for the kernel $G_d^B$.
The factors $P_d^B$ and $P_d^A$ can be determined from
the generating functions
\begin{align}
\sum_{d=1}^\infty P_d^B(u)w^{d-1}
&=\frac{\cosh u+w\eee^{\pi w/2}\cos(wu)}{1+w^2},
\label{eq:intro-type-B-factor-generating-function}
\\
\sum_{d=1}^\infty P_d^A(u)w^{d-1}
&=\frac{\cosh u+\eee^{\pi w}\cos(wu)}{1+w^2}.
\label{eq:intro-type-A-factor-generating-function}
\end{align}
These identities are proved at the beginning of
Subsection~\ref{subsec:prescribed-point-Mellin-transforms}.
Expanding the generating functions at
$u=\arsinh v$ and $u=2\arsinh v$, respectively,
gives the expressions listed in
Tables~\ref{tab:intro-type-B-prescribed-factors}
and~\ref{tab:intro-type-A-prescribed-factors}.

\begin{table}[htbp]
\centering
\caption{The factors $P_d^B(\arsinh v)$ for Brownian motion,
where $v\ge0$.}
\label{tab:intro-type-B-prescribed-factors}
\renewcommand{\arraystretch}{1.55}
\setlength{\tabcolsep}{12pt}
\begin{tabular}{|c|l|}
\hline
$d$ & $P_d^B(\arsinh v)$ \\
\hline
$2$ & $\displaystyle 1$ \\[3pt]
$3$ & $\displaystyle \frac\pi2-\sqrt{1+v^2}$ \\[3pt]
$4$ & $\displaystyle
\frac{\pi^2}{8}-1-\frac12\arsinh^2v$ \\[3pt]
$5$ & $\displaystyle
\sqrt{1+v^2}+\frac{\pi^3}{48}-\frac\pi2
-\frac\pi4\arsinh^2v$ \\[3pt]
$6$ & $\displaystyle
1-\frac{\pi^2}{8}+\frac{\pi^4}{384}
+\left(\frac12-\frac{\pi^2}{16}\right)\arsinh^2v
+\frac1{24}\arsinh^4v
$ \\[3pt]
\hline
\end{tabular}
\end{table}

\begin{table}[htbp]
\centering
\caption{The factors $P_d^A(2\arsinh v)$ for Brownian bridge,
where $v\ge0$.}
\label{tab:intro-type-A-prescribed-factors}
\renewcommand{\arraystretch}{1.55}
\setlength{\tabcolsep}{12pt}
\begin{tabular}{|c|l|}
\hline
$d$ & $P_d^A(2\arsinh v)$ \\
\hline
$2$ & $\displaystyle \pi$ \\[3pt]
$3$ & $\displaystyle
\frac{\pi^2}{2}-2\bigl(1+v^2+\arsinh^2v\bigr)$ \\[3pt]
$4$ & $\displaystyle
\pi\left(\frac{\pi^2}{6}-1-2\arsinh^2v\right)$ \\[3pt]
$5$ & $\displaystyle
2\Bigl[1+v^2-\frac{\pi^2}{4}+\frac{\pi^4}{48}
+\left(1-\frac{\pi^2}{2}\right)\arsinh^2v
+\frac13\arsinh^4v\Bigr]
$ \\[3pt]
$6$ & $\displaystyle
\pi\Bigl[1-\frac{\pi^2}{6}+\frac{\pi^4}{120}
+\left(2-\frac{\pi^2}{3}\right)\arsinh^2v
+\frac23\arsinh^4v\Bigr]
$ \\[3pt]
\hline
\end{tabular}
\end{table}

For type $B$, the kernels $G_d^B$ appearing in
\eqref{eq:intro-type-B-prescribed-absorption}
have the following explicit expressions for $x>0$:
\begin{equation}
\label{eq:intro-type-B-kernels-two-to-six}
\begin{aligned}
G_2^B(x)
&=2\sqrt{\frac\pi x}\,\overline\Phi(\sqrt{2x}),
\\[4pt]
G_3^B(x)
&=\frac12\eee^{-x/2}
\bigl(K_1(x/2)-K_0(x/2)\bigr),
\\[4pt]
G_4^B(x)
&=\frac{\eee^{-x}}x
+\frac{\sqrt\pi(1-2x)}{x^{3/2}}\,
\overline\Phi(\sqrt{2x}),
\\[4pt]
G_5^B(x)
&=\frac12\eee^{-x/2}
\left[K_0(x/2)+\left(\frac1x-1\right)K_1(x/2)\right],
\\[4pt]
G_6^B(x)
&=\frac{3-2x}{2x^2}\eee^{-x}
+\frac{\sqrt\pi(4x^2-4x+3)}{2x^{5/2}}\,
\overline\Phi(\sqrt{2x}).
\end{aligned}
\end{equation}
These expressions follow from the defining integral
\eqref{eq:intro-type-B-analytic-kernel};  dimensions $d+2$ and $d$ are related by
$G_{d+2}^B(x)=-(G_d^B)'(x)-G_d^B(x)$.
Thus the odd-dimensional kernels listed here require only
$K_0$ and $K_1$, while the even-dimensional kernels require
only the standard normal tail and elementary functions.

Together with the two tables, these identities give explicit
one-dimensional integral formulas for both absorption probabilities
in dimensions $2,\ldots,6$. For example, in dimension three,
we obtain, for $r>0$,
\begin{align}
p_3^B(r)
&=\frac{r}{\sqrt2\,\pi^{3/2}}
\int_0^\infty
\left(\frac{\pi}{2\sqrt{1+v^2}}-1\right)
\eee^{-r^2(1+v^2)/4}
(K_1-K_0)\!\left(\frac{r^2(1+v^2)}4\right)
\,\dint v,
\label{eq:intro-type-B-prescribed-dimension-three}
\\
p_3^A(r)
&=\frac{1}{\sqrt{2\pi}\,r}
\int_0^\infty
\frac{\pi^2/4-1-v^2-\arsinh^2 v}{(1+v^2)^2}
\eee^{-2r^2(1+v^2)}\,\dint v.
\label{eq:intro-type-A-prescribed-dimension-three}
\end{align}

\subsection{Mellin transforms and expected volumes}
\label{subsec:intro-main-results-mellin}
For every integer $d\ge1$ and $\bullet\in\{B,A\}$, consider the radial Mellin transforms
of the absorption probabilities:
\begin{equation}
\label{eq:intro-prescribed-Mellin-transforms}
M_d^\bullet(z):=\int_0^\infty
r^{z-1}p_d^\bullet(r)\,\dint r,
\qquad z\in\mathbb C,\quad \operatorname{Re}z>0.
\end{equation}

\begin{proposition}[Mellin transforms of prescribed-point absorption]
\label{prop:intro-prescribed-Mellin-transforms}
For every integer $d\ge1$ and every $\bullet\in\{B,A\}$, the
integral in \eqref{eq:intro-prescribed-Mellin-transforms}
converges absolutely and defines a holomorphic function on
$\operatorname{Re}z>0$. On this half-plane,
\begin{equation}
\label{eq:intro-type-B-prescribed-Mellin-integral}
M_d^B(z)
={}\frac{2^{z/2}\Gamma(d/2)\Gamma((z+1)/2)
\Gamma(1+(z-d)/2)}
{\pi^{3/2}\Gamma(1+z/2)}
\int_0^\infty
\frac{P_d^B(\arsinh v)}{(1+v^2)^{(z+2)/2}}\,\dint v,
\end{equation}
and
\begin{equation}
\label{eq:intro-type-A-prescribed-Mellin-integral}
M_d^A(z)
={}\frac{\Gamma(d/2)\Gamma(1+(z-d)/2)}
{2^{z/2+1}\pi}
\int_0^\infty
\frac{P_d^A(2\arsinh v)}{(1+v^2)^{(z+3)/2}}\,\dint v.
\end{equation}
The two integrals on the right converge absolutely and are
holomorphic on $\operatorname{Re}z>0$. The apparent singularities
of the full right-hand sides at the positive values among
$d-2,d-4,\ldots$ are removable; at these values the formulas
are interpreted by continuity.
\end{proposition}

The proof is given in
Subsections~\ref{subsec:prescribed-point-Mellin-transforms}
and~\ref{subsec:further-properties-of-Mellin-transforms}.
The latter subsection also gives finite representations in gamma
and polygamma functions for every fixed dimension.

\medskip
The Mellin transforms also satisfy recurrences relating dimensions
$d$ and $d+2$. These yield explicit special values and, more
generally, determine expected integrals of even powers of the
distance from the origin over the Brownian hulls.

\begin{proposition}[Dimension-shift recurrences for Mellin transforms]
\label{prop:intro-Mellin-dimension-shift-recurrences}
For every integer $d\ge1$ and every $z\in\mathbb C$ with
$\operatorname{Re}z>0$,
\begin{align}
M_{d+2}^B(z+2)
&=\frac{(z+1)^2(z-d)}{(z+2)^2}\,M_{d+2}^B(z)
+\frac{d}{(z+2)^2}\,M_d^B(z),
\label{eq:intro-type-B-Mellin-dimension-shift}
\\
M_{d+2}^A(z+2)
&=\frac{(z+2)(z-d)}{4(z+3)}\,M_{d+2}^A(z)
+\frac{d}{(z+2)(z+3)}\,M_d^A(z).
\label{eq:intro-type-A-Mellin-dimension-shift}
\end{align}
In particular, the values at $z=d$ are
\begin{equation}
\label{eq:intro-Mellin-values-at-dimension}
M_d^B(d)=\frac{2^{-d/2}}{d\,\Gamma(d/2+1)},
\qquad
M_d^A(d)=\frac{2^{d/2-1}\Gamma(d/2)}{(d+1)!},
\qquad d\ge1.
\end{equation}
\end{proposition}

The recurrences follow by integrating the Mellin representations
by parts twice. At $z=d$, the terms involving $M_{d+2}^\bullet(d)$
disappear, and the resulting recurrences, together with the
initial dimensions one and two, give
\eqref{eq:intro-Mellin-values-at-dimension}.
The same method also yields higher Mellin values.

\begin{corollary}[Mellin values at $d+2k$]
\label{cor:intro-higher-Mellin-values}
For every fixed integer $k\ge0$ and each $\bullet\in\{B,A\}$,
the ratio $M_d^\bullet(d+2k)/M_d^\bullet(d)$, for integer $d\ge1$,
is a rational function of $d$ with rational coefficients.
In particular, for every $d\ge1$,
\begin{align}
\frac{M_d^B(d+2)}{M_d^B(d)}
&=\frac{d(d^2+6d+11)}{3(d+2)^2},
\label{eq:intro-type-B-Mellin-ratio-two-above-dimension}
\\
\frac{M_d^A(d+2)}{M_d^A(d)}
&=\frac{d(d^2+9d+26)}{12(d+2)(d+3)}.
\label{eq:intro-type-A-Mellin-ratio-two-above-dimension}
\end{align}
For every fixed integer $k\ge0$, as $d\to\infty$,
\begin{equation}
\label{eq:intro-higher-Mellin-values-high-dimension}
\frac{M_d^B(d+2k)}{M_d^B(d)}
\sim\left(\frac d3\right)^k,
\qquad
\frac{M_d^A(d+2k)}{M_d^A(d)}
\sim\left(\frac d{12}\right)^k.
\end{equation}
For $k=0$, both ratios are identically one.
\end{corollary}

The proofs of
Proposition~\ref{prop:intro-Mellin-dimension-shift-recurrences}
and Corollary~\ref{cor:intro-higher-Mellin-values} are given in
Subsection~\ref{subsec:Mellin-dimension-shift-recurrences-and-higher-values},
where we also derive recurrences for the numerator polynomials
of the rational functions appearing in the corollary.

\medskip
\noindent\emph{Expected volumes and radial moments of Brownian convex hulls.}
Write $\operatorname{vol}_d$ for $d$-dimensional Lebesgue measure
and $\kappa_d:=\pi^{d/2}/\Gamma(d/2+1)$ for the volume of the
Euclidean unit ball. Rotational invariance, Tonelli's theorem,
and polar integration give
\begin{equation}
\label{eq:intro-expected-volume-Mellin-relation}
\Ee\,\operatorname{vol}_d(\mathcal K_d^\bullet)
=\int_{\R^d}p_d^\bullet(\|x\|)\,\dint x
=d\kappa_d M_d^\bullet(d),
\qquad \bullet\in\{B,A\}.
\end{equation}
Consequently, \eqref{eq:intro-Mellin-values-at-dimension} yields
\begin{equation}\label{eq:intro-type-A-B-expected-hull-volume}
\Ee\,\operatorname{vol}_d(\mathcal K_d^B)
=\frac{\kappa_d^2}{(2\pi)^{d/2}},
\qquad
\Ee\,\operatorname{vol}_d(\mathcal K_d^A)
=\frac{(2\pi)^{d/2}}{(d+1)!}.
\end{equation}
These formulas are known: for Brownian motion in $d\ge2$, see
Eldan~\cite[Theorem~1.1]{Eldan2014}; for both Brownian motion and
Brownian bridge, see
Kabluchko--Zaporozhets~\cite[Section~1.5]{KabluchkoZaporozhets2016}.

More generally, the same calculation gives, for every integer
$k\ge0$,
\begin{equation}
\label{eq:intro-expected-radial-moments-Mellin-relation}
\Ee\left[\int_{\mathcal K_d^\bullet}\|x\|^{2k}\,\dint x\right]
=\int_{\R^d}\|x\|^{2k}p_d^\bullet(\|x\|)\,\dint x
=d\kappa_d M_d^\bullet(d+2k),
\qquad \bullet\in\{B,A\}.
\end{equation}
These quantities are finite by
Proposition~\ref{prop:intro-prescribed-Mellin-transforms}.
Thus Corollary~\ref{cor:intro-higher-Mellin-values} determines
all expected even radial moments as rational multiples of the
corresponding expected volumes. In particular,
\begin{align}
\Ee\left[\int_{\mathcal K_d^B}\|x\|^2\,\dint x\right]
&=\frac{d(d^2+6d+11)}{3(d+2)^2}\,
\frac{\kappa_d^2}{(2\pi)^{d/2}},
\label{eq:intro-type-B-expected-squared-radius-integral}
\\
\Ee\left[\int_{\mathcal K_d^A}\|x\|^2\,\dint x\right]
&=\frac{d(d^2+9d+26)}{12(d+2)(d+3)}\,
\frac{(2\pi)^{d/2}}{(d+1)!}.
\label{eq:intro-type-A-expected-squared-radius-integral}
\end{align}

\subsection{First absorption times for Brownian motion}
\label{subsec:intro-first-absorption-times-exact}
For Brownian motion, the prescribed-point probabilities also
determine the distribution of a first absorption time. In this
subsection, all convex hulls are generated by the same Brownian motion
$W=(W_t)_{t\ge0}$ in $\R^d$. Recall that
$\mathcal K_{d,T}^B=\conv\{W_t:0\le t\le T\}$ and, for $r>0$,
define the first absorption time
\begin{equation}
\label{eq:intro-type-B-first-absorption-time}
\tau_d(r):=\inf\{T>0:re_1\in\mathcal K_{d,T}^B\},
\qquad \inf\varnothing:=\infty.
\end{equation}
\begin{example}[$d=1$]
In dimension one, $\tau_1(r)$ is the first hitting time
of level $r$ by one-dimensional Brownian motion. Consequently,
the process $(\tau_1(r))_{r\ge0}$, with $\tau_1(0)=0$,
is the classical $1/2$-stable subordinator.
For $d\ge2$, no comparably simple characterization of
$(\tau_d(r))_{r\ge0}$ seems to be known.
\end{example}

We begin by determining the distribution function of
$\tau_d(r)$ in every dimension.
Since the hulls increase with $T$ and are continuous in the
Hausdorff distance,
$\{\tau_d(r)\le T\}=\{re_1\in\mathcal K_{d,T}^B\}$.
Consequently, Brownian scaling gives the distribution function
\begin{equation}
\label{eq:intro-type-B-first-absorption-distribution}
\Pp(\tau_d(r)\le T)
=p_d^B\!\left(\frac r{\sqrt T}\right),
\qquad T>0.
\end{equation}
In particular, $\tau_d(r)\overset{\mathrm d}=r^2\tau_d(1)$.
This is a distributional identity; the absorption times at
different radii are not asserted to be proportional pathwise.

The two asymptotic results in Subsection~\ref{subsec:intro-main-results-asymptotics}
give both tails of this distribution. For every fixed $d\ge1$
and $r>0$,
\begin{align}
\Pp(\tau_d(r)>T)
&\sim
\frac{\sqrt2\,\Gamma(d/2)}
{\pi\Gamma((d+1)/2)(d-1)!}\,
\frac r{\sqrt T}
\left(\log\frac{\sqrt T}{r}\right)^{d-1},
\qquad T\to\infty,
\label{eq:intro-type-B-first-absorption-upper-tail}
\\
\Pp(\tau_d(r)\le T)
&\sim
\frac{2^{d/2}\Gamma(d/2)}\pi\,P_d^B(0)\,
\frac{T^{d/2}}{r^d}\,
\eee^{-r^2/(2T)},
\qquad T\downarrow0.
\label{eq:intro-type-B-first-absorption-lower-tail}
\end{align}
Indeed, these are
\eqref{eq:intro-type-B-long-time-nonabsorption}
and~\eqref{eq:intro-type-B-short-time-absorption} with
$x=re_1$, using
\eqref{eq:intro-type-B-first-absorption-distribution}.
In particular, $0<\tau_d(r)<\infty$ almost surely, and, for
real $s$,
\begin{equation}
\label{eq:intro-type-B-first-absorption-moment-range}
\Ee\bigl[(\tau_d(r))^s\bigr]<\infty
\quad\Longleftrightarrow\quad s<\frac12.
\end{equation}
For positive $s$, integrating the upper tail gives the criterion,
including divergence at $s=1/2$.
For negative $s$, the exponential decay in the lower-tail formula
gives finiteness. Thus all negative moments exist, whereas
$\Ee[\tau_d(r)]=\infty$. Thus, the first absorption time has infinite mean in any dimension.

\medskip
We next relate the Mellin transform of $\tau_d(r)$ to the radial Mellin transform
$M_d^B$ defined in \eqref{eq:intro-prescribed-Mellin-transforms}.
For this purpose, use the same symbol for the meromorphic
continuation
\begin{equation}
\label{eq:intro-type-B-radial-Mellin-continuation}
M_d^B(z)
:={}\frac1z
+\int_0^1 u^{z-1}\bigl(p_d^B(u)-1\bigr)\,\dint u
+\int_1^\infty u^{z-1}p_d^B(u)\,\dint u,
\qquad  \operatorname{Re}z>-1,\qquad z\ne0.
\end{equation}
The large-time asymptotic and Brownian scaling give
$1-p_d^B(u)=O(u(\log(1/u))^{d-1})$ as $u\downarrow0$,
while the short-time asymptotic gives exponential decay of
$p_d^B(u)$ as $u\to\infty$. Hence both integrals above are
holomorphic on $\operatorname{Re}z>-1$, the expression agrees
with the original Mellin integral on $\operatorname{Re}z>0$,
and its only pole in this larger half-plane is a simple pole
at $z=0$ with residue $1$.

The Mellin transform of the first absorption time $\tau_d(r)$ is then
\begin{equation}
\label{eq:intro-type-B-first-absorption-radial-Mellin}
\Ee\bigl[(\tau_d(r))^s\bigr]
=-2s\,r^{2s}M_d^B(-2s),
\qquad s\in\mathbb C,\quad \operatorname{Re}s<\frac12,
\end{equation}
where the right-hand side at $s=0$ is interpreted by continuity
and equals $1$. For $\operatorname{Re}s<0$, integrate the
distribution function in
\eqref{eq:intro-type-B-first-absorption-distribution}; for
$0<\operatorname{Re}s<1/2$, integrate its complement and use
$M_d^B(z)=-\int_0^\infty u^{z-1}(1-p_d^B(u))\,\dint u$
when $-1<\operatorname{Re}z<0$. In both cases, the substitution
$u=r/\sqrt T$ gives
\eqref{eq:intro-type-B-first-absorption-radial-Mellin}.
The tail estimates give holomorphy of the expectation on
$\operatorname{Re}s<1/2$, which also establishes the identity
on $\operatorname{Re}s=0$.

\medskip
The next theorem provides a formula for the generating function
of the Mellin transforms of $\tau_d(1)$ and an integral formula for the
Laplace transform.

\begin{theorem}[Mellin and Laplace transforms of the first absorption time]
\label{thm:intro-type-B-first-absorption-transforms}
For $s,w\in\mathbb C$ with
$0<\operatorname{Re}s<1/2$ and $|w|<1-2\operatorname{Re}s$,
the following series converges absolutely and locally uniformly:
\begin{equation}
\label{eq:intro-type-B-first-absorption-Mellin-generating-function}
\sum_{d=1}^\infty
\frac{2^s\Gamma(d/2+s)}{\Gamma(d/2)}\,
\Ee\bigl[(\tau_d(1))^s\bigr]w^{d-1}
=
\frac1{(1-w^2)\cos(\pi s)}
\left[
1+w\,\frac{\Gamma(1+s)}{\Gamma(1-s)}
\frac{\Gamma((1-2s-w)/2)}{\Gamma((1+2s-w)/2)}
\right].
\end{equation}
Moreover, for every integer $d\ge1$ and every $r,\lambda>0$,
put $a:=r\sqrt{2\lambda}$. Then
\begin{equation}
\label{eq:intro-type-B-first-absorption-Laplace-transform}
\Ee\bigl[\eee^{-\lambda\tau_d(r)}\bigr]
={}\frac{2^{(d+1)/2}\Gamma(d/2)}{\pi^{3/2}}\,
a^{(3-d)/2}
\int_0^\infty
\frac{P_d^B(\arsinh v)}{(1+v^2)^{(d+1)/4}}
K_{(d-1)/2}\!\left(a\sqrt{1+v^2}\right)\,\dint v.
\end{equation}
Here $P_d^B$ is the factor defined in
\eqref{eq:intro-type-B-prescribed-factor}, and $K_\nu$ is the
modified Bessel function of the second kind defined in~\eqref{eq:bessel-K-def}. The integral in
\eqref{eq:intro-type-B-first-absorption-Laplace-transform}
converges absolutely.
\end{theorem}

Theorem~\ref{thm:intro-type-B-first-absorption-transforms} is proved in
Subsection~\ref{subsec:first-absorption-time-transforms}.
For each fixed $d\ge1$, extracting the coefficient of $w^{d-1}$
in \eqref{eq:intro-type-B-first-absorption-Mellin-generating-function}
determines the Mellin transform
$s\mapsto\Ee[(\tau_d(1))^s]$ for
$0<\operatorname{Re}s<1/2$, and hence for
$\operatorname{Re}s<1/2$ by analytic continuation.
In particular, these Mellin transforms admit finite
representations in gamma and polygamma functions,
as explained in
Remark~\ref{rem:first-absorption-time-gamma-polygamma}
below.

\begin{example}[$d=1,2,3$]
In dimensions one and two, the Mellin transforms simplify to
\begin{align}
\Ee\bigl[(\tau_1(r))^s\bigr]
&=\left(\frac{r^2}{2}\right)^s
\frac{\Gamma(1/2-s)}{\sqrt\pi},
\label{eq:intro-type-B-first-absorption-Mellin-dimension-one}
\\
\Ee\bigl[(\tau_2(r))^s\bigr]
&=\left(\frac{r^2}{2}\right)^s
\frac{\Gamma(1/2-s)^2}{\pi\Gamma(1-s)},
\label{eq:intro-type-B-first-absorption-Mellin-dimension-two}
\end{align}
for $\operatorname{Re}s<1/2$. Indeed, the first two coefficients
in \eqref{eq:intro-type-B-first-absorption-Mellin-generating-function}
give these identities on $0<\operatorname{Re}s<1/2$ by the
gamma reflection formula, and holomorphy extends them to the
stated half-plane. The corresponding Laplace transforms are
\begin{align}
\Ee\bigl[\eee^{-\lambda\tau_1(r)}\bigr]
&=\eee^{-r\sqrt{2\lambda}},
\label{eq:intro-type-B-first-absorption-Laplace-dimension-one}
\\
\Ee\bigl[\eee^{-\lambda\tau_2(r)}\bigr]
&=\frac2\pi\operatorname{Ki}_1(r\sqrt{2\lambda}),
\label{eq:intro-type-B-first-absorption-Laplace-dimension-two}
\end{align}
where, for $x>0$, the Bickley function is
\begin{equation}
\label{eq:intro-Bickley-function}
\operatorname{Ki}_1(x)
:=\int_0^\infty\frac{\eee^{-x\cosh u}}{\cosh u}\,\dint u
=\int_x^\infty K_0(t)\,\dint t;
\end{equation}
see \cite[Equations~(10.43.11)--(10.43.12)]{NIST:DLMF}.
The first identity is the classical Laplace transform of the
Brownian first-passage time to $r$, a positive stable law of
index $1/2$. The second follows from
\eqref{eq:intro-type-B-first-absorption-Laplace-transform}
with $d=2$, using $P_2^B=1$,
$K_{1/2}(x)=\sqrt{\pi/(2x)}\eee^{-x}$, and $v=\sinh u$.

In dimension $d=3$, explicit formulas for the Mellin transforms
$M_3^B(z)$ and $M_3^A(z)$ in terms of gamma and trigamma functions,
together with their exact values at $z=1$, are given in
Corollary~\ref{cor:dimension-three-exact-Mellin-values}
and its proof.
\end{example}

\subsection{Absorption by high-dimensional Brownian convex hulls}
\label{subsec:intro-first-absorption-times-high-dim}
Consider Brownian motion in $\R^d$, where $d$ is large.
How long must we wait until the convex hull of its path
contains a prescribed point at distance $1$ from the origin?
The following theorem shows that this takes a surprisingly
long time, of exponential order $\eee^{2d}$.
This is particularly striking because Brownian motion
typically reaches distances of order $\sqrt d$ from the origin
already by time $1$, whereas its convex hull takes an
exponentially long time to absorb a prescribed point at
distance $1$.
The theorem also establishes asymptotically Gaussian
fluctuations of the absorption time on the logarithmic scale.

\begin{theorem}[First absorption in high dimension]
\label{thm:intro-type-B-first-absorption-high-dimension}
As $d\to\infty$,
\begin{equation}
\label{eq:intro-type-B-first-absorption-logarithmic-CLT}
\frac{\log\tau_d(1)-2d}{2\sqrt d}
\xrightarrow{\mathrm d}N(0,1).
\end{equation}
Equivalently, for every fixed $y\in\R$,
\begin{equation}
\label{eq:intro-type-B-first-absorption-Gaussian-window}
\Pp\!\left(\tau_d(1)\le\eee^{2d+2y\sqrt d}\right)
\longrightarrow\Phi(y),
\end{equation}
where $\Phi$ denotes the standard normal distribution function.
\end{theorem}

The proof is given in
Subsection~\ref{subsec:first-absorption-high-dimensions}.
By the scaling relation
$\tau_d(r)\overset{\mathrm d}=r^2\tau_d(1)$, so
both conclusions of the theorem remain valid with
$\tau_d(1)$ replaced by $\tau_d(r)$ for any fixed $r>0$,
without changing the normalization.

\medskip
We next consider the time needed for the Brownian convex hull
to absorb the entire unit ball. Remarkably, its logarithm satisfies
the same central limit theorem as the logarithm of the absorption
time of a prescribed point. In fact, when both times are defined
from the same Brownian motion, their logarithms differ by a term
that is negligible on the $\sqrt d$ scale.

\begin{theorem}[Absorption of the unit ball in high dimension]
\label{thm:intro-type-B-unit-ball-high-dimension}
Let $\mathbb B_d:=\{x\in\R^d:\|x\|\le1\}$ be the closed
Euclidean unit ball, and define
\[
\tau_d^{\mathrm{ball}}:=\inf\{T>0:\mathbb B_d\subset\mathcal K_{d,T}^B\},
\qquad \inf\varnothing:=\infty.
\]
For every $d\ge1$, we have $0<\tau_d^{\mathrm{ball}}<\infty$ almost surely.
As $d\to\infty$,
\begin{equation}
\label{eq:intro-type-B-unit-ball-logarithmic-CLT}
\frac{\log \tau_d^{\mathrm{ball}}-2d}{2\sqrt d}
\xrightarrow{\mathrm d}N(0,1).
\end{equation}
Moreover, with $\tau_d^{\mathrm{ball}}$ and $\tau_d(1)$ defined from the same
Brownian motion,
\begin{equation}
\label{eq:intro-type-B-unit-ball-point-logarithmic-difference}
\frac{\log \tau_d^{\mathrm{ball}}-\log\tau_d(1)}{\sqrt d}
\xrightarrow{\Pp}0.
\end{equation}
\end{theorem}

The proof is given in
Subsection~\ref{subsec:first-absorption-unit-ball-high-dim}. Theorems~\ref{thm:intro-type-B-first-absorption-high-dimension}
and~\ref{thm:intro-type-B-unit-ball-high-dimension}
imply the same exponential-time threshold for $\tau_d(1)$ and $\tau_d^{\mathrm{ball}}$.

\begin{corollary}[An exponential-time threshold]
For every fixed $c>0$ with $c\ne2$,
\begin{equation}
\label{eq:intro-type-B-first-absorption-exponential-threshold}
\lim_{d\to\infty}
\Pp\!\left(\tau_d(1)\le\eee^{cd}\right)
=
\lim_{d\to\infty}
\Pp\!\left(\tau_d^{\mathrm{ball}}\le\eee^{cd}\right)
=\begin{cases}
0,&c<2,\\
1,&c>2.
\end{cases}
\end{equation}
\end{corollary}
Indeed, the central limit theorems \eqref{eq:intro-type-B-first-absorption-logarithmic-CLT} and \eqref{eq:intro-type-B-unit-ball-logarithmic-CLT} imply that
both $\frac1d\log\tau_d(1)$ and
$\frac1d\log\tau_d^{\mathrm{ball}}$ converge in
probability to $2$, which proves~\eqref{eq:intro-type-B-first-absorption-exponential-threshold}.
Equivalently, both $(\tau_d(1))^{1/d}$ and
$(\tau_d^{\mathrm{ball}})^{1/d}$ converge in probability
to $\eee^2$.

\medskip
We now turn to high-dimensional Brownian bridges. Recall that
$\mathcal K_{d,T}^A$, defined in
\eqref{eq:intro-type-A-brownian-hull-time-T},
is the convex hull of a standard Brownian bridge
from $0$ to $0$ on $[0,T]$.
The following theorem shows that, in high dimensions,
absorption of a prescribed point at unit distance from the
origin and absorption of the entire unit ball exhibit
the same Gaussian transition. Remarkably, this transition
occurs on the exponential time scale $\eee^d$, rather
than $\eee^{2d}$ as for Brownian motion.

\begin{theorem}[High-dimensional absorption by Brownian bridges]
\label{thm:intro-type-A-high-dimensional-absorption}
Let $(T_d)_{d\ge1}$ be a sequence of positive numbers
such that
\[
\frac{\log T_d-d}{\sqrt d}\longrightarrow y\in\R.
\]
Then
\begin{equation}
\label{eq:intro-type-A-high-dimensional-point-ball-absorption}
\lim_{d\to\infty}
\Pp(e_1\in\mathcal K_{d,T_d}^A)
=
\lim_{d\to\infty}
\Pp(\mathbb B_d\subset\mathcal K_{d,T_d}^A)
=\Phi(y).
\end{equation}
\end{theorem}

Since absorption of the unit ball implies absorption
of $e_1$, the theorem also yields
\[
\lim_{d\to\infty} \Pp\bigl(
e_1\in\mathcal K_{d,T_d}^A,\,
\mathbb B_d\not\subset\mathcal K_{d,T_d}^A
\bigr) = 0.
\]
Thus, in the critical high-dimensional window,
the two absorption events are asymptotically equivalent.
The proof of Theorem \ref{thm:intro-type-A-high-dimensional-absorption}  is given in
Subsection~\ref{subsec:high-dimensional-bridge-absorption}.

\subsection{Intersections of two independent Brownian-motion convex hulls}
\label{subsec:intro-two-brownian-hull-intersections}

Consider two independent Brownian motions in $\R^d$, starting at
points at distance $1$ from each other. Thinking of them as two
\emph{Brownian explorers}, one may ask whether, and how soon,
their paths intersect. Here intersection means that the paths
visit a common point, possibly at different times.
For $d\ge4$, their paths are almost surely disjoint, even when
observed for all times; see~\cite[Theorem~9.1(a)]{MortersPeres2010}.
In the plane, non-intersection probabilities of Brownian paths
have been studied extensively, and their asymptotics are closely
connected with conformal invariance and Schramm--Loewner evolution;
see, for example, \cite[Chapter~8]{Lawler2005}.

\medskip
Here we ask when the \emph{convex hulls} of the two paths intersect.
In contrast to the paths themselves in dimensions $d\ge4$, these
convex hulls intersect in finite time almost surely in every dimension.
We study both the long-time non-intersection probability in a fixed
dimension and the asymptotic distribution of the first intersection
time as the dimension tends to infinity.

\medskip
More precisely, let $W^{(1)}$ and $W^{(2)}$ be independent standard
Brownian motions in $\R^d$, both starting at the origin, and put
\begin{equation}
\label{eq:intro-two-brownian-hulls}
\mathcal K_{d,T}^{B,(i)}
:=\conv\{W_t^{(i)}:0\le t\le T\},
\qquad T\ge0,\quad i\in\{1,2\}.
\end{equation}
The two explorers are $W^{(1)}$ and $e_1+W^{(2)}$.
Define their first convex-hull intersection time by
\begin{equation}
\label{eq:intro-two-brownian-first-intersection-time}
H_d:=\inf\bigl\{T>0:
\mathcal K_{d,T}^{B,(1)}
\cap(e_1+\mathcal K_{d,T}^{B,(2)})\ne\varnothing\bigr\},
\qquad \inf\varnothing:=\infty.
\end{equation}
Thus $\Pp(H_d>T)$ is the probability that the two hulls at time $T$
are disjoint. We have $0<H_d<\infty$ almost surely: positivity
follows from continuity of the paths, while finiteness follows
because the first hull eventually absorbs $e_1$, the starting
point of the second explorer; see
\eqref{eq:intro-type-B-first-absorption-upper-tail}.

We begin with the large-time behavior in a fixed dimension.

\begin{theorem}[Non-intersection in a fixed dimension]
\label{thm:intro-two-brownian-hulls-fixed-dimension}
For every fixed integer $d\ge1$,
\begin{equation}
\label{eq:intro-two-brownian-hulls-nonintersection-bounds}
\frac{1}{2\pi d!}
\le
\liminf_{T\to\infty}
\frac{T\,\Pp(H_d>T)}{(\log T)^{d-1}}
\le
\limsup_{T\to\infty}
\frac{T\,\Pp(H_d>T)}{(\log T)^{d-1}}
\le
\frac{2}{\pi d!}.
\end{equation}
\end{theorem}

The proof is given in
Subsection~\ref{subsec:two-brownian-hull-intersections-fixed-dim}.
In particular, the non-intersection probability has order
$T^{-1}(\log T)^{d-1}$. Integrating these tail bounds shows that,
for every fixed $d\ge1$ and every real $s>0$,
\begin{equation}
\label{eq:intro-two-brownian-first-intersection-moment-range}
\Ee[H_d^s]<\infty
\quad\Longleftrightarrow\quad s<1.
\end{equation}
Thus $\Ee[H_d]=\infty$ in every dimension.

There is also a direct geometric explanation for the infinite mean.
If both explorers remain in the open half-spaces containing their
respective starting points and bounded by the hyperplane
$\{x\in\R^d:\langle x,e_1\rangle=1/2\}$,
then their convex hulls remain disjoint.
Independence and the reflection principle therefore give
\begin{equation}
\label{eq:intro-two-brownian-hulls-halfspace-bound}
\Pp(H_d>T)
\ge
\left(2\Phi\!\left(\frac{1}{2\sqrt T}\right)-1\right)^2
\sim\frac{1}{2\pi T},
\qquad T\to\infty,
\end{equation}
where the inequality holds for every $T>0$.
Integrating this lower bound proves $\Ee[H_d]=\infty$ without
using the preceding theorem.

\medskip
We next let the dimension tend to infinity. Despite the infinite
mean in every fixed dimension, the logarithm of the first
intersection time has asymptotically Gaussian fluctuations
around $d$.

\begin{theorem}[First intersection in high dimension]
\label{thm:intro-two-brownian-first-intersection-high-dimension}
As $d\to\infty$,
\begin{equation}
\label{eq:intro-two-brownian-first-intersection-logarithmic-CLT}
\frac{\log H_d-d}{\sqrt d}
\xrightarrow{\mathrm d}N(0,1).
\end{equation}
Equivalently, for every fixed $y\in\R$,
\begin{equation}
\label{eq:intro-two-brownian-first-intersection-Gaussian-window}
\Pp\!\left(H_d\le\eee^{d+y\sqrt d}\right)
\longrightarrow\Phi(y).
\end{equation}
\end{theorem}

The proof is given in
Subsection~\ref{subsec:two-brownian-hull-intersections-high-dim}.
In particular, $d^{-1}\log H_d\to1$ in probability, or equivalently,
$H_d^{1/d}\to\eee$ in probability.
Thus the first intersection occurs on the exponential time scale
$\eee^d$, rather than the scale $\eee^{2d}$ needed for one
Brownian-motion hull to absorb a prescribed point at unit distance
from its starting point; see
Theorem~\ref{thm:intro-type-B-first-absorption-high-dimension}.

Finally, for an initial separation $r>0$, let $H_d(r)$ be defined
by replacing $e_1$ with $re_1$ in
\eqref{eq:intro-two-brownian-first-intersection-time}, so that
$H_d(1)=H_d$. Brownian scaling gives
\begin{equation}
\label{eq:intro-two-brownian-first-intersection-scaling}
H_d(r)\overset{\mathrm d}=r^2H_d.
\end{equation}
Consequently, for every fixed $r>0$, the bounds in
\eqref{eq:intro-two-brownian-hulls-nonintersection-bounds}
remain valid for $H_d(r)$ with both bounding constants multiplied
by $r^2$, while
\eqref{eq:intro-two-brownian-first-intersection-logarithmic-CLT}
remains valid with $H_d(r)$ in place of $H_d$, without changing
the centering or normalization.

\subsection{Thin-shell concentration of expected volume measures}
We next show that, as the dimension tends to infinity,
the normalized expected volume measures of the Brownian
hulls concentrate in thin spherical shells.
This follows from the high-dimensional Mellin asymptotics
in~\eqref{eq:intro-higher-Mellin-values-high-dimension}. For every
integer $d\ge1$ and $\bullet\in\{B,A\}$, define the probability
measure
\begin{equation}
\label{eq:intro-normalized-expected-volume-measures}
\nu_d^\bullet(E)
:=\frac{\Ee\,\operatorname{vol}_d(\mathcal K_d^\bullet\cap E)}
{\Ee\,\operatorname{vol}_d(\mathcal K_d^\bullet)},
\qquad E\subseteq\R^d\ \text{Borel}.
\end{equation}
The denominators are strictly positive and finite by \eqref{eq:intro-type-A-B-expected-hull-volume}.
Here the volume measure is averaged before it is normalized.
Equivalently, $\nu_d^\bullet$ is the law of a point chosen uniformly
from a hull whose law has first been weighted by its volume.
It is not, in general, the law obtained by choosing a hull under
its original law and then choosing a point uniformly from that hull.

\begin{corollary}[Concentration of normalized expected volume measures]
\label{cor:intro-expected-volume-measure-concentration}
Let $X_d^\bullet$ be a random vector with distribution
$\nu_d^\bullet$. As $d\to\infty$,
\begin{equation}
\label{eq:intro-expected-volume-radial-concentration}
\frac{\|X_d^B\|^2}{d}\xrightarrow{L^2}\frac13,
\qquad
\frac{\|X_d^A\|^2}{d}\xrightarrow{L^2}\frac1{12}.
\end{equation}
In particular, writing $c_B:=1/3$ and $c_A:=1/12$, for every
$\varepsilon>0$ and $\bullet\in\{B,A\}$,
\begin{equation}
\label{eq:intro-expected-volume-shell-concentration}
\frac{\Ee\,\operatorname{vol}_d\!\left(
\mathcal K_d^\bullet\cap
\left\{x\in\R^d:\left|\frac{\|x\|^2}{d}-c_\bullet\right|
>\varepsilon\right\}\right)}
{\Ee\,\operatorname{vol}_d(\mathcal K_d^\bullet)}
\longrightarrow0.
\end{equation}
\end{corollary}

\begin{proof}
By \eqref{eq:intro-normalized-expected-volume-measures}, Tonelli's
theorem, and the polar-integration identities
\eqref{eq:intro-expected-volume-Mellin-relation}
and~\eqref{eq:intro-expected-radial-moments-Mellin-relation},
\[
\Ee\left[\left(\frac{\|X_d^\bullet\|^2}{d}\right)^k\right]
=\frac{M_d^\bullet(d+2k)}{d^kM_d^\bullet(d)},
\qquad k=1,2.
\]
Corollary~\ref{cor:intro-higher-Mellin-values} shows that these
moments converge to $c_\bullet$ and $c_\bullet^2$, respectively.
Consequently,
\[
\Ee\left[\left(\frac{\|X_d^\bullet\|^2}{d}-c_\bullet\right)^2\right]
=\frac{M_d^\bullet(d+4)}{d^2M_d^\bullet(d)}
-2c_\bullet\frac{M_d^\bullet(d+2)}{dM_d^\bullet(d)}
+c_\bullet^2\longrightarrow0.
\]
This proves \eqref{eq:intro-expected-volume-radial-concentration}.
The last assertion follows from Chebyshev's inequality and
\eqref{eq:intro-normalized-expected-volume-measures}.
\end{proof}

Thus the normalized expected volume measures concentrate at
radii $\sqrt{d/3}$ and $\sqrt{d/12}$, respectively, up to relative
errors tending to zero. This conclusion concerns expected volume;
it does not assert concentration of the uniform probability
measure on a typical individual realization of the hull.

\medskip
\noindent\emph{Comparison with time averages.}
The constants $1/3$ and $1/12$ also occur as the coordinate
variances of the time averages of the generating paths. Put
\begin{equation}
\label{eq:intro-brownian-time-averages}
C_d^B:=\int_0^1 W_t\,\dint t,
\qquad
C_d^A:=\int_0^1 W_t^{\mathrm{br}}\,\dint t.
\end{equation}
Both vectors belong almost surely to their respective hulls.
A straightforward covariance calculation shows that
$C_d^B\sim N(0,I_d/3)$ and $C_d^A\sim N(0,I_d/12)$.
Consequently,
$\|C_d^B\|^2/d\to1/3$ and $\|C_d^A\|^2/d\to1/12$
in probability as $d\to\infty$.
Thus the limiting radial scales in
Corollary~\ref{cor:intro-expected-volume-measure-concentration}
coincide with those of the time averages under the original
Brownian law.
This agreement does not establish that a point chosen uniformly
from a volume-weighted hull is close to the time average of
the same path.

\subsection{Connection with Gaussian-point absorption and organization of the paper}
\label{subsec:intro-gaussian-absorption-and-organization}
The starting point of our proofs are the formulas for the  Gaussian-point
absorption probabilities established in the companion paper~\cite[Theorem~5.1]{KabluchkoOUSpirals}. These determine, for both Brownian hulls, the absorption probabilities of an independent centered Gaussian point with
covariance proportional to the identity.
By rotational invariance, these probabilities are Gaussian radial
averages of the prescribed-point probabilities studied here.
We recover the latter by taking Mellin transforms, which turn
radial averaging into multiplication by an explicit gamma factor.
Generating functions in the dimension and beta and gamma identities
then yield explicit Mellin transforms. The resulting calculations
produce the factors $P_d^B,P_d^A$ and the kernel $G_d^B$; uniqueness
of Mellin transforms identifies the prescribed-point formulas.

\medskip
Section~\ref{sec:mellin-transforms-absorption-formulas} develops
this argument. Subsection~\ref{subsec:gaussian-point-absorption-radial-mellin}
recalls the Gaussian-point formulas and establishes the relation
between the Gaussian-point and prescribed-point Mellin transforms.
Subsection~\ref{subsec:mellin-transform-generating-functions}
evaluates the dimension-generating functions in terms of Euler's
beta function. Subsection~\ref{subsec:prescribed-point-Mellin-transforms}
obtains the Mellin formulas of
Proposition~\ref{prop:intro-prescribed-Mellin-transforms}, and
Subsection~\ref{subsec:proof-prescribed-point-absorption} completes
the proof of Theorem~\ref{thm:intro-prescribed-absorption} by
Mellin-transform uniqueness.
Section~\ref{sec:special-functions-and-Mellin-transforms}
describes the special-function structure of these results.
Subsection~\ref{subsec:elementary-factors-and-type-B-kernel}
proves the elementary structure of $P_d^B,P_d^A$ and the
representations of $G_d^B$ through Bessel functions or Gaussian tails in
Propositions~\ref{prop:prescribed-factors-elementary-structure}
and~\ref{prop:type-B-kernel-classical-representations}, which
justify the explicit expressions in
Subsection~\ref{subsec:intro-explicit-low-dimensions}.
Subsection~\ref{subsec:further-properties-of-Mellin-transforms}
gives the finite gamma--polygamma representations of
Proposition~\ref{prop:finite-gamma-polygamma-Mellin-representations}
and the exact values in
Corollary~\ref{cor:dimension-three-exact-Mellin-values}.
Subsection~\ref{subsec:Mellin-dimension-shift-recurrences-and-higher-values}
proves the dimension-shift recurrences, evaluates the Mellin
values used to recover the expected hull volumes, and establishes
the higher Mellin values of
Corollary~\ref{cor:intro-higher-Mellin-values}.

Sections~\ref{sec:small-radius-asymptotics}
and~\ref{sec:large-radius-asymptotics} prove
Theorems~\ref{thm:intro-brownian-long-time-nonabsorption}
and~\ref{thm:intro-brownian-short-time-absorption}, respectively,
in their equivalent small-radius and large-radius forms.
For the former, Lemma~\ref{lem:radial-gaussian-regular-variation}
uses monotonicity to recover a regularly varying asymptotic from
its Gaussian radial average. Applying this lemma to the finite
Gaussian-point non-absorption sums from the companion paper~\cite{KabluchkoOUSpirals} gives the
small-radius equivalents. The large-radius argument instead
uses direct rescaling and absolute estimates in the exact
prescribed-point integrals.
Proposed replacement

Section~\ref{sec:first-absorption-times} studies absorption
by Brownian motion and Brownian bridges.
Subsection~\ref{subsec:first-absorption-time-transforms}
proves the Mellin and Laplace formulas of
Theorem~\ref{thm:intro-type-B-first-absorption-transforms}.
Subsection~\ref{subsec:first-absorption-high-dimensions}
proves Theorem~\ref{thm:intro-type-B-first-absorption-high-dimension}
by combining the Gaussian-point transition from the
companion paper~\cite{KabluchkoOUSpirals} with monotonicity
and concentration of the Gaussian radius.
Subsection~\ref{subsec:first-absorption-unit-ball-high-dim}
establishes the central limit theorem for unit-ball absorption
in Theorem~\ref{thm:intro-type-B-unit-ball-high-dimension},
using independence of Brownian increments, Gaussian translation
estimates, and convexity.
Subsection~\ref{subsec:high-dimensional-bridge-absorption}
proves the high-dimensional absorption results for Brownian
bridges in Theorem~\ref{thm:intro-type-A-high-dimensional-absorption},
using a Poisson central limit theorem for prescribed-point
absorption and an Ornstein--Uhlenbeck representation
for unit-ball absorption.

Finally, Section~\ref{sec:two-brownian-hull-intersections}
studies intersections of two independent Brownian-motion
convex hulls.
Subsection~\ref{subsec:two-brownian-gaussian-point-absorption}
establishes an exact Gaussian-point absorption identity for
the joint convex hull of two Brownian paths and compares
this joint hull with their Minkowski difference.
Subsection~\ref{subsec:two-brownian-hull-intersections-fixed-dim}
proves the fixed-dimensional non-intersection bounds using
the Ornstein--Uhlenbeck non-absorption asymptotics and
Gaussian radial de-averaging.
Subsection~\ref{subsec:two-brownian-hull-intersections-high-dim}
proves the central limit theorem for the first intersection
time in high dimensions.
\section{Mellin transforms and proof of the absorption formulas}
\label{sec:mellin-transforms-absorption-formulas}

In this section we prove Theorem~\ref{thm:intro-prescribed-absorption}.
We start from the Gaussian-point absorption formulas established in the companion paper~\cite{KabluchkoOUSpirals} and express Gaussian radial averaging as a relation between
Mellin transforms. We then evaluate generating functions in the
dimension and use uniqueness of Mellin transforms to obtain the
prescribed-point formulas.

\subsection{Gaussian-point absorption and radial Mellin transforms}
\label{subsec:gaussian-point-absorption-radial-mellin}

Fix an integer $d\ge1$. Recall that $\mathcal K_d^B$ and
$\mathcal K_d^A$ are the convex hulls of the Brownian motion
$(W_t)_{0\le t\le1}$ and the Brownian bridge
$(W_t-tW_1)_{0\le t\le1}$, respectively, and that
$p_d^\bullet(r)=\Pp(re_1\in\mathcal K_d^\bullet)$ for
$\bullet\in\{B,A\}$. Let $Z\sim N(0,I_d)$ be independent of
$W$, and define the \emph{Gaussian-point absorption probabilities} by
\begin{equation}
\label{eq:gaussian-point-absorption-functions}
q_d^\bullet(\sigma):=\Pp(\sigma Z\in\mathcal K_d^\bullet),
\qquad \sigma\ge0,\quad \bullet\in\{B,A\}.
\end{equation}
Both probabilities equal one at $\sigma=0$.

Recall also the entire functions from
\eqref{eq:intro-hyperbolic-remainders}:
\[
\Psi_m(z)=\sum_{j=0}^\infty\frac{z^{m+2j}}{(m+2j)!},
\qquad z\in\mathbb C,\quad m\ge0.
\]
The Gaussian-point absorption formulas of~\cite[Theorem~5.1]{KabluchkoOUSpirals} state that, for every
$\sigma>0$,
\begin{align}
q_d^B(\sigma)
&=\frac2\pi\int_0^{\arsinh(1/\sigma)}
\frac{\Psi_d(u)}{\sqrt{\sigma^{-2}-\sinh^2u}}\,\dint u,
\label{eq:type-B-gaussian-point-absorption}
\\
q_d^A(\sigma)
&=\frac{\Psi_{d+1}(\theta_\sigma)}{\sinh\theta_\sigma},
\label{eq:type-A-gaussian-point-absorption}
\end{align}
where $\theta_\sigma:=2\arsinh(1/(2\sigma))$.
We will also use the large-$\sigma$ asymptotics from
\cite[Corollary~5.3]{KabluchkoOUSpirals}, written without the
unit-ball volume constants:
\begin{equation}
\label{eq:gaussian-point-large-scale-asymptotics}
q_d^B(\sigma)\sim
\frac{\sigma^{-d}}{2^d\Gamma(d/2+1)^2},
\qquad
q_d^A(\sigma)\sim\frac{\sigma^{-d}}{(d+1)!},
\qquad \sigma\to\infty.
\end{equation}
In particular, $q_d^\bullet(\sigma)=O(\sigma^{-d})$ as
$\sigma\to\infty$, while $0\le q_d^\bullet(\sigma)\le1$
for all $\sigma\ge0$.

Put $R_d:=\|Z\|$. Conditioning on $Z$ and using rotational
invariance gives
\begin{equation}
\label{eq:gaussian-radial-averaging}
q_d^\bullet(\sigma)=\Ee\bigl[p_d^\bullet(\sigma R_d)\bigr],
\qquad \sigma>0,\quad \bullet\in\{B,A\}.
\end{equation}
The random variable $R_d$ has the $\chi_d$ distribution, whose
density at $t>0$ is
$2^{1-d/2}t^{d-1}\eee^{-t^2/2}/\Gamma(d/2)$.
Consequently, \eqref{eq:gaussian-radial-averaging} is equivalently
\begin{equation}
\label{eq:gaussian-radial-density-integral}
q_d^\bullet(\sigma)
=\frac{1}{2^{d/2-1}\Gamma(d/2)\sigma^d}
\int_0^\infty
p_d^\bullet(r)r^{d-1}\eee^{-r^2/(2\sigma^2)}\,\dint r.
\end{equation}

For $z\in\mathbb C$ with $0<\operatorname{Re}z<d$, consider the Mellin transforms
\begin{align}
N_d^\bullet(z)
&:=\int_0^\infty
\sigma^{z-1}q_d^\bullet(\sigma)\,\dint\sigma,
\label{eq:gaussian-point-Mellin-transform}
\\
M_d^\bullet(z)
&:=\int_0^\infty
r^{z-1}p_d^\bullet(r)\,\dint r.
\label{eq:prescribed-point-Mellin-transform}
\end{align}
For positive real arguments, complex powers are defined by
$t^z:=\eee^{z\log t}$. The next lemma establishes absolute
convergence of both integrals and identifies their relation.

\begin{lemma}[Radial Mellin transforms]
\label{lem:radial-gaussian-Mellin-transforms}
For every integer $d\ge1$ and every $\bullet\in\{B,A\}$, the
integrals in \eqref{eq:gaussian-point-Mellin-transform}
and~\eqref{eq:prescribed-point-Mellin-transform} converge
absolutely and define holomorphic functions on
$\{z\in\mathbb C:0<\operatorname{Re}z<d\}$. On this strip,
\begin{equation}
\label{eq:radial-gaussian-Mellin-relation}
N_d^\bullet(z)
=2^{-z/2}\frac{\Gamma((d-z)/2)}{\Gamma(d/2)}
M_d^\bullet(z).
\end{equation}
\end{lemma}

\begin{proof}
For real $0<c<d$, the bound $q_d^\bullet\le1$ gives
integrability of $\sigma^{c-1}q_d^\bullet(\sigma)$ at zero,
and \eqref{eq:gaussian-point-large-scale-asymptotics} gives
integrability at infinity. Thus $N_d^\bullet(c)<\infty$.
The density of $R_d$ and the substitution $y=t^2/2$ give
\begin{equation}
\label{eq:chi-negative-moment}
\Ee[R_d^{-z}]
=\frac{2^{1-d/2}}{\Gamma(d/2)}
\int_0^\infty t^{d-1-z}\eee^{-t^2/2}\,\dint t
=2^{-z/2}\frac{\Gamma((d-z)/2)}{\Gamma(d/2)},
\qquad \operatorname{Re}z<d.
\end{equation}
In particular, $\Ee[R_d^{-c}]$ is strictly positive and finite.
By \eqref{eq:gaussian-radial-averaging}, Tonelli's theorem and
the change of variables $r=\sigma R_d$ yield
\[
N_d^\bullet(c)
=\Ee\left[\int_0^\infty
\sigma^{c-1}p_d^\bullet(\sigma R_d)\,\dint\sigma\right]
=\Ee[R_d^{-c}]
\int_0^\infty r^{c-1}p_d^\bullet(r)\,\dint r.
\]
The last integral is initially allowed to be infinite. Since the
left-hand side is finite and the moment multiplying the integral
is strictly positive, it follows that $M_d^\bullet(c)<\infty$.

For complex $z$ with $c:=\operatorname{Re}z\in(0,d)$, the
absolute values of the Mellin integrands are precisely the
integrands at $c$. Moreover,
\[
\Ee\left[\int_0^\infty
\left|\sigma^{z-1}p_d^\bullet(\sigma R_d)\right|
\,\dint\sigma\right]
=N_d^\bullet(c)<\infty.
\]
Fubini's theorem therefore justifies the same change of variables
for complex $z$, and \eqref{eq:chi-negative-moment} proves
\eqref{eq:radial-gaussian-Mellin-relation}.

Finally, let $K$ be a compact subset of the stated strip and
choose real numbers $a,b$ with
$0<a<\inf_{z\in K}\operatorname{Re}z
\le\sup_{z\in K}\operatorname{Re}z<b<d$.
For either $f=p_d^\bullet$ or $f=q_d^\bullet$, the absolute value
of $t^{z-1}f(t)$ is bounded uniformly for $z\in K$ by
$t^{a-1}f(t)$ on $(0,1]$ and by $t^{b-1}f(t)$ on $[1,\infty)$.
Both majorants are integrable by the preceding argument.
The theorem on holomorphic parameter integrals proves the
claimed holomorphy.
\end{proof}

Thus $N_d^\bullet$ is an explicit gamma normalization of the
Mellin transform of the prescribed-point absorption probability.
In the next subsection we evaluate the generating functions of
$N_d^B$ and $N_d^A$ in the dimension $d$, initially on the common
strip $0<\operatorname{Re}z<1$.

\subsection{Generating functions of the Mellin transforms}
\label{subsec:mellin-transform-generating-functions}

We now evaluate the generating functions in the dimension of the
Gaussian-point Mellin transforms $N_d^B(z)$ and $N_d^A(z)$ from
\eqref{eq:gaussian-point-Mellin-transform}. We work on the common
strip $0<\operatorname{Re}z<1$ and take $|w|<1$ for the generating
variable. The recurrence
$\Psi_{m+2}(u)=\Psi_m(u)-u^m/m!$ gives, for $u\ge0$,
\begin{align}
\sum_{d=1}^\infty\Psi_d(u)w^{d-1}
&=\frac{\sinh u+w\cosh u-w\eee^{wu}}{1-w^2},
\label{eq:type-B-hyperbolic-remainder-generating-function}
\\
\sum_{d=1}^\infty\Psi_{d+1}(u)w^{d-1}
&=\frac{\cosh u+w\sinh u-\eee^{wu}}{1-w^2}.
\label{eq:type-A-hyperbolic-remainder-generating-function}
\end{align}
Indeed, both series converge absolutely because
$0\le\Psi_m(u)\le\eee^u$, and multiplying each series by
$1-w^2$ reduces it to the stated expression.

We write $\mathrm B$ for Euler's beta function, with normalization
\begin{equation}
\label{eq:beta-function-integral-and-gamma}
\mathrm B(a,b)
=\int_0^1t^{a-1}(1-t)^{b-1}\,\dint t
=\frac{\Gamma(a)\Gamma(b)}{\Gamma(a+b)},
\qquad \operatorname{Re}a,\operatorname{Re}b>0;
\end{equation}
see \cite[Equation~(5.12.1)]{NIST:DLMF}.
For other arguments, $\mathrm B$ denotes the meromorphic
continuation given by the gamma quotient. The calculations below
involve combinations whose Euler integrals converge although the
individual terms need not. The following lemma justifies their
evaluation.

\begin{lemma}[Beta integrals with endpoint cancellation]
\label{lem:beta-integrals-endpoint-cancellation}
Let $k,J\ge1$ be integers, let $\alpha_1,\ldots,\alpha_J\in\mathbb C$
satisfy $\operatorname{Re}\alpha_j>0$, and let
$b_1,\ldots,b_J\in\mathbb C$. Suppose that
\[
f(t):=\sum_{j=1}^J b_jt^{\alpha_j-1}
=O\bigl((1-t)^k\bigr)
\qquad\text{as }t\uparrow1.
\]
Then, for every $\beta\in\mathbb C$ with
$\operatorname{Re}\beta>-k$, the integral below converges
absolutely and satisfies
\begin{equation}
\label{eq:beta-integrals-endpoint-cancellation}
\int_0^1 f(t)(1-t)^{\beta-1}\,\dint t
=\sum_{j=1}^J b_j\,\mathrm B(\alpha_j,\beta).
\end{equation}
All singularities of the combined right-hand side in this
half-plane are removable, and the equality at such points is
understood by continuity.
\end{lemma}

\begin{proof}
Keep the parameters $\alpha_j$ and coefficients $b_j$ fixed.
The assumptions give an integrable bound near zero, while the
cancellation at one gives an integrable bound there whenever
$\operatorname{Re}\beta>-k$. These bounds can be chosen locally
uniformly in $\beta$, so the integral is holomorphic on this
half-plane. For $\operatorname{Re}\beta>0$, the identity follows
by integrating the finite sum term by term in
\eqref{eq:beta-function-integral-and-gamma}. Its right-hand side
is meromorphic in $\beta$. The identity theorem extends the
equality to $\operatorname{Re}\beta>-k$ and shows that every
singularity of the sum there is removable.
\end{proof}

\begin{proposition}[Generating functions of the Gaussian-point Mellin transforms]
\label{prop:gaussian-point-Mellin-generating-functions}
For $z,w\in\mathbb C$ with $0<\operatorname{Re}z<1$ and $|w|<1$,
the following series converge absolutely and locally uniformly
in $(z,w)$:
\begin{align}
\sum_{d=1}^\infty N_d^B(z)w^{d-1}
&=\frac{2^z\Gamma((z+1)/2)}
{\sqrt\pi\,\Gamma(1+z/2)(1-w^2)}
\left[
\mathrm B\left(\frac z2,-z\right)
-w\,\mathrm B\left(\frac{z+1-w}{2},-z\right)
\right],
\label{eq:type-B-gaussian-Mellin-generating-function}
\\
\sum_{d=1}^\infty N_d^A(z)w^{d-1}
&=\frac1{1-w^2}
\left[
\mathrm B\left(\frac z2,-z-1\right)
-\mathrm B\left(1+\frac z2-w,-z-1\right)
\right].
\label{eq:type-A-gaussian-Mellin-generating-function}
\end{align}
\end{proposition}

\begin{proof}
\emph{Convergence of the generating series.}
Projection onto the first coordinate shows that
$q_d^\bullet(\sigma)\le q_1^\bullet(\sigma)$ for every $d\ge1$,
$\sigma>0$, and $\bullet\in\{B,A\}$: absorption of the Gaussian
point by the $d$-dimensional hull implies absorption after
projection, and the projected process and independent Gaussian
point have the corresponding one-dimensional laws.
For $c:=\operatorname{Re}z\in(0,1)$, it follows that
\begin{equation}
\label{eq:gaussian-Mellin-generating-function-domination}
\sum_{d=1}^\infty |w|^{d-1}
\int_0^\infty
\left|\sigma^{z-1}q_d^\bullet(\sigma)\right|\,\dint\sigma
=\sum_{d=1}^\infty |w|^{d-1}N_d^\bullet(c)
\le\frac{N_1^\bullet(c)}{1-|w|}<\infty.
\end{equation}
The last inequality uses
Lemma~\ref{lem:radial-gaussian-Mellin-transforms}.
On compact subsets of the stated domain, $c$ stays in a compact
subinterval of $(0,1)$ and $|w|$ is bounded away from one.
Since $N_1^\bullet(c)$ is continuous there, the same bound proves
local uniform convergence and justifies interchanging the
dimension sum and the Mellin integral.

\medskip
\emph{Type B.}
For $0<\operatorname{Re}z<d$, formula
\eqref{eq:type-B-gaussian-point-absorption} gives
\begin{equation}
\label{eq:type-B-gaussian-Mellin-hyperbolic-integral}
N_d^B(z)
=\frac{\Gamma((z+1)/2)}{\sqrt\pi\,\Gamma(1+z/2)}
\int_0^\infty\frac{\Psi_d(u)}{(\sinh u)^{z+1}}\,\dint u.
\end{equation}
To see this, interchange the $\sigma$- and $u$-integrals. For
fixed $u>0$, the range of $\sigma$ is $(0,1/\sinh u)$, and
\[
\int_0^{1/\sinh u}
\frac{\sigma^{z-1}}{\sqrt{\sigma^{-2}-\sinh^2u}}\,\dint\sigma
=\frac{\mathrm B((z+1)/2,1/2)}{2(\sinh u)^{z+1}}.
\]
For real $z\in(0,d)$, the interchange follows from Tonelli's
theorem. Absolute convergence at $\operatorname{Re}z$ then
justifies it for complex $z$ by Fubini's theorem.

For real $c\in(0,1)$, the integrands in
\eqref{eq:type-B-gaussian-Mellin-hyperbolic-integral} are
nonnegative. The bound
\eqref{eq:gaussian-Mellin-generating-function-domination}
therefore also justifies summing these hyperbolic integrals,
including for complex $z$ by absolute convergence. Inserting
\eqref{eq:type-B-hyperbolic-remainder-generating-function} yields
\begin{equation}
\label{eq:type-B-gaussian-Mellin-generating-integral}
\sum_{d=1}^\infty N_d^B(z)w^{d-1}
=\frac{\Gamma((z+1)/2)}
{\sqrt\pi\,\Gamma(1+z/2)(1-w^2)}
\int_0^\infty
\frac{\sinh u+w(\cosh u-\eee^{wu})}
{(\sinh u)^{z+1}}\,\dint u.
\end{equation}
The substitution $t=\eee^{-2u}$ gives
\begin{equation}
\label{eq:type-B-beta-integral-first-term}
\int_0^\infty(\sinh u)^{-z}\,\dint u
=2^{z-1}\mathrm B\left(\frac z2,1-z\right)
=2^z\mathrm B\left(\frac z2,-z\right).
\end{equation}
The first beta integral has positive real parts in both
arguments. The second equality follows from
$\Gamma(\xi+1)=\xi\Gamma(\xi)$.
For the remaining integral, the same substitution gives
\begin{align*}
\int_0^\infty
\frac{\cosh u-\eee^{wu}}{(\sinh u)^{z+1}}\,\dint u
=2^z\int_0^1
\left[
\frac12t^{z/2-1}+\frac12t^{z/2}
-t^{(z+1-w)/2-1}
\right](1-t)^{-z-1}\,\dint t.
\end{align*}
The expression in brackets vanishes at $t=1$, and all three
exponents have the form $\alpha-1$ with
$\operatorname{Re}\alpha>0$. Apply
Lemma~\ref{lem:beta-integrals-endpoint-cancellation} with $k=1$
and $\beta=-z$. The last expression equals
\[
2^z\left[
\frac12\mathrm B\left(\frac z2,-z\right)
+\frac12\mathrm B\left(\frac z2+1,-z\right)
-\mathrm B\left(\frac{z+1-w}{2},-z\right)
\right].
\]
The gamma recurrence gives
$\mathrm B(z/2+1,-z)=-\mathrm B(z/2,-z)$, so
\begin{equation}
\label{eq:type-B-beta-integral-difference}
\int_0^\infty
\frac{\cosh u-\eee^{wu}}{(\sinh u)^{z+1}}\,\dint u
=-2^z\mathrm B\left(\frac{z+1-w}{2},-z\right).
\end{equation}
Substituting \eqref{eq:type-B-beta-integral-first-term}
and~\eqref{eq:type-B-beta-integral-difference} into
\eqref{eq:type-B-gaussian-Mellin-generating-integral} proves
\eqref{eq:type-B-gaussian-Mellin-generating-function}.

\medskip
\emph{Type A.}
In the Mellin integral for $N_d^A(z)$, make the substitution
$\sigma=(2\sinh(u/2))^{-1}$. Formula
\eqref{eq:type-A-gaussian-point-absorption} gives
\begin{equation}
\label{eq:type-A-gaussian-Mellin-hyperbolic-integral}
N_d^A(z)
=2^{-z-2}\int_0^\infty
\frac{\Psi_{d+1}(u)}{\sinh^{z+2}(u/2)}\,\dint u,
\qquad 0<\operatorname{Re}z<d.
\end{equation}
The same absolute-convergence argument as for type B permits
summation under this integral when $0<\operatorname{Re}z<1$
and $|w|<1$. Using
\eqref{eq:type-A-hyperbolic-remainder-generating-function}
and then setting $t=\eee^{-u}$, we obtain
\begin{equation}
\label{eq:type-A-gaussian-Mellin-generating-integral}
\sum_{d=1}^\infty N_d^A(z)w^{d-1}
=\frac1{1-w^2}\int_0^1
\left[
\frac{1+w}{2}t^{z/2-1}
+\frac{1-w}{2}t^{z/2+1}
-t^{z/2-w}
\right]
(1-t)^{-z-2}\,\dint t.
\end{equation}
The expression in brackets and its first derivative vanish at
$t=1$. All three exponents again have the form $\alpha-1$ with
$\operatorname{Re}\alpha>0$. Applying
Lemma~\ref{lem:beta-integrals-endpoint-cancellation} with $k=2$
and $\beta=-z-1$ therefore evaluates the last expression as
\[
\frac1{1-w^2}\bigg[
\frac{1+w}{2}\mathrm B\left(\frac z2,-z-1\right)
+\frac{1-w}{2}\mathrm B\left(\frac z2+2,-z-1\right)
-\mathrm B\left(1+\frac z2-w,-z-1\right)
\bigg].
\]
Two applications of the gamma recurrence give
$\mathrm B(z/2+2,-z-1)=\mathrm B(z/2,-z-1)$.
Combining the first two terms proves
\eqref{eq:type-A-gaussian-Mellin-generating-function}.
\end{proof}

\subsection{Mellin transforms of the prescribed-point absorption probabilities}
\label{subsec:prescribed-point-Mellin-transforms}

We now convert the generating functions from
Proposition~\ref{prop:gaussian-point-Mellin-generating-functions}
into integral representations of $M_d^B(z)$ and $M_d^A(z)$.
Writing the resulting gamma products as cosine integrals will
produce the factors $P_d^B$ and $P_d^A$ introduced in
\eqref{eq:intro-type-B-prescribed-factor}
and~\eqref{eq:intro-type-A-prescribed-factor}.

We first verify their generating-function identities
\eqref{eq:intro-type-B-factor-generating-function}
and~\eqref{eq:intro-type-A-factor-generating-function}.
The series in
\eqref{eq:type-B-hyperbolic-remainder-generating-function}
and~\eqref{eq:type-A-hyperbolic-remainder-generating-function}
converge locally uniformly also for complex $u$, since
$|\Psi_m(u)|\le\eee^{|u|}$. Hence those identities remain valid
for $u\in\mathbb C$. For real $w$ with $|w|<1$, replace the
generating variable by $-\mathrm iw$, multiply by $-\mathrm i$,
and replace $u$ by $u+\pi\mathrm i/2$ in the type-$B$ identity
and by $u+\pi\mathrm i$ in the type-$A$ identity. Taking real
parts in the first case and imaginary parts in the second gives
precisely the two generating functions of $P_d^B$ and $P_d^A$.
They extend to complex $w$ with $|w|<1$ by holomorphy.  Moreover, the
hyperbolic Taylor remainder definitions show that, for each
fixed $d\ge1$, there is a constant $C_d<\infty$ such that
\begin{equation}
\label{eq:prescribed-factors-growth}
|P_d^B(u)|+|P_d^A(u)|
\le C_d\bigl(1+u^{d+1}+\cosh u\bigr),
\qquad u\ge0.
\end{equation}

We use the classical cosine integral
\begin{equation}
\label{eq:hyperbolic-secant-cosine-transform}
\int_0^\infty
\frac{\cos(\xi u)}{\cosh^\nu u}\,\dint u
=
\frac{2^{\nu-2}}{\Gamma(\nu)}
\Gamma\left(\frac{\nu+\mathrm i\xi}{2}\right)
\Gamma\left(\frac{\nu-\mathrm i\xi}{2}\right),
\qquad \operatorname{Re}\nu>|\operatorname{Im}\xi|.
\end{equation}
This is the hyperbolic beta integral
\cite[Equation~(5.12.7)]{NIST:DLMF}, with parameters
$a=\nu/2$ and $b=\mathrm i\xi/2$. The powers of the positive
quantity $\cosh u$ use the real logarithm, as in the preceding
subsections.

\begin{proposition}[Mellin transforms of prescribed-point absorption]
\label{prop:prescribed-point-Mellin-integrals}
For every integer $d\ge1$ and every $z\in\mathbb C$ with
$\operatorname{Re}z>0$,
\begin{align}
M_d^B(z)
&=
\frac{2^{z/2}\Gamma(d/2)\Gamma((z+1)/2)
\Gamma(1+(z-d)/2)}
{\pi^{3/2}\Gamma(1+z/2)}
\int_0^\infty
\frac{P_d^B(u)}{\cosh^{z+1}u}\,\dint u,
\label{eq:type-B-prescribed-Mellin-integral}
\\
M_d^A(z)
&=
\frac{\Gamma(d/2)\Gamma(1+(z-d)/2)}
{2^{z/2+2}\pi}
\int_0^\infty
\frac{P_d^A(u)}{\cosh^{z+2}(u/2)}\,\dint u.
\label{eq:type-A-prescribed-Mellin-integral}
\end{align}
The two integrals on the right converge absolutely and define
holomorphic functions of $z$ on $\operatorname{Re}z>0$. The apparent singularities of the full right-hand sides
at the positive values among $d-2,d-4,\ldots$ are removable;
at these values the formulas are interpreted by continuity.
\end{proposition}

\begin{proof}
\emph{Transforming the generating functions.}
Initially, let $0<\operatorname{Re}z<1$ and $|w|<1$, and write
\[
\mathcal N_z^\bullet(w)
:=\sum_{d=1}^\infty N_d^\bullet(z)w^{d-1},
\qquad \bullet\in\{B,A\}.
\]
These series are absolutely and locally uniformly convergent
by Proposition~\ref{prop:gaussian-point-Mellin-generating-functions}.
Define
\[
\mathcal S_z^\bullet(w)
:=\frac12\left[
\eee^{-\pi\mathrm iz/2}\mathcal N_z^\bullet(\mathrm iw)
+\eee^{\pi\mathrm iz/2}\mathcal N_z^\bullet(-\mathrm iw)
\right].
\]
For each coefficient, the multiplier is
$\tfrac12(\eee^{-\pi\mathrm iz/2}\mathrm i^{d-1}
+\eee^{\pi\mathrm iz/2}(-\mathrm i)^{d-1})
=\sin(\pi(d-z)/2)$. Consequently,
\begin{equation}
\label{eq:gaussian-Mellin-sine-weighted-series}
\mathcal S_z^\bullet(w)
=\sum_{d=1}^\infty
\sin\left(\frac{\pi(d-z)}2\right)N_d^\bullet(z)w^{d-1}.
\end{equation}

We evaluate these combinations using the beta-function formulas
from the preceding subsection. The reflection and duplication
formulas for the gamma function
\cite[Equations~(5.5.3) and~(5.5.5)]{NIST:DLMF}
rewrite \eqref{eq:type-B-gaussian-Mellin-generating-function} as
\[
\mathcal N_z^B(w)
=\frac1{1-w^2}\bigg[
\frac1{z\cos(\pi z/2)}
+\frac{
 w\Gamma((z+1+w)/2)\Gamma((z+1-w)/2)
 \cos(\pi(z+w)/2)}
 {\Gamma(1+z/2)^2\sin(\pi z)}
\bigg].
\]
Indeed, the constant term follows by simplifying
$2^z\Gamma((z+1)/2)\mathrm B(z/2,-z)/
(\sqrt\pi\,\Gamma(1+z/2))$. For the second term, use
\[
\frac1{\Gamma((1-z-w)/2)}
=\frac{\Gamma((1+z+w)/2)}\pi
\cos\left(\frac{\pi(z+w)}2\right)
\]
and $\Gamma(-z)=-\pi/(\Gamma(1+z)\sin(\pi z))$, followed
by duplication. Substitution of $\mathrm iw$ and $-\mathrm iw$
now gives
\begin{equation}
\label{eq:type-B-transformed-Mellin-generating-function}
\begin{aligned}
\mathcal S_z^B(w)
=\frac1{1+w^2}\bigg[
\frac1z
+\frac{w\eee^{\pi w/2}}{2\Gamma(1+z/2)^2}
\Gamma\left(\frac{z+1+\mathrm iw}{2}\right)
\Gamma\left(\frac{z+1-\mathrm iw}{2}\right)
\bigg].
\end{aligned}
\end{equation}
The trigonometric simplification uses
\[
\eee^{-\mathrm i\alpha}\cos(\alpha+\mathrm ib)
-\eee^{\mathrm i\alpha}\cos(\alpha-\mathrm ib)
=-\mathrm i\eee^b\sin(2\alpha),
\qquad \alpha,b\in\mathbb C,
\]
with $\alpha=\pi z/2$ and $b=\pi w/2$.

For type $A$, put $a:=1+z/2$. Reflection in the two beta
quotients of \eqref{eq:type-A-gaussian-Mellin-generating-function}
gives
\[
\mathcal N_z^A(w)
=\frac{
\Gamma(a+1)\Gamma(a-1)\sin(\pi z/2)
+\Gamma(a+w)\Gamma(a-w)\sin(\pi z/2+\pi w)
}{\Gamma(z+2)\sin(\pi z)(1-w^2)}.
\]
Here one uses
$\Gamma(-z-1)=\pi/(\Gamma(z+2)\sin(\pi z))$.
In the combination defining $\mathcal S_z^A$, the identity
\[
\eee^{-\mathrm i\alpha}\sin(\alpha+\mathrm ib)
+\eee^{\mathrm i\alpha}\sin(\alpha-\mathrm ib)
=\eee^b\sin(2\alpha)
\]
with $\alpha=\pi z/2$ and $b=\pi w$ yields
\begin{equation}
\label{eq:type-A-transformed-Mellin-generating-function}
\mathcal S_z^A(w)
=\frac{
\Gamma(a+1)\Gamma(a-1)
+\eee^{\pi w}\Gamma(a+\mathrm iw)\Gamma(a-\mathrm iw)
}{2\Gamma(z+2)(1+w^2)}.
\end{equation}
All gamma factors in these formulas are retained as products,
so the identities apply to complex $z$ and $w$.

\medskip
\emph{Cosine integrals and coefficient extraction.}
Applying \eqref{eq:hyperbolic-secant-cosine-transform} with
$\nu=z+1$ and $\xi=w$, and using duplication, gives
\[
\frac{\Gamma((z+1)/2)}{\sqrt\pi\,\Gamma(1+z/2)}
\int_0^\infty\frac{\cos(wu)}{\cosh^{z+1}u}\,\dint u
=
\frac{
\Gamma((z+1+\mathrm iw)/2)\Gamma((z+1-\mathrm iw)/2)
}{2\Gamma(1+z/2)^2}.
\]
The same formula with $\xi=\mathrm i$, followed by the gamma
recurrence, gives
\[
\frac{\Gamma((z+1)/2)}{\sqrt\pi\,\Gamma(1+z/2)}
\int_0^\infty\frac{\cosh u}{\cosh^{z+1}u}\,\dint u
=\frac1z.
\]
It follows from
\eqref{eq:type-B-transformed-Mellin-generating-function} that
\begin{equation}
\label{eq:type-B-transformed-Mellin-generating-integral}
\mathcal S_z^B(w)
=\frac{\Gamma((z+1)/2)}{\sqrt\pi\,\Gamma(1+z/2)}
\int_0^\infty
\frac{\cosh u+w\eee^{\pi w/2}\cos(wu)}
{(1+w^2)\cosh^{z+1}u}\,\dint u.
\end{equation}
For type $A$, replace $u$ by $u/2$ in the cosine integral and
use $\nu=z+2$, first with $\xi=2w$ and then with
$\xi=2\mathrm i$. Equation
\eqref{eq:type-A-transformed-Mellin-generating-function} becomes
\begin{equation}
\label{eq:type-A-transformed-Mellin-generating-integral}
\mathcal S_z^A(w)
=2^{-z-2}\int_0^\infty
\frac{\cosh u+\eee^{\pi w}\cos(wu)}
{(1+w^2)\cosh^{z+2}(u/2)}\,\dint u.
\end{equation}
The conditions in
\eqref{eq:hyperbolic-secant-cosine-transform} are satisfied
in all four applications, since $\operatorname{Re}z>0$
and $|w|<1$.

We may extract coefficients under these integrals. To justify
this, fix $\delta>0$ and $0<\rho<1$, restrict $z$ to a compact
subset of $0<\operatorname{Re}z<1$ with
$\operatorname{Re}z\ge\delta$, and take $|w|\le\rho$.
The absolute values of the two integrands, including their
prefactors, are bounded respectively by
\[
C\bigl(\eee^{-\delta u}
+\eee^{-(1+\delta-\rho)u}\bigr),
\qquad
C\bigl(\eee^{-\delta u/2}
+\eee^{-(1+\delta/2-\rho)u}\bigr).
\]
Both bounds are integrable on $[0,\infty)$ and independent of
$z,w$ in these sets. Thus the integrals are holomorphic in $w$,
and Cauchy's coefficient formula permits coefficient extraction
under the integral sign. Using the generating identities
\eqref{eq:intro-type-B-factor-generating-function}
and~\eqref{eq:intro-type-A-factor-generating-function}, and
comparing with \eqref{eq:gaussian-Mellin-sine-weighted-series},
we obtain
\begin{align}
\sin\left(\frac{\pi(d-z)}2\right)N_d^B(z)
&=\frac{\Gamma((z+1)/2)}{\sqrt\pi\,\Gamma(1+z/2)}
\int_0^\infty
\frac{P_d^B(u)}{\cosh^{z+1}u}\,\dint u,
\label{eq:type-B-gaussian-Mellin-factor-integral}
\\
\sin\left(\frac{\pi(d-z)}2\right)N_d^A(z)
&=2^{-z-2}\int_0^\infty
\frac{P_d^A(u)}{\cosh^{z+2}(u/2)}\,\dint u.
\label{eq:type-A-gaussian-Mellin-factor-integral}
\end{align}
So far, these identities hold for every $d\ge1$ on
$0<\operatorname{Re}z<1$.

\medskip
\emph{Continuation in the Mellin parameter.}
For each fixed $d$, the growth bound
\eqref{eq:prescribed-factors-growth} implies that the two
integrals in
\eqref{eq:type-B-gaussian-Mellin-factor-integral}
and~\eqref{eq:type-A-gaussian-Mellin-factor-integral} converge
absolutely and locally uniformly on $\operatorname{Re}z>0$.
The same bound gives integrable majorants on compact subsets
of this half-plane, so both integrals are holomorphic there.
By Lemma~\ref{lem:radial-gaussian-Mellin-transforms}, the
left-hand sides are holomorphic on
$0<\operatorname{Re}z<d$. The identity theorem therefore
extends both identities to that entire strip. We apply this
continuation before dividing by the sine factor.

Now assume $\max\{0,d-2\}<\operatorname{Re}z<d$. Then
$0<\operatorname{Re}((d-z)/2)<1$, and reflection gives
\[
\frac1{\Gamma((d-z)/2)\sin(\pi(d-z)/2)}
=\frac{\Gamma(1+(z-d)/2)}\pi.
\]
Substitute the radial relation
\eqref{eq:radial-gaussian-Mellin-relation} into
\eqref{eq:type-B-gaussian-Mellin-factor-integral}
and~\eqref{eq:type-A-gaussian-Mellin-factor-integral}, and
use this identity. The resulting formulas are exactly
\eqref{eq:type-B-prescribed-Mellin-integral}
and~\eqref{eq:type-A-prescribed-Mellin-integral}.

Finally, the Mellin transforms $M_d^B(z)$ and $M_d^A(z)$
are holomorphic on $\operatorname{Re}z>0$. Indeed,
the defining integrals converge absolutely and locally
uniformly there, by the Gaussian tail bounds for
$p_d^B$ and $p_d^A$.
The integrals on the right-hand sides of
\eqref{eq:type-B-prescribed-Mellin-integral}
and~\eqref{eq:type-A-prescribed-Mellin-integral}
are also holomorphic on this half-plane, as shown above.
Since their gamma prefactors are meromorphic, the
identity theorem extends both formulas throughout
$\operatorname{Re}z>0$, away from the poles of these
prefactors. The apparent singularities at the positive
values among $d-2,d-4,\ldots$ are removable because the
resulting expressions agree with the holomorphic
functions $M_d^B$ and $M_d^A$.
\end{proof}

These Mellin identities were derived from the Gaussian-point
probabilities, without using the prescribed-point integral
formulas. In the next subsection we identify their inverse
transforms and complete the proof of
Theorem~\ref{thm:intro-prescribed-absorption}.

\subsection{Proof of the prescribed-point absorption formulas}
\label{subsec:proof-prescribed-point-absorption}

We now identify the Mellin transforms in
Proposition~\ref{prop:prescribed-point-Mellin-integrals} with
those of the integral expressions in
Theorem~\ref{thm:intro-prescribed-absorption}. The absolute
estimates below justify the comparison without relying on
cancellations in the factors $P_d^B$ and $P_d^A$.

\begin{proof}[Proof of Theorem~\ref{thm:intro-prescribed-absorption}]
Fix an integer $d\ge1$. The substitutions $v=\sinh u$ for
type $B$ and $v=\sinh(u/2)$ for type $A$ show that the asserted
identities are equivalent to
\begin{align}
p_d^B(r)
&=\frac{\sqrt2\,r}{\pi^{3/2}}
\int_0^\infty
P_d^B(u)\,G_d^B\!\left(\frac{r^2\cosh^2u}{2}\right)
\,\dint u,
\label{eq:type-B-prescribed-absorption-hyperbolic}
\\
p_d^A(r)
&=\frac{\Gamma(d/2)}{2^{d/2}\pi}\,r^{2-d}
\int_0^\infty
\frac{P_d^A(u)}{\cosh^d(u/2)}
\eee^{-2r^2\cosh^2(u/2)}\,\dint u,
\qquad r>0.
\label{eq:type-A-prescribed-absorption-hyperbolic}
\end{align}
Denote the right-hand sides by $\widetilde p_d^B(r)$ and
$\widetilde p_d^A(r)$, respectively. We first establish their
convergence and regularity, and then compare their Mellin
transforms with those of the absorption probabilities.

\medskip
\emph{Absolute convergence and smoothness.}
Differentiation under the integral in
\eqref{eq:intro-type-B-analytic-kernel} gives, for every integer
$j\ge0$ and every $x>0$,
\begin{equation}
\label{eq:type-B-analytic-kernel-derivative-integral}
(G_d^B)^{(j)}(x)
=(-1)^j\int_0^\infty
t^{d/2-1}(1+t)^{j-1/2}\eee^{-x(1+t)}\,\dint t.
\end{equation}
These differentiations are locally dominated for $x>0$.
For every $x_0>0$, splitting the exponential into two equal
factors shows that
$|(G_d^B)^{(j)}(x)|\le C_{d,j,x_0}\eee^{-x/2}$ for $x\ge x_0$,
where $C_{d,j,x_0}<\infty$.

Fix $0<a<b<\infty$ and an integer $j\ge0$. This bound and
\eqref{eq:prescribed-factors-growth}, together with the product
and chain rules, control the $j$th derivative with respect to
$r$ of the type-$B$ integrand, including its prefactor $r$.
For type $A$, differentiate its exponential and the prefactor
$r^{2-d}$ directly. In both cases, for $a\le r\le b$, the
absolute value of the resulting integrand is bounded by
\[
C\bigl(1+u^{d+1}+\cosh u\bigr)\cosh^{2j}u\,
\eee^{-(a^2/4)\cosh^2(u/2)},
\qquad u\ge0,
\]
with a constant depending only on $d,j,a,b$. This majorant is
integrable on $[0,\infty)$. Taking $j=0$ proves absolute
convergence, and the bounds for all $j$ justify differentiation
under the integrals. Hence both $\widetilde p_d^B$ and
$\widetilde p_d^A$ belong to $C^\infty((0,\infty))$.

\medskip
\emph{Absolute integrability and Mellin transforms.}
Let $z\in\mathbb C$ satisfy
$\max\{0,d-2\}<\operatorname{Re}z<d$, and put
$c:=\operatorname{Re}z$. We first take absolute values and use
Tonelli's theorem to verify the integrability required for the
Mellin-transform calculations.

For type $B$, insert the defining integral
\eqref{eq:intro-type-B-analytic-kernel} into the expression for
$\widetilde p_d^B$. The gamma integral in $r$, followed by
Euler's beta integral in $t$, gives
\begin{equation}
\label{eq:type-B-candidate-absolute-Mellin-integral}
\begin{aligned}
&\frac{\sqrt2}{\pi^{3/2}}
\int_0^\infty\int_0^\infty\int_0^\infty
r^c|P_d^B(u)|\frac{t^{d/2-1}}{\sqrt{1+t}}\,
\eee^{-r^2(1+t)\cosh^2u/2}
\,\dint r\,\dint t\,\dint u
\\
&\qquad=
\frac{2^{c/2}\Gamma((c+1)/2)}{\pi^{3/2}}
\mathrm B\!\left(\frac d2,1+\frac{c-d}{2}\right)
\int_0^\infty
\frac{|P_d^B(u)|}{\cosh^{c+1}u}\,\dint u
<\infty.
\end{aligned}
\end{equation}
Both beta parameters are positive because $c>d-2$, and the
last integral is finite by
\eqref{eq:prescribed-factors-growth} and $c>0$. In particular,
\eqref{eq:type-B-candidate-absolute-Mellin-integral} bounds
$\int_0^\infty r^{c-1}|\widetilde p_d^B(r)|\,\dint r$.
It also justifies Fubini's theorem with $z$ in place of $c$.
Removing the absolute values and using
\[
\mathrm B\!\left(\frac d2,1+\frac{z-d}{2}\right)
=\frac{\Gamma(d/2)\Gamma(1+(z-d)/2)}{\Gamma(1+z/2)}
\]
therefore gives exactly the right-hand side of
\eqref{eq:type-B-prescribed-Mellin-integral} as the Mellin
transform of $\widetilde p_d^B$.

For type $A$, the gamma integral in $r$ gives
\begin{equation}
\label{eq:type-A-candidate-absolute-Mellin-integral}
\begin{aligned}
&\frac{\Gamma(d/2)}{2^{d/2}\pi}
\int_0^\infty\int_0^\infty
r^{c+1-d}\frac{|P_d^A(u)|}{\cosh^d(u/2)}\,
\eee^{-2r^2\cosh^2(u/2)}\,\dint r\,\dint u
\\
&\qquad=
\frac{\Gamma(d/2)\Gamma(1+(c-d)/2)}{2^{c/2+2}\pi}
\int_0^\infty
\frac{|P_d^A(u)|}{\cosh^{c+2}(u/2)}\,\dint u
<\infty.
\end{aligned}
\end{equation}
Here $c>d-2$ ensures integrability at $r=0$, while
\eqref{eq:prescribed-factors-growth} and $c>0$ give finiteness
of the last integral. Thus
$\int_0^\infty r^{c-1}|\widetilde p_d^A(r)|\,\dint r<\infty$,
and Fubini's theorem again applies. Replacing $c$ by $z$ and
removing the absolute values gives precisely the right-hand
side of \eqref{eq:type-A-prescribed-Mellin-integral}.

Proposition~\ref{prop:prescribed-point-Mellin-integrals}
therefore implies, for both types,
\begin{equation}
\label{eq:prescribed-candidates-Mellin-identification}
\int_0^\infty r^{z-1}\widetilde p_d^\bullet(r)\,\dint r
=M_d^\bullet(z),
\qquad
\max\{0,d-2\}<\operatorname{Re}z<d,
\quad \bullet\in\{B,A\}.
\end{equation}

\medskip
\emph{Identification at every positive radius.}
Fix a real $c$ with $\max\{0,d-2\}<c<d$.
The preceding absolute estimates and
Lemma~\ref{lem:radial-gaussian-Mellin-transforms} show that
\[
x\longmapsto
\eee^{cx}\bigl(\widetilde p_d^\bullet(\eee^x)
-p_d^\bullet(\eee^x)\bigr)
\]
belongs to $L^1(\R)$. Evaluating
\eqref{eq:prescribed-candidates-Mellin-identification}
at $z=c+\mathrm it$ and substituting $r=\eee^x$ yields
\[
\int_{\R}\eee^{\mathrm itx}\eee^{cx}
\bigl(\widetilde p_d^\bullet(\eee^x)
-p_d^\bullet(\eee^x)\bigr)\,\dint x=0,
\qquad t\in\R.
\]
Uniqueness of Fourier transforms on $L^1(\R)$ implies
$\widetilde p_d^\bullet=p_d^\bullet$ almost everywhere on
$(0,\infty)$. This is the Mellin uniqueness step; no contour
inversion is needed.

Finally, $p_d^\bullet$ is nonincreasing because
$\mathcal K_d^\bullet$ is convex and contains the origin,
whereas $\widetilde p_d^\bullet$ is continuous. For any $r>0$,
choose $r_n\uparrow r$ and $s_n\downarrow r$ outside the null
set on which equality has not yet been established. Monotonicity
gives
\[
\widetilde p_d^\bullet(r_n)
\ge p_d^\bullet(r)
\ge\widetilde p_d^\bullet(s_n).
\]
Letting $n\to\infty$ proves equality at $r$. Thus
\eqref{eq:type-B-prescribed-absorption-hyperbolic}
and~\eqref{eq:type-A-prescribed-absorption-hyperbolic} hold
for every $r>0$. The substitutions at the beginning of the
proof give
\eqref{eq:intro-type-B-prescribed-absorption}
and~\eqref{eq:intro-type-A-prescribed-absorption}.
Their absolute convergence and smoothness were established
above, completing the proof.
\end{proof}

\section{Special-function representations and Mellin transforms}
\label{sec:special-functions-and-Mellin-transforms}

We now examine the functions appearing in the absorption formulas of
Theorem~\ref{thm:intro-prescribed-absorption}. We first describe the
structure of the elementary factors $P_d^B,P_d^A$ and reduce the
kernel $G_d^B$ to fixed families of classical special functions in
every integer dimension. We then use the Mellin identities from
Proposition~\ref{prop:prescribed-point-Mellin-integrals} to obtain
finite expressions in gamma and polygamma functions and selected
exact Mellin values.

\subsection{Elementary factors and the type-B kernel}
\label{subsec:elementary-factors-and-type-B-kernel}

Recall that the factors $P_d^B$ and $P_d^A$ are defined by
\eqref{eq:intro-type-B-prescribed-factor}
and~\eqref{eq:intro-type-A-prescribed-factor} in terms of the
hyperbolic Taylor remainders $\Psi_m$. Their elementary structure
and the entries in
Tables~\ref{tab:intro-type-B-prescribed-factors}
and~\ref{tab:intro-type-A-prescribed-factors} follow from the
recurrences below.

\begin{proposition}[Structure of the elementary factors]
\label{prop:prescribed-factors-elementary-structure}
For $u\ge0$, the initial factors are
\begin{equation}
\label{eq:prescribed-factors-initial-values}
\begin{aligned}
P_1^B(u)&=\cosh u, & P_2^B(u)&=1,\\
P_1^A(u)&=1+\cosh u, & P_2^A(u)&=\pi.
\end{aligned}
\end{equation}
For every integer $d\ge1$,
\begin{align}
P_{d+2}^B(u)
&=-P_d^B(u)
+\frac{\operatorname{Re}(\pi/2+\mathrm iu)^d}{d!},
\label{eq:type-B-factor-dimension-recurrence}
\\
P_{d+2}^A(u)
&=-P_d^A(u)
+\frac{\operatorname{Re}(\pi+\mathrm iu)^{d+1}}{(d+1)!}.
\label{eq:type-A-factor-dimension-recurrence}
\end{align}
Both families also satisfy
\begin{equation}
\label{eq:prescribed-factors-differential-recurrence}
(P_{d+2}^\bullet)''(u)=-P_d^\bullet(u),
\qquad \bullet\in\{B,A\}.
\end{equation}
For every $m\ge1$, the functions $P_{2m}^B$ and $P_{2m}^A$
are even polynomials of degree at most $2m-2$. In odd dimensions,
$P_{2m+1}^B(u)-(-1)^m\cosh u$ is an even polynomial of degree
at most $2m-2$ for $m\ge1$, whereas
$P_{2m+1}^A(u)-(-1)^m\cosh u$ is an even polynomial of degree
at most $2m$ for $m\ge0$.

The coefficients of the even-dimensional polynomials belong to
$\Q[\pi^2]$ in type $B$ and to $\pi\Q[\pi^2]$ in type $A$.
For the odd-dimensional polynomial remainders, these two
coefficient classes are interchanged. Here $\Q[\pi^2]$ denotes
the set of polynomials in $\pi^2$ with rational coefficients.
\end{proposition}

\begin{proof}
The initial values follow directly from
$\Psi_1=\sinh$, $\Psi_2=\cosh-1$, and $\Psi_3(z)=\sinh z-z$.
Substitute $z=u+\pi\mathrm i/2$ into
$\Psi_{d+2}(z)=\Psi_d(z)-z^d/d!$, multiply by
$\mathrm i^{-(d+2)}$, and take real parts. Since
$u+\pi\mathrm i/2=\mathrm i(\pi/2-\mathrm iu)$, this gives
\eqref{eq:type-B-factor-dimension-recurrence}. Similarly,
substitute $z=u+\pi\mathrm i$ into
$\Psi_{d+3}(z)=\Psi_{d+1}(z)-z^{d+1}/(d+1)!$ and take the
imaginary part after multiplication by $\mathrm i^{-(d+2)}$.
The identity $u+\pi\mathrm i=\mathrm i(\pi-\mathrm iu)$ gives
\eqref{eq:type-A-factor-dimension-recurrence}. Differentiating
the definitions twice and using $\Psi_{m+2}''=\Psi_m$ proves
\eqref{eq:prescribed-factors-differential-recurrence}.

For real $a,u$ and every integer $n\ge0$,
\[
\operatorname{Re}(a+\mathrm iu)^n
=\sum_{j=0}^{\lfloor n/2\rfloor}
(-1)^j\binom n{2j}a^{n-2j}u^{2j}.
\]
Thus the added terms in the two dimension recurrences are even
polynomials. Their degrees and powers of $\pi$, together with
\eqref{eq:prescribed-factors-initial-values}, give all the
assertions about the polynomial structure by induction. The
coefficient of $\cosh u$ changes sign at each step, and hence
equals $(-1)^m$ in dimension $2m+1$.
\end{proof}

We next consider the type-$B$ kernel. Recall from
\eqref{eq:intro-type-B-analytic-kernel} that
\[
G_d^B(x)=\int_0^\infty
\frac{t^{d/2-1}}{\sqrt{1+t}}\,
\eee^{-x(1+t)}\,\dint t,
\qquad x>0.
\]
For the analytic identities below, the same integral defines
$G_d^B$ for every real $d>0$.
Equivalently, $G_d^B(x)=\Gamma(d/2)\eee^{-x}
U(d/2,(d+1)/2,x)$, where $U$ is Tricomi's confluent
hypergeometric function; see
\cite[Equation~(13.4.4)]{NIST:DLMF}.
For integer dimensions, the following proposition reduces this
representation to Bessel functions or the standard normal tail.

\begin{proposition}[Classical representations of the type-B kernel]
\label{prop:type-B-kernel-classical-representations}
For every real $d>0$ and every $x>0$,
\begin{equation}
\label{eq:type-B-kernel-differential-dimension-recurrence}
G_{d+2}^B(x)=-(G_d^B)'(x)-G_d^B(x).
\end{equation}
The initial kernels are
\begin{align}
G_1^B(x)&=\eee^{-x/2}K_0(x/2),
\label{eq:type-B-kernel-initial-odd}
\\
G_2^B(x)&=2\sqrt{\frac\pi x}\,
\overline\Phi(\sqrt{2x}),
\label{eq:type-B-kernel-initial-even}
\end{align}
where $K_\nu$ is the modified Bessel function of the second kind
and $\overline\Phi$ is the standard normal tail.
For every integer $m\ge0$, there are polynomials
$a_m^B,b_m^B,\widetilde a_m^B,\widetilde b_m^B$ with rational
coefficients such that, for $x>0$,
\begin{align}
G_{2m+1}^B(x)
&=\eee^{-x/2}\bigl[
 a_m^B(1/x)K_0(x/2)+b_m^B(1/x)K_1(x/2)\bigr],
\label{eq:type-B-kernel-odd-classical-representation}
\\
G_{2m+2}^B(x)
&=\widetilde a_m^B(1/x)\eee^{-x}
+2\widetilde b_m^B(1/x)\sqrt{\frac\pi x}\,
\overline\Phi(\sqrt{2x}).
\label{eq:type-B-kernel-even-classical-representation}
\end{align}
\end{proposition}

\begin{proof}
For fixed $d>0$, differentiation under the defining integral
is locally dominated on $x>0$. Expanding the factor $1+t$ in
the differentiated integrand gives
\[
(G_d^B)'(x)
=-\int_0^\infty
\frac{t^{d/2-1}(1+t)}{\sqrt{1+t}}\,
\eee^{-x(1+t)}\,\dint t
=-G_d^B(x)-G_{d+2}^B(x),
\]
which proves
\eqref{eq:type-B-kernel-differential-dimension-recurrence}.

For $d=1$, put $t=\sinh^2(u/2)$ in the defining integral.
The measure $\dint t/\sqrt{t(1+t)}$ becomes $\dint u$, while
$1+t=(1+\cosh u)/2$. Hence
\[
G_1^B(x)=\eee^{-x/2}
\int_0^\infty\eee^{-(x/2)\cosh u}\,\dint u
=\eee^{-x/2}K_0(x/2),
\]
by the Bessel integral
\cite[Equation~(10.32.9)]{NIST:DLMF}.
For $d=2$, the substitution $y=\sqrt{2x(1+t)}$ gives
\[
G_2^B(x)=\sqrt{\frac2x}
\int_{\sqrt{2x}}^\infty\eee^{-y^2/2}\,\dint y
=2\sqrt{\frac\pi x}\,\overline\Phi(\sqrt{2x}).
\]

It remains to check that the two forms in
\eqref{eq:type-B-kernel-odd-classical-representation}
and~\eqref{eq:type-B-kernel-even-classical-representation}
are preserved by the operation $f\mapsto-f'-f$.
For the odd-dimensional form, use
\[
\frac{\dint}{\dint x}K_0(x/2)=-\frac12K_1(x/2),
\qquad
\frac{\dint}{\dint x}K_1(x/2)
=-\frac12K_0(x/2)-\frac1xK_1(x/2),
\]
which follow from the standard Bessel derivative identities;
see \cite[Section~10.29]{NIST:DLMF}.
Differentiating a polynomial in $1/x$ again gives a polynomial
in $1/x$, so the required form and rational coefficients are
preserved. For the even-dimensional form, differentiation of
\eqref{eq:type-B-kernel-initial-even} yields
\[
(G_2^B)'(x)=-\frac1{2x}G_2^B(x)-\frac{\eee^{-x}}x.
\]
Thus the operation $f\mapsto-f'-f$ also preserves sums of
$\eee^{-x}$ and $G_2^B(x)$ multiplied by polynomials in $1/x$
with rational coefficients. The initial kernels have the
claimed forms, and induction using
\eqref{eq:type-B-kernel-differential-dimension-recurrence}
proves both representations.
\end{proof}

Successive applications of
\eqref{eq:type-B-kernel-differential-dimension-recurrence} to
\eqref{eq:type-B-kernel-initial-odd}
and~\eqref{eq:type-B-kernel-initial-even} give the explicit
kernels $G_2^B,\ldots,G_6^B$ listed in
\eqref{eq:intro-type-B-kernels-two-to-six}.
The two propositions above also give a uniform description
of the functions needed in the absorption integrals.

\begin{corollary}[One-dimensional integrals of classical functions]
\label{cor:absorption-classical-integrands}
For every integer $d\ge1$ and every $r>0$, the absorption
probabilities $p_d^B(r)$ and $p_d^A(r)$ admit absolutely
convergent one-dimensional integral representations of the
following form. For type $B$ in odd dimensions, the integrand
is a linear combination of $K_0$ and $K_1$ evaluated at
$r^2(1+v^2)/4$, with elementary coefficients. In even dimensions,
it is a linear combination of $\overline\Phi(r\sqrt{1+v^2})$
and $\eee^{-r^2(1+v^2)/2}$, again with elementary coefficients.
For type $A$, the integrand is elementary in every dimension.
\end{corollary}

\begin{proof}
Use \eqref{eq:intro-type-B-prescribed-absorption}
and~\eqref{eq:intro-type-A-prescribed-absorption}.
Proposition~\ref{prop:prescribed-factors-elementary-structure}
and the identities $\cosh(\arsinh v)=\sqrt{1+v^2}$ and
$\cosh(2\arsinh v)=1+2v^2$ show that
$P_d^B(\arsinh v)$ and $P_d^A(2\arsinh v)$ are elementary
functions of $v$. In type $B$, substitute the corresponding
representation of $G_d^B$ from
Proposition~\ref{prop:type-B-kernel-classical-representations};
in type $A$, the remaining factors are already elementary.
Absolute convergence follows from
Theorem~\ref{thm:intro-prescribed-absorption}.
\end{proof}

The corollary concerns the integrands in one-dimensional
representations; it does not assert that the probabilities
themselves are finite combinations of the indicated functions
without a remaining integration.

\subsection{Finite gamma–polygamma representations of Mellin transforms}
\label{subsec:further-properties-of-Mellin-transforms}

The polynomial structure of $P_d^B,P_d^A$ allows us to evaluate
both integrals in Proposition~\ref{prop:prescribed-point-Mellin-integrals}
as finite expressions in gamma and polygamma functions. We first
note that the defining integral
$M_d^\bullet(z)=\int_0^\infty r^{z-1}p_d^\bullet(r)\,\dint r$
in fact converges absolutely and is holomorphic throughout
$\operatorname{Re}z>0$. Indeed, projection onto the first coordinate
and \eqref{eq:intro-prescribed-absorption-dimension-one} give
\[
0\le p_d^B(r)\le2\overline\Phi(r),
\qquad
0\le p_d^A(r)\le\eee^{-2r^2},
\qquad r\ge0.
\]
These bounds give integrability at infinity for every Mellin
parameter, while $p_d^\bullet\le1$ gives integrability at zero
when $\operatorname{Re}z>0$. The same bounds give locally uniform
domination on this half-plane and hence holomorphy.

Write $\psi_0(a):=\Gamma'(a)/\Gamma(a)$ for the digamma function
and $\psi_k(a):=\psi_0^{(k)}(a)$ for $k\ge1$ for the polygamma
functions; see \cite[Section~5.15]{NIST:DLMF}. To describe the
finite expressions explicitly, put $\varepsilon_d:=\sin(\pi d/2)$
and write
\begin{equation}
\label{eq:prescribed-factors-polynomial-coefficients}
P_d^\bullet(u)=\varepsilon_d\cosh u+
\sum_{j\ge0}b_{d,j}^\bullet u^{2j},
\qquad \bullet\in\{B,A\}.
\end{equation}
By Proposition~\ref{prop:prescribed-factors-elementary-structure},
these sums are finite, and their coefficients are determined by
\eqref{eq:type-B-factor-dimension-recurrence}
and~\eqref{eq:type-A-factor-dimension-recurrence}. In particular,
all $b_{1,j}^B$ are zero. The coefficients belong to $\Q[\pi]$,
with the more precise parity-dependent descriptions given in that
proposition.

For $\operatorname{Re}a>0$, define $c_0(a):=1$ and, recursively,
\begin{equation}
\label{eq:gamma-product-polygamma-coefficient-recurrence}
c_j(a):=\frac1j\sum_{k=1}^j
\frac{\psi_{2k-1}(a)}{(2k-1)!}\,c_{j-k}(a),
\qquad j\ge1.
\end{equation}
For example,
\[
c_1(a)=\psi_1(a),
\qquad
c_2(a)=\frac12\psi_1(a)^2+\frac1{12}\psi_3(a).
\]
For every $j\ge1$, this recurrence expresses $c_j(a)$ as a
polynomial with rational coefficients in the variables
$\psi_1(a),\psi_3(a),\ldots,\psi_{2j-1}(a)$, of total degree
at most $j$.

\begin{proposition}[Finite gamma--polygamma representations]
\label{prop:finite-gamma-polygamma-Mellin-representations}
For every integer $d\ge1$ and every $z\in\mathbb C$ with
$\operatorname{Re}z>0$,
\begin{equation}
\label{eq:type-B-finite-gamma-polygamma-Mellin}
M_d^B(z)
={}\frac{2^{z/2-1}\Gamma(d/2)\Gamma(1+(z-d)/2)}{\pi}
\biggl[
\frac{2\varepsilon_d}{z}
+\frac{\Gamma((z+1)/2)^2}{\Gamma(1+z/2)^2}
\sum_{j\ge0}\frac{(2j)!}{4^j}\,
 b_{d,j}^B c_j\!\left(\frac{z+1}{2}\right)
\biggr],
\end{equation}
and
\begin{equation}
\label{eq:type-A-finite-gamma-polygamma-Mellin}
M_d^A(z)
={}\frac{\Gamma(d/2)\Gamma(1+(z-d)/2)\Gamma(1+z/2)}
{2^{z/2+2}\sqrt\pi\,\Gamma((z+3)/2)}
\biggl[
\varepsilon_d\frac{z+2}{z}
+\sum_{j\ge0}(2j)!\,b_{d,j}^A c_j\!\left(1+\frac z2\right)
\biggr].
\end{equation}
All sums are finite. The possible singularities at the positive
values among $d-2,d-4,\ldots$ are removable, and the expressions
at those values are interpreted by limits.
\end{proposition}

\begin{proof}
For $\operatorname{Re}a>0$ and $w$ sufficiently close to zero,
Taylor expansion of $\log\Gamma$ gives
\begin{equation}
\label{eq:symmetric-gamma-product-polygamma-expansion}
\begin{aligned}
\frac{\Gamma(a+w)\Gamma(a-w)}{\Gamma(a)^2}
&=\exp\left\{
2\sum_{k=1}^\infty\frac{\psi_{2k-1}(a)}{(2k)!}w^{2k}
\right\}
=\sum_{j=0}^\infty c_j(a)w^{2j}.
\end{aligned}
\end{equation}
The odd Taylor terms cancel. Logarithmic differentiation of
the exponential series gives exactly
\eqref{eq:gamma-product-polygamma-coefficient-recurrence},
so its coefficients are the $c_j$ defined above.

For $\operatorname{Re}z>0$, the cosine integral
\eqref{eq:hyperbolic-secant-cosine-transform}, with an imaginary
frequency, gives
\begin{align}
T^B(z,w)
&:=\int_0^\infty
\frac{\cosh(wu)}{\cosh^{z+1}u}\,\dint u
=\frac{2^{z-1}}{\Gamma(z+1)}
\Gamma\!\left(\frac{z+1+w}{2}\right)
\Gamma\!\left(\frac{z+1-w}{2}\right),
\label{eq:type-B-hyperbolic-Mellin-auxiliary-transform}
\\
T^A(z,w)
&:=\int_0^\infty
\frac{\cosh(wu)}{\cosh^{z+2}(u/2)}\,\dint u
=\frac{2^{z+1}}{\Gamma(z+2)}
\Gamma\!\left(1+\frac z2+w\right)
\Gamma\!\left(1+\frac z2-w\right).
\label{eq:type-A-hyperbolic-Mellin-auxiliary-transform}
\end{align}
The first identity holds when $|\operatorname{Re}w|<1+\operatorname{Re}z$,
and the second when $|\operatorname{Re}w|<1+\operatorname{Re}z/2$.
In particular, both include $w=0$ and $w=1$.
On compact subsets of these domains, the integrands and their
fixed-order derivatives in $w$ have integrable exponential
majorants. We may therefore differentiate under the integrals.

Applying \eqref{eq:symmetric-gamma-product-polygamma-expansion}
with $a=(z+1)/2$ and $w/2$ in type $B$, and with $a=1+z/2$
and $w$ in type $A$, yields
\begin{align}
\int_0^\infty\frac{u^{2j}}{\cosh^{z+1}u}\,\dint u
&=\frac{(2j)!}{4^j}\,T^B(z,0)
 c_j\!\left(\frac{z+1}{2}\right),
\label{eq:type-B-sech-even-moments}
\\
\int_0^\infty\frac{u^{2j}}{\cosh^{z+2}(u/2)}\,\dint u
&=(2j)!\,T^A(z,0)c_j\!\left(1+\frac z2\right),
\qquad j\ge0.
\label{eq:type-A-sech-even-moments}
\end{align}
The values needed for the constant and hyperbolic-cosine terms
are, by the gamma recurrence and duplication formula,
\[
\begin{aligned}
T^B(z,0)
&=\frac{\sqrt\pi\,\Gamma((z+1)/2)}{2\Gamma(1+z/2)},
&
T^B(z,1)
&=\frac{\sqrt\pi\,\Gamma(z/2)}{2\Gamma((z+1)/2)},
\\
T^A(z,0)
&=\frac{\sqrt\pi\,\Gamma(1+z/2)}{\Gamma((z+3)/2)},
&
T^A(z,1)
&=\frac{z+2}{z}\,T^A(z,0).
\end{aligned}
\]

First restrict to $\max\{0,d-2\}<\operatorname{Re}z<d$.
Insert the finite decomposition
\eqref{eq:prescribed-factors-polynomial-coefficients} into
\eqref{eq:type-B-prescribed-Mellin-integral}
and~\eqref{eq:type-A-prescribed-Mellin-integral}.
Evaluate each polynomial term using
\eqref{eq:type-B-sech-even-moments}
and~\eqref{eq:type-A-sech-even-moments}, and each $\cosh u$
term using $T^\bullet(z,1)$. The displayed values of
$T^\bullet(z,0)$ and $T^\bullet(z,1)$ then give
\eqref{eq:type-B-finite-gamma-polygamma-Mellin}
and~\eqref{eq:type-A-finite-gamma-polygamma-Mellin}.

The resulting expressions are meromorphic on
$\operatorname{Re}z>0$. Their only possible poles in this
half-plane come from $\Gamma(1+(z-d)/2)$ and occur at the
positive values among $d-2,d-4,\ldots$. Since $M_d^\bullet$
is holomorphic throughout this half-plane, as shown at the
start of the subsection, the identity theorem extends the
two formulas to the whole half-plane and shows that all these
possible singularities are removable.
\end{proof}

\begin{remark}[Alternative proof]
Proposition~\ref{prop:finite-gamma-polygamma-Mellin-representations}
can also be proved by coefficient extraction from the transformed
Mellin generating functions
\eqref{eq:type-B-transformed-Mellin-generating-function}
and~\eqref{eq:type-A-transformed-Mellin-generating-function}.
Indeed, expanding the symmetric gamma products using
\eqref{eq:symmetric-gamma-product-polygamma-expansion}
and comparing coefficients with the generating functions
\eqref{eq:intro-type-B-factor-generating-function}
and~\eqref{eq:intro-type-A-factor-generating-function}
yields the finite sums involving $b_{d,j}^\bullet c_j$.
The radial Mellin relation
\eqref{eq:radial-gaussian-Mellin-relation}
and the gamma reflection formula then give
\eqref{eq:type-B-finite-gamma-polygamma-Mellin}
and~\eqref{eq:type-A-finite-gamma-polygamma-Mellin}
for $0<\operatorname{Re}z<1$.
Analytic continuation completes the argument.
\end{remark}

\begin{remark}[Mellin transforms of first absorption times]
\label{rem:first-absorption-time-gamma-polygamma}
Recall from
\eqref{eq:intro-type-B-first-absorption-radial-Mellin}
that
\[
\Ee\bigl[(\tau_d(r))^s\bigr]
=-2s\,r^{2s}M_d^B(-2s),
\qquad \operatorname{Re}s<\frac12,
\]
where $r>0$ and the value at $s=0$ is understood by
continuity.
Although Proposition~\ref{prop:finite-gamma-polygamma-Mellin-representations}
is stated for $\operatorname{Re}z>0$, the right-hand side of
\eqref{eq:type-B-finite-gamma-polygamma-Mellin}
extends meromorphically to $\operatorname{Re}z>-1$
and agrees there with the continuation of $M_d^B(z)$
defined in \eqref{eq:intro-type-B-radial-Mellin-continuation}.
Consequently, for every $d\ge1$ and $r>0$, the Mellin transform
$s\mapsto\Ee[(\tau_d(r))^s]$ admits a finite representation
in gamma and polygamma functions throughout
$\operatorname{Re}s<1/2$.
Away from removable singularities, only polygamma functions
of odd order, evaluated at $1/2-s$, are required.
At the removable singularities, the values are obtained
by taking limits.
\end{remark}

The explicit formulas of
Proposition~\ref{prop:finite-gamma-polygamma-Mellin-representations}
show that, away from the removable singularities, the only
polygamma functions required are of odd order, evaluated at
$(z+1)/2$ for type $B$ and at $1+z/2$ for type $A$.
The dependence of these expressions on the indicated polygamma
functions is polynomial, with no polygamma functions in denominators.
For type $A$, a single gamma prefactor multiplies such a
polynomial together with an elementary rational term. The same
is true in even dimensions for type $B$; in odd dimensions,
the contribution of $\varepsilon_d\cosh u$ gives the additional
gamma term visible in
\eqref{eq:type-B-finite-gamma-polygamma-Mellin}.
At removable singularities, taking limits may introduce the
digamma function or even-order polygamma functions. The next
two evaluations illustrate this point.

\begin{corollary}[Two exact Mellin values in dimension three]
\label{cor:dimension-three-exact-Mellin-values}
The first Mellin values of the three-dimensional absorption
probabilities are
\begin{align}
M_3^B(1)=\int_0^\infty p_3^B(r)\,\dint r
&=\sqrt{\frac2\pi}\,(2\log2-1),
\label{eq:type-B-dimension-three-first-Mellin-value}
\\
M_3^A(1)=\int_0^\infty p_3^A(r)\,\dint r
&=\frac{\sqrt\pi}{8\sqrt2}\bigl(7\zeta(3)-6\bigr),
\label{eq:type-A-dimension-three-first-Mellin-value}
\end{align}
where $\zeta$ denotes the Riemann zeta function.
\end{corollary}

\begin{proof}
In dimension three, the factors are
$P_3^B(u)=\pi/2-\cosh u$ and
$P_3^A(u)=\pi^2/2-1-\cosh u-u^2/2$.
Substitution into
\eqref{eq:type-B-finite-gamma-polygamma-Mellin}
and~\eqref{eq:type-A-finite-gamma-polygamma-Mellin} gives
\begin{align}
M_3^B(z)
&=\frac{2^{z/2-2}\Gamma((z-1)/2)}{\sqrt\pi}
\left[
\frac\pi2\frac{\Gamma((z+1)/2)^2}{\Gamma(1+z/2)^2}
-\frac2z
\right],
\label{eq:type-B-dimension-three-Mellin-explicit}
\\
M_3^A(z)
&=\frac{2^{-z/2-1}\Gamma(1+z/2)}{z^2-1}
\left[
\frac{\pi^2}{2}-\psi_1\!\left(1+\frac z2\right)
-\frac{2(z+1)}z
\right],
\label{eq:type-A-dimension-three-Mellin-explicit}
\end{align}
for $\operatorname{Re}z>0$, with the removable singularities
at $z=1$ interpreted by limits.

For type $B$, the square bracket in
\eqref{eq:type-B-dimension-three-Mellin-explicit} vanishes at
$z=1$. The standard values of the digamma function and its
recurrence give
$\psi_0(1)-\psi_0(3/2)=2\log2-2$; see
\cite[Equations~(5.4.12), (5.4.13), and~(5.5.2)]{NIST:DLMF}.
Differentiating the bracket at $z=1$ therefore gives
$2(2\log2-1)$. Since
$\Gamma((z-1)/2)\sim2/(z-1)$ as $z\to1$, taking the limit
in \eqref{eq:type-B-dimension-three-Mellin-explicit} proves
\eqref{eq:type-B-dimension-three-first-Mellin-value}.

For type $A$, the identities
\[
\psi_1(3/2)=\frac{\pi^2}{2}-4,
\qquad
\psi_2(3/2)=16-14\zeta(3)
\]
follow from the half-integer values and recurrence
\cite[Equations~(5.15.3) and~(5.15.5)]{NIST:DLMF}.
Hence the square bracket in
\eqref{eq:type-A-dimension-three-Mellin-explicit} also vanishes
at $z=1$, and its derivative there is
$-\psi_2(3/2)/2+2=7\zeta(3)-6$.
Using $z^2-1\sim2(z-1)$ and $\Gamma(3/2)=\sqrt\pi/2$
in \eqref{eq:type-A-dimension-three-Mellin-explicit} gives
\eqref{eq:type-A-dimension-three-first-Mellin-value}.
\end{proof}

\subsection{Dimension-shift recurrences for Mellin transforms}
\label{subsec:Mellin-dimension-shift-recurrences-and-higher-values}

The differential recurrence for the factors $P_d^B,P_d^A$ yields
recurrences for the Mellin transforms in both the dimension and
the Mellin parameter. We first prove
Proposition~\ref{prop:intro-Mellin-dimension-shift-recurrences},
including the values $M_d^\bullet(d)$. We then give polynomial
representations of the ratios
$M_d^\bullet(d+2k)/M_d^\bullet(d)$ and prove
Corollary~\ref{cor:intro-higher-Mellin-values}.
Throughout this subsection, the Mellin transforms are holomorphic
on $\operatorname{Re}z>0$, as established above.

\begin{proof}[Proof of Proposition~\ref{prop:intro-Mellin-dimension-shift-recurrences}]
For $\operatorname{Re}z>0$, write
\[
J_d^B(z):=\int_0^\infty
\frac{P_d^B(u)}{\cosh^{z+1}u}\,\dint u,
\qquad
J_d^A(z):=\int_0^\infty
\frac{P_d^A(u)}{\cosh^{z+2}(u/2)}\,\dint u.
\]
These integrals converge absolutely and are holomorphic on this
half-plane. Put $w_z^B(u):=\cosh^{-(z+1)}u$ and
$w_z^A(u):=\cosh^{-(z+2)}(u/2)$. Differentiation gives
\begin{align*}
(w_z^B)''
&=(z+1)^2w_z^B-(z+1)(z+2)w_{z+2}^B,\\
(w_z^A)''
&=\frac{(z+2)^2}{4}w_z^A
-\frac{(z+2)(z+3)}4w_{z+2}^A.
\end{align*}
By \eqref{eq:prescribed-factors-differential-recurrence},
$(P_{d+2}^\bullet)''=-P_d^\bullet$. Integrating by parts twice
therefore yields
\[
\int_0^\infty P_{d+2}^\bullet(u)(w_z^\bullet)''(u)\,\dint u
=-J_d^\bullet(z).
\]
There are no boundary terms at zero, because the factors and
weights are even. At infinity, the factors and their derivatives
are sums of polynomials and hyperbolic functions, while the
weights and their derivatives decay faster than $\eee^{-u}$
when $\operatorname{Re}z>0$. The boundary terms therefore
vanish there as well. Consequently,
\begin{align}
J_{d+2}^B(z+2)
&=\frac{z+1}{z+2}J_{d+2}^B(z)
+\frac{J_d^B(z)}{(z+1)(z+2)},
\label{eq:type-B-hyperbolic-Mellin-dimension-shift}
\\
J_{d+2}^A(z+2)
&=\frac{z+2}{z+3}J_{d+2}^A(z)
+\frac{4J_d^A(z)}{(z+2)(z+3)}.
\label{eq:type-A-hyperbolic-Mellin-dimension-shift}
\end{align}

First take $\operatorname{Re}z>d$, so that all gamma factors in
the relevant Mellin representations are finite. Substitution of
\eqref{eq:type-B-prescribed-Mellin-integral}
and~\eqref{eq:type-A-prescribed-Mellin-integral} into these two
identities, followed by $\Gamma(a+1)=a\Gamma(a)$, gives
\eqref{eq:intro-type-B-Mellin-dimension-shift}
and~\eqref{eq:intro-type-A-Mellin-dimension-shift}.
Both sides of the resulting recurrences are holomorphic on
$\operatorname{Re}z>0$. The identity theorem extends them to
this entire half-plane, including the parameters at which the
individual gamma prefactors have poles.

Setting $z=d$ in the recurrences gives
\begin{equation}
\label{eq:Mellin-diagonal-dimension-recurrences}
M_{d+2}^B(d+2)=\frac{d}{(d+2)^2}M_d^B(d),
\qquad
M_{d+2}^A(d+2)=\frac{d}{(d+2)(d+3)}M_d^A(d).
\end{equation}
The initial factors in \eqref{eq:prescribed-factors-initial-values}
and the hyperbolic-secant integral
\eqref{eq:hyperbolic-secant-cosine-transform} at zero frequency
give
\[
J_1^B(1)=\frac\pi2,\quad J_2^B(2)=\frac\pi4,
\qquad
J_1^A(1)=2\pi,\quad J_2^A(2)=\frac{4\pi}{3}.
\]
Using the Mellin prefactors, we obtain
\[
M_1^B(1)=\sqrt{\frac2\pi},\quad M_2^B(2)=\frac14,
\qquad
M_1^A(1)=\frac{\sqrt\pi}{2\sqrt2},\quad M_2^A(2)=\frac16.
\]
Induction in each parity using
\eqref{eq:Mellin-diagonal-dimension-recurrences} now proves
\eqref{eq:intro-Mellin-values-at-dimension}. This derivation
uses neither the known expected-volume formulas nor the
large-scale Gaussian-point asymptotics.
\end{proof}

For the higher values, put
\begin{equation}
\label{eq:higher-Mellin-ratios-definition}
R_{d,k}^\bullet:=\frac{M_d^\bullet(d+2k)}{M_d^\bullet(d)},
\qquad d\ge1,\quad k\ge0,\quad \bullet\in\{B,A\}.
\end{equation}
In particular, $R_{d,0}^\bullet=1$. The following proposition
gives a constructive version of the rationality assertion in
Corollary~\ref{cor:intro-higher-Mellin-values}. Empty products
are understood to equal one.

\begin{proposition}[Polynomial representations of higher Mellin values]
\label{prop:higher-Mellin-values-polynomial-representations}
There are unique sequences of polynomials
$Q_k^B,Q_k^A\in\Q[x]$, starting with $Q_0^B=Q_0^A=1$, such that,
for every integer $k\ge1$,
\begin{align}
(x+2k)Q_k^B(x+2)-xQ_k^B(x)
&=2k(x+2k+1)^2Q_{k-1}^B(x+2),
\label{eq:type-B-higher-Mellin-polynomial-recurrence}
\\
(x+2k)Q_k^A(x+2)-xQ_k^A(x)
&=\frac{k}{2}(x+2k+2)^2Q_{k-1}^A(x+2).
\label{eq:type-A-higher-Mellin-polynomial-recurrence}
\end{align}
Each $Q_k^\bullet$ has degree $2k$. For all integers $d\ge1$
and $k\ge1$,
\begin{align}
R_{d,k}^B
&=\frac{d\,Q_k^B(d)}
{(d+2k)^2\prod_{j=1}^{k-1}(d+2j)},
\label{eq:type-B-higher-Mellin-rational-representation}
\\
R_{d,k}^A
&=\frac{d\,Q_k^A(d)}
{(d+2k)\prod_{j=1}^{k}(d+2j+1)}.
\label{eq:type-A-higher-Mellin-rational-representation}
\end{align}
The leading coefficients of $Q_k^B$ and $Q_k^A$ are $3^{-k}$
and $12^{-k}$, respectively.
\end{proposition}

\begin{proof}
\emph{Rational dependence on the dimension.}
We first establish that each ratio is given by a single rational
function of $d$, rather than by separate expressions for even
and odd dimensions. For this step we use the exact Gaussian-point
formulas \eqref{eq:type-B-gaussian-point-absorption}
and~\eqref{eq:type-A-gaussian-point-absorption}.

For fixed $d\ge1$, all the moments $M_d^\bullet(d+2j)$ are finite.
Expanding the exponential in
\eqref{eq:gaussian-radial-density-integral} gives, for every fixed
$k\ge0$,
\begin{equation}
\label{eq:Gaussian-scale-expansion-higher-Mellin-ratios}
\frac{2^{d/2-1}\Gamma(d/2)}{t^{d/2}M_d^\bullet(d)}
q_d^\bullet(t^{-1/2})
=
\sum_{j=0}^k\frac{(-t/2)^j}{j!}R_{d,j}^\bullet
+O(t^{k+1}),
\qquad t\downarrow0.
\end{equation}
Indeed, for $a\ge0$, the remainder after order $k$ in the Taylor
expansion of $\eee^{-a}$ has absolute value at most
$a^{k+1}/(k+1)!$. The error after integration is therefore
bounded by a constant times $t^{k+1}M_d^\bullet(d+2k+2)$.
No interchange of an infinite series and an improper integral
is needed here.

Define the analytic function near zero
\[
h(t):=\frac{\arsinh\sqrt t}{\sqrt t},
\qquad h(0):=1.
\]
Its Taylor coefficients are rational. We write $[t^j]f(t)$ for
the coefficient of $t^j$, and
$(a)_n:=a(a+1)\cdots(a+n-1)$, with $(a)_0:=1$.
For a variable exponent $x$, the coefficients of $h(t)^x$ are
polynomials in $x$ with rational coefficients: this follows by
writing $h(t)^x=\exp(x\log h(t))$, with the logarithm vanishing
at zero. The same holds after replacing $t$ by $t/4$.

In the type-$B$ Gaussian-point formula, the substitution
$y=t^{-1/2}\sinh u$ gives
\[
q_d^B(t^{-1/2})
=\frac2\pi\int_0^1
\frac{\Psi_d(\arsinh(\sqrt t\,y))}
{\sqrt{1-y^2}\sqrt{1+ty^2}}\,\dint y.
\]
Insert the series for $\Psi_d$, factor out $t^{d/2}$, and expand
in $t$. The expansions are locally uniform for $y\in[0,1]$,
with an integrable majorant proportional to $(1-y^2)^{-1/2}$.
The ratio of the beta integrals arising at orders $k$ and zero is
\[
\frac{\displaystyle\int_0^1
 y^{d+2k}(1-y^2)^{-1/2}\,\dint y}
{\displaystyle\int_0^1
 y^d(1-y^2)^{-1/2}\,\dint y}
=\frac{((d+1)/2)_k}{(1+d/2)_k};
\]
see \cite[Equation~(5.12.1)]{NIST:DLMF}.
Comparison with
\eqref{eq:Gaussian-scale-expansion-higher-Mellin-ratios}, using
the value $M_d^B(d)$ already proved above for the normalization,
yields
\begin{equation}
\label{eq:type-B-higher-Mellin-coefficient-formula}
R_{d,k}^B
=(-2)^k k!\,\frac{((d+1)/2)_k}{(1+d/2)_k}
\sum_{j=0}^k
\frac{[t^{k-j}]\bigl(h(t)^{d+2j}(1+t)^{-1/2}\bigr)}
{(d+1)_{2j}}.
\end{equation}
For type $A$, at $\sigma=t^{-1/2}$ we have
$\theta_\sigma=\sqrt t\,h(t/4)$ and
$\sinh\theta_\sigma=\sqrt t\sqrt{1+t/4}$. Expanding the exact
formula for $q_d^A$ and comparing coefficients in the same way gives
\begin{equation}
\label{eq:type-A-higher-Mellin-coefficient-formula}
R_{d,k}^A
=(-2)^k k!
\sum_{j=0}^k
\frac{[t^{k-j}]\bigl(h(t/4)^{d+1+2j}(1+t/4)^{-1/2}\bigr)}
{(d+2)_{2j}}.
\end{equation}
For every fixed $k$, both right-hand sides are rational functions
of $d$ with rational coefficients. We use $R_{x,k}^\bullet$ for
these rational extensions; no probabilistic interpretation at
noninteger dimensions is needed.

These rational functions are regular at $x=0$ and at every
positive real $x$. Moreover,
\begin{equation}
\label{eq:higher-Mellin-rational-zero-dimension}
R_{0,k}^B=R_{0,k}^A=0,\qquad k\ge1.
\end{equation}
To check this, set $d=0$ in the coefficient expressions. The
sums then extract the coefficient of $t^k$ from, respectively,
\[
\frac{\cosh(\arsinh\sqrt t)}{\sqrt{1+t}}=1,
\qquad
\frac{\sinh(2\arsinh(\sqrt t/2))}
{\sqrt t\sqrt{1+t/4}}=1.
\]
Terms with $j>k$ do not affect these coefficients. This proves
\eqref{eq:higher-Mellin-rational-zero-dimension}.

\medskip
\emph{Polynomial recurrences and identification.}
For a fixed $k\ge1$, the operator
\[
Q(x)\longmapsto(x+2k)Q(x+2)-xQ(x)
\]
preserves the space of polynomials of degree at most $2k$.
On the monomial $x^n$ its coefficient of $x^n$ is $2(n+k)$,
and all other terms have smaller degree. Hence its matrix in
the monomial basis is triangular with nonzero diagonal entries.
Starting from $Q_0^B=Q_0^A=1$, this proves existence and
uniqueness of polynomial solutions of
\eqref{eq:type-B-higher-Mellin-polynomial-recurrence}
and~\eqref{eq:type-A-higher-Mellin-polynomial-recurrence},
with rational coefficients and degree at most $2k$.
There cannot be another polynomial solution of higher degree,
because the operator preserves the degree of every nonzero
polynomial.

It remains to identify these polynomial solutions with the
actual Mellin ratios. Substitute $z=d+2k$ into
\eqref{eq:intro-type-B-Mellin-dimension-shift}
and~\eqref{eq:intro-type-A-Mellin-dimension-shift}, and normalize
using \eqref{eq:Mellin-diagonal-dimension-recurrences}. For
$k\ge1$, this gives
\begin{align}
R_{d+2,k}^B
&=\frac{2k(d+2k+1)^2}{(d+2k+2)^2}R_{d+2,k-1}^B
+\frac{(d+2)^2}{(d+2k+2)^2}R_{d,k}^B,
\label{eq:type-B-normalized-higher-Mellin-recurrence}
\\
R_{d+2,k}^A
&=\frac{k(d+2k+2)}{2(d+2k+3)}R_{d+2,k-1}^A
+\frac{(d+2)(d+3)}{(d+2k+2)(d+2k+3)}R_{d,k}^A.
\label{eq:type-A-normalized-higher-Mellin-recurrence}
\end{align}
Since the ratios have rational extensions, these equalities
hold as rational identities after replacing $d$ by $x$.

For $k\ge1$, define $\widehat Q_k^\bullet(x)$ by solving
\eqref{eq:type-B-higher-Mellin-rational-representation}
and~\eqref{eq:type-A-higher-Mellin-rational-representation}
for the numerator polynomial, with $d$ replaced by $x$ and
$R_{d,k}^\bullet$ by $R_{x,k}^\bullet$. Initially these are
rational functions. They are regular at $x=0$ by
\eqref{eq:higher-Mellin-rational-zero-dimension}, and also at
all positive real $x$. Set $\widehat Q_0^\bullet:=1$.
Substitution into the normalized recurrences shows that
$\widehat Q_k^\bullet$ satisfies the same recurrence as
$Q_k^\bullet$.

We prove their equality by induction on $k$. Once equality is
known at $k-1$, the difference
$H(x):=\widehat Q_k^\bullet(x)-Q_k^\bullet(x)$ is a rational
function, regular at all nonnegative real $x$, satisfying
\[
(x+2k)H(x+2)=xH(x).
\]
At $x=0$ this gives $H(2)=0$, and induction gives
$H(2n)=0$ for every integer $n\ge1$. A rational function with
infinitely many zeros is identically zero. Thus
$\widehat Q_k^\bullet=Q_k^\bullet$, proving the two rational
representations simultaneously for even and odd dimensions.

Finally, let $c_k^\bullet$ be the coefficient of $x^{2k}$ in
$Q_k^\bullet$. The leading coefficient of the left-hand side
of either polynomial recurrence is $6k c_k^\bullet$.
The leading coefficients on the right are $2k c_{k-1}^B$
and $(k/2)c_{k-1}^A$, respectively. Since $c_0^B=c_0^A=1$,
\[
c_k^B=\frac13c_{k-1}^B=3^{-k},
\qquad
c_k^A=\frac1{12}c_{k-1}^A=12^{-k}.
\]
In particular, both polynomials have degree exactly $2k$.
\end{proof}

The polynomial recurrences provide an explicit procedure for
computing the ratios for every fixed $k$, using only rational
arithmetic. The first polynomials are
\begin{equation}
\label{eq:higher-Mellin-first-numerator-polynomials}
Q_1^B(x)=\frac{x^2+6x+11}{3},
\qquad
Q_1^A(x)=\frac{x^2+9x+26}{12}.
\end{equation}
Their substitution into the rational representations gives
\eqref{eq:intro-type-B-Mellin-ratio-two-above-dimension}
and~\eqref{eq:intro-type-A-Mellin-ratio-two-above-dimension}.
The next pair is
\begin{align}
R_{d,2}^B
&=\frac{d(5d^4+92d^3+646d^2+2068d+2589)}
{45(d+2)(d+4)^2},
\label{eq:type-B-higher-Mellin-ratio-four-above-dimension}
\\
R_{d,2}^A
&=\frac{d(5d^4+122d^3+1171d^2+5278d+9624)}
{720(d+3)(d+4)(d+5)}.
\label{eq:type-A-higher-Mellin-ratio-four-above-dimension}
\end{align}
These expressions follow by applying the polynomial recurrences
at $k=2$ and substituting the resulting polynomials into
\eqref{eq:type-B-higher-Mellin-rational-representation}
and~\eqref{eq:type-A-higher-Mellin-rational-representation}.

\begin{proof}[Proof of Corollary~\ref{cor:intro-higher-Mellin-values}]
The rationality assertion follows from
Proposition~\ref{prop:higher-Mellin-values-polynomial-representations},
and \eqref{eq:higher-Mellin-first-numerator-polynomials} gives
the two stated ratios at $k=1$. For fixed $k\ge1$, the
numerators of the rational representations have degree $2k+1$
and leading coefficients $3^{-k}$ and $12^{-k}$, while both
denominators are monic of degree $k+1$. Hence, as $d\to\infty$,
\[
R_{d,k}^B=3^{-k}d^k+O(d^{k-1}),
\qquad
R_{d,k}^A=12^{-k}d^k+O(d^{k-1}).
\]
The constants in these error terms may depend on $k$.
This proves \eqref{eq:intro-higher-Mellin-values-high-dimension}.
For $k=0$, both ratios equal one identically.
\end{proof}

By \eqref{eq:intro-expected-radial-moments-Mellin-relation},
these evaluations determine the expected integrals of every
even power of the distance from the origin over the Brownian
hulls. In particular, the values at $k=0$ recover the expected
volumes, while those at $k=1$ give the expected squared-radius
integrals in
\eqref{eq:intro-type-B-expected-squared-radius-integral}
and~\eqref{eq:intro-type-A-expected-squared-radius-integral}.

\section{Small-radius asymptotics}
\label{sec:small-radius-asymptotics}

In this section we prove
Theorem~\ref{thm:intro-brownian-long-time-nonabsorption} in its
equivalent small-radius form. We first show how monotonicity
allows us to recover a regularly varying asymptotic from a
Gaussian radial average. We then derive the small-scale
Gaussian-point non-absorption asymptotics directly from the
finite sums established in~\cite{KabluchkoOUSpirals} and apply this general result
to both Brownian hulls. The argument does not use the
prescribed-point integral representations.

\subsection{Gaussian radial de-averaging}
\label{subsec:gaussian-radial-deaveraging}

The following lemma identifies the asymptotic factor lost under
Gaussian radial averaging. We give a direct proof using Gaussian
scaling and monotonicity.

\begin{lemma}[Removing Gaussian radial averaging]
\label{lem:radial-gaussian-regular-variation}
Let $d\ge1$ be an integer, let $R_d$ have the $\chi_d$
distribution, and let $f:[0,\infty)\to[0,\infty)$ be bounded,
measurable, and nondecreasing. Put
$\overline f(\sigma):=\Ee[f(\sigma R_d)]$ for $\sigma>0$.
Suppose that $\overline f$ is positive for all sufficiently small
arguments and is regularly varying at zero with index $\beta>0$;
that is,
\begin{equation}
\label{eq:radial-gaussian-regular-variation-assumption}
\frac{\overline f(a\sigma)}{\overline f(\sigma)}
\longrightarrow a^\beta
\qquad\text{as }\sigma\downarrow0,
\qquad a>0.
\end{equation}
Then
\begin{equation}
\label{eq:radial-gaussian-deaveraging-limit}
\frac{f(\sigma)}{\overline f(\sigma)}
\longrightarrow
\frac{\Gamma(d/2)}{2^{\beta/2}\Gamma((d+\beta)/2)}
=\frac1{\Ee[R_d^\beta]}
\qquad\text{as }\sigma\downarrow0.
\end{equation}
\end{lemma}

\begin{proof}
For real $\eta>0$, let $\mu_\eta$ be the probability measure
on $[0,\infty)$ with density
$2^{1-\eta/2}t^{\eta-1}\eee^{-t^2/2}/\Gamma(\eta/2)$ for
$t>0$. Thus $\mu_d$ is the law of $R_d$. For sufficiently
small $\sigma>0$, define the probability measure
\begin{equation}
\label{eq:radial-gaussian-reweighted-measure}
\nu_\sigma(\dint t)
:=\frac{f(\sigma t)}{\overline f(\sigma)}\,\mu_d(\dint t).
\end{equation}
For $\lambda\ge0$, the substitution
$y=\sqrt{1+\lambda}\,t$ in the Gaussian density gives
\begin{equation}
\label{eq:radial-gaussian-reweighted-transform}
\int_0^\infty\eee^{-\lambda t^2/2}\,\nu_\sigma(\dint t)
=(1+\lambda)^{-d/2}
\frac{\overline f(\sigma/\sqrt{1+\lambda})}
{\overline f(\sigma)}
\longrightarrow(1+\lambda)^{-(d+\beta)/2}.
\end{equation}
The limit is the same transform under $\mu_{d+\beta}$.
It follows that
\begin{equation}
\label{eq:radial-gaussian-reweighted-weak-limit}
\nu_\sigma\Longrightarrow\mu_{d+\beta}
\qquad\text{as }\sigma\downarrow0.
\end{equation}
Indeed, for $R,\lambda>0$,
\[
\nu_\sigma((R,\infty))
\le
\frac{1-\int_0^\infty\eee^{-\lambda t^2/2}\,
\nu_\sigma(\dint t)}
{1-\eee^{-\lambda R^2/2}}.
\]
Taking the upper limit as $\sigma\downarrow0$, then choosing
$\lambda$ small and $R$ large, proves tightness on $[0,\infty)$.
Every subsequential weak limit has the transforms in
\eqref{eq:radial-gaussian-reweighted-transform}, because their
integrands are bounded and continuous. Uniqueness of Laplace
transforms, applied after the change of variables $y=t^2/2$,
identifies that limit with $\mu_{d+\beta}$ and proves
\eqref{eq:radial-gaussian-reweighted-weak-limit}.

For $0<\varepsilon<1$, monotonicity of $f$ gives
\begin{equation}
\label{eq:radial-gaussian-monotone-squeeze}
\frac{\nu_\sigma([1-\varepsilon,1])}
{\mu_d([1-\varepsilon,1])}
\le\frac{f(\sigma)}{\overline f(\sigma)}
\le
\frac{\nu_\sigma([1,1+\varepsilon])}
{\mu_d([1,1+\varepsilon])}.
\end{equation}
The intervals are continuity sets of $\mu_{d+\beta}$.
First let $\sigma\downarrow0$ and use
\eqref{eq:radial-gaussian-reweighted-weak-limit}. Then let
$\varepsilon\downarrow0$. Both bounds tend to the ratio of the
$\mu_{d+\beta}$ and $\mu_d$ densities at $1$, namely
$2^{-\beta/2}\Gamma(d/2)/\Gamma((d+\beta)/2)$.
Finally, integration against the density of $R_d$ gives
$\Ee[R_d^\beta]=2^{\beta/2}\Gamma((d+\beta)/2)/\Gamma(d/2)$.
This proves \eqref{eq:radial-gaussian-deaveraging-limit}.
\end{proof}

\subsection{Proof of the large-time non-absorption theorem}
\label{subsec:proof-large-time-nonabsorption}

Fix an integer $d\ge1$. Recall that $q_d^\bullet(\sigma)$ is the
probability that the unit-time Brownian hull $\mathcal K_d^\bullet$ absorbs an
independent centered Gaussian point with covariance matrix
$\sigma^2 I_d$.
To determine $1-q_d^\bullet(\sigma)$ as $\sigma\downarrow0$, we
use the finite complementary formula from the companion paper~\cite[Corollary~5.2]{KabluchkoOUSpirals}:
\begin{equation}
\label{eq:brownian-gaussian-nonabsorption-finite-sum}
1-q_d^\bullet(\sigma)
=2\sum_{j=0}^{\lfloor(d-1)/2\rfloor}
L_{d-1-2j}^\bullet(\sigma),
\qquad \sigma>0,\quad \bullet\in\{B,A\}.
\end{equation}
The coefficient functions occurring here have the explicit
expressions from \cite[Theorem~4.1]{KabluchkoOUSpirals}:
\begin{equation}
\label{eq:small-radius-zeroth-coefficients}
L_0^B(\sigma)=\frac{\arctan\sigma}{\pi},
\qquad
L_0^A(\sigma)=\frac12\left(1-\frac1{\sqrt{1+4\sigma^2}}\right),
\end{equation}
and, for $k\ge1$,
\begin{align}
L_k^B(\sigma)
&=\frac{\sigma}{\pi(k-1)!}
\int_0^{\arsinh(1/\sigma)}
\frac{u^{k-1}}{\sqrt{1-\sigma^2\sinh^2u}}\,\dint u,
\label{eq:type-B-small-radius-positive-coefficients}
\\
L_k^A(\sigma)
&=\frac{\theta_\sigma^k}{2k!\sinh\theta_\sigma},
\label{eq:type-A-small-radius-positive-coefficients}
\end{align}
where $\theta_\sigma:=2\arsinh(1/(2\sigma))$, as in
\eqref{eq:type-A-gaussian-point-absorption}.
Only these explicit expressions are needed below.

\begin{theorem}[Small-radius non-absorption]
\label{thm:brownian-small-radius-nonabsorption}
For every fixed integer $d\ge1$, as $r\downarrow0$,
\begin{align}
1-p_d^B(r)
&\sim
\frac{\sqrt2\,\Gamma(d/2)}
{\pi\Gamma((d+1)/2)(d-1)!}\,
 r\bigl(\log(1/r)\bigr)^{d-1},
\label{eq:type-B-small-radius-nonabsorption}
\\
1-p_d^A(r)
&\sim\frac{2^d}{d!}\,
 r^2\bigl(\log(1/r)\bigr)^{d-1}.
\label{eq:type-A-small-radius-nonabsorption}
\end{align}
For $d=1$, the logarithmic factors are interpreted as one.
In particular, $p_d^B(r)\to1$ and $p_d^A(r)\to1$ as $r\downarrow0$.
\end{theorem}

\begin{proof}
\emph{Asymptotics of the coefficient functions.}
We first prove that, for every fixed integer $k\ge0$,
\begin{align}
L_k^B(\sigma)
&\sim\frac{\sigma}{\pi k!}
\bigl(\log(1/\sigma)\bigr)^k,
\label{eq:type-B-coefficients-small-scale}
\\
L_k^A(\sigma)
&\sim\frac{2^k\sigma^2}{k!}
\bigl(\log(1/\sigma)\bigr)^k,
\label{eq:type-A-coefficients-small-scale}
\end{align}
as $\sigma\downarrow0$. At $k=0$, both statements follow directly
from \eqref{eq:small-radius-zeroth-coefficients}.

For type $B$ and $k\ge1$, put $b_\sigma:=\arsinh(1/\sigma)$.
For $0\le v\le b_\sigma$,
\[
0\le\sigma\sinh(b_\sigma-v)
=\cosh v-\sqrt{1+\sigma^2}\sinh v
\le\eee^{-v}.
\]
Consequently, the substitution $v=b_\sigma-u$ gives
\begin{equation}
\label{eq:type-B-small-scale-endpoint-bound}
0\le{}
\int_0^{b_\sigma}u^{k-1}
\left(\frac1{\sqrt{1-\sigma^2\sinh^2u}}-1\right)\dint u
\le{}
b_\sigma^{k-1}\int_0^\infty
\left(\frac1{\sqrt{1-\eee^{-2v}}}-1\right)\dint v.
\end{equation}
The last integral is finite: its integrand is $O(v^{-1/2})$
near zero and $O(\eee^{-2v})$ at infinity. Thus
\[
\int_0^{b_\sigma}
\frac{u^{k-1}}{\sqrt{1-\sigma^2\sinh^2u}}\,\dint u
=\frac{b_\sigma^k}{k}+O(b_\sigma^{k-1}).
\]
Since $b_\sigma\sim\log(1/\sigma)$, substitution into
\eqref{eq:type-B-small-radius-positive-coefficients} proves
\eqref{eq:type-B-coefficients-small-scale}.

For type $A$ and $k\ge1$, the elementary identity
$\sinh\theta_\sigma=\sqrt{1+4\sigma^2}/(2\sigma^2)$ gives
\[
L_k^A(\sigma)
=\frac{\sigma^2\theta_\sigma^k}
{k!\sqrt{1+4\sigma^2}}.
\]
Together with $\theta_\sigma\sim2\log(1/\sigma)$, this proves
\eqref{eq:type-A-coefficients-small-scale}.

\medskip
\emph{Gaussian-point non-absorption.}
In \eqref{eq:brownian-gaussian-nonabsorption-finite-sum}, the
largest index is $d-1$. By
\eqref{eq:type-B-coefficients-small-scale}
and~\eqref{eq:type-A-coefficients-small-scale}, every term with
a smaller index is negligible compared with that term.
Consequently,
\begin{align}
1-q_d^B(\sigma)
&\sim\frac{2}{\pi(d-1)!}\,
\sigma\bigl(\log(1/\sigma)\bigr)^{d-1},
\label{eq:type-B-gaussian-nonabsorption-small-scale}
\\
1-q_d^A(\sigma)
&\sim\frac{2^d}{(d-1)!}\,
\sigma^2\bigl(\log(1/\sigma)\bigr)^{d-1},
\label{eq:type-A-gaussian-nonabsorption-small-scale}
\end{align}
as $\sigma\downarrow0$. This includes $d=1$, for which the
finite sum consists only of the zeroth coefficient.

\medskip
\emph{Prescribed-point non-absorption.}
For $\bullet\in\{B,A\}$, set $f_d^\bullet(r):=1-p_d^\bullet(r)$.
These functions are bounded, nonnegative, and nondecreasing,
because the hulls are convex and contain the origin. By
\eqref{eq:gaussian-radial-averaging},
\[
\Ee[f_d^\bullet(\sigma R_d)]=1-q_d^\bullet(\sigma).
\]
The equivalents just proved are regularly varying at zero
with index $1$ in type $B$ and index $2$ in type $A$.
Lemma~\ref{lem:radial-gaussian-regular-variation} therefore yields
\begin{align*}
\frac{1-p_d^B(r)}{1-q_d^B(r)}
&\longrightarrow
\frac{\Gamma(d/2)}{\sqrt2\,\Gamma((d+1)/2)},
\\
\frac{1-p_d^A(r)}{1-q_d^A(r)}
&\longrightarrow
\frac{\Gamma(d/2)}{2\Gamma((d+2)/2)}=\frac1d,
\end{align*}
as $r\downarrow0$. Combining these limits with
\eqref{eq:type-B-gaussian-nonabsorption-small-scale}
and~\eqref{eq:type-A-gaussian-nonabsorption-small-scale}
proves \eqref{eq:type-B-small-radius-nonabsorption}
and~\eqref{eq:type-A-small-radius-nonabsorption}.
\end{proof}

\begin{proof}[Proof of Theorem~\ref{thm:intro-brownian-long-time-nonabsorption}]
For fixed $x\in\R^d\setminus\{0\}$, substitute
$r=\|x\|/\sqrt T$ into
\eqref{eq:type-B-small-radius-nonabsorption}
and~\eqref{eq:type-A-small-radius-nonabsorption}, and use
\eqref{eq:intro-brownian-absorption-scaling}. Since
$\log(1/r)=\log(\sqrt T/\|x\|)$, the resulting equivalents are
exactly \eqref{eq:intro-type-B-long-time-nonabsorption}
and~\eqref{eq:intro-type-A-long-time-nonabsorption}.
\end{proof}

\section{Large-radius asymptotics}
\label{sec:large-radius-asymptotics}

We prove Theorem~\ref{thm:intro-brownian-short-time-absorption}
by deriving large-radius asymptotics from the exact absorption
formulas. The substitutions already used in
\eqref{eq:intro-type-B-prescribed-absorption}
and~\eqref{eq:intro-type-A-prescribed-absorption} make the
rescaling particularly simple. Absolute estimates control the
signed factors throughout the integration range and give a
relative error of order $r^{-2}$ in both types.

\begin{theorem}[Large-radius absorption]
\label{thm:brownian-large-radius-absorption}
For every fixed integer $d\ge1$, as $r\to\infty$,
\begin{align}
p_d^B(r)
&=\frac{2^{d/2}\Gamma(d/2)}{\pi}\,P_d^B(0)\,
 r^{-d}\eee^{-r^2/2}\bigl(1+O(r^{-2})\bigr),
\label{eq:type-B-large-radius-absorption}
\\
p_d^A(r)
&=\frac{\Gamma(d/2)}{2^{(d+1)/2}\sqrt\pi}\,P_d^A(0)\,
 r^{1-d}\eee^{-2r^2}\bigl(1+O(r^{-2})\bigr).
\label{eq:type-A-large-radius-absorption}
\end{align}
The constants implicit in the error terms may depend on $d$.
Both leading constants are strictly positive.
\end{theorem}

\begin{proof}
\emph{The leading constants and bounds on the factors.}
We first verify the integral representations
\eqref{eq:intro-type-B-factor-at-zero}
and~\eqref{eq:intro-type-A-factor-at-zero}.
The defining series for $\Psi_m$ shows that the functions
$t\mapsto\mathrm i^{-d}\Psi_d(\mathrm it)$ and
$t\mapsto\mathrm i^{-(d+1)}\Psi_{d+1}(\mathrm it)$ are real-valued
for real $t$. Their derivatives at zero of orders
$0,\ldots,d-1$ vanish, and their $d$th derivatives are
$\cos t$ and $\sin t$, respectively. Moreover, the definitions
of $P_d^B,P_d^A$ give
\[
P_d^B(0)=\mathrm i^{-d}\Psi_d(\pi\mathrm i/2),
\qquad
P_d^A(0)=\mathrm i^{-(d+1)}\Psi_{d+1}(\pi\mathrm i).
\]
Taylor's formula with integral remainder, evaluated at
$t=\pi/2$ and $t=\pi$, therefore gives exactly
\eqref{eq:intro-type-B-factor-at-zero}
and~\eqref{eq:intro-type-A-factor-at-zero}. Their integrands
are strictly positive in the interiors of the respective
integration intervals. Hence $P_d^B(0)>0$ and $P_d^A(0)>0$,
including when $d=1$.

By \eqref{eq:prescribed-factors-growth}, the inequality
$\arsinh v\le v$, and the identities
$\cosh(\arsinh v)=\sqrt{1+v^2}$ and
$\cosh(2\arsinh v)=1+2v^2$, the two functions
\[
v\longmapsto
\frac{P_d^B(\arsinh v)}{(1+v^2)^{(d+1)/2}},
\qquad
v\longmapsto
\frac{P_d^A(2\arsinh v)}{(1+v^2)^{(d+1)/2}}
\]
are bounded on $[0,\infty)$.
Proposition~\ref{prop:prescribed-factors-elementary-structure}
also shows that they extend to even analytic functions in a
neighborhood of zero. Their differences from their values at
zero are therefore $O(v^2)$ as $v\downarrow0$. Boundedness
extends this estimate to all $v\ge0$: there is a constant
$C_d<\infty$ such that
\begin{equation}
\label{eq:large-radius-factor-quadratic-bounds}
\begin{aligned}
\left|
\frac{P_d^B(\arsinh v)}{(1+v^2)^{(d+1)/2}}-P_d^B(0)
\right|&\le C_d v^2,\\
\left|
\frac{P_d^A(2\arsinh v)}{(1+v^2)^{(d+1)/2}}-P_d^A(0)
\right|&\le C_d v^2,
\qquad v\ge0.
\end{aligned}
\end{equation}

\medskip
\emph{Type $B$.}
The substitution $t=y/x$ in
\eqref{eq:intro-type-B-analytic-kernel} gives
\[
x^{d/2}\eee^xG_d^B(x)
=\int_0^\infty
 y^{d/2-1}\eee^{-y}(1+y/x)^{-1/2}\,\dint y.
\]
Since $0\le1-(1+a)^{-1/2}\le a/2$ for $a\ge0$, the gamma
integral yields the explicit bound
\begin{equation}
\label{eq:type-B-kernel-large-argument-bound}
0\le
1-\frac{x^{d/2}\eee^xG_d^B(x)}{\Gamma(d/2)}
\le\frac{d}{4x},
\qquad x>0.
\end{equation}
Thus $G_d^B(x)=\Gamma(d/2)x^{-d/2}\eee^{-x}
(1+O(x^{-1}))$. At $x=r^2(1+v^2)/2$, the absolute relative
error in this approximation is at most $d/(2r^2)$,
uniformly for $v\ge0$.

Insert this approximation into
\eqref{eq:intro-type-B-prescribed-absorption} and set $v=y/r$.
The boundedness proved above and
\eqref{eq:type-B-kernel-large-argument-bound}, applied with
absolute values to the error integral, give
\begin{equation}
\label{eq:type-B-large-radius-Gaussian-rescaling}
\begin{aligned}
r^d\eee^{r^2/2}p_d^B(r)
={}&\frac{2^{(d+1)/2}\Gamma(d/2)}{\pi^{3/2}}
\int_0^\infty
\frac{P_d^B(\arsinh(y/r))}
{(1+y^2/r^2)^{(d+1)/2}}
\eee^{-y^2/2}\,\dint y
+O(r^{-2}).
\end{aligned}
\end{equation}
Indeed, after the rescaling the absolute kernel error is
bounded by a constant times
$r^{-2}\int_0^\infty\eee^{-y^2/2}\,\dint y$.
By \eqref{eq:large-radius-factor-quadratic-bounds}, replacing
the remaining fraction in
\eqref{eq:type-B-large-radius-Gaussian-rescaling} by $P_d^B(0)$
produces an error bounded by a constant times
$r^{-2}\int_0^\infty y^2\eee^{-y^2/2}\,\dint y$.
Consequently,
\[
r^d\eee^{r^2/2}p_d^B(r)
=\frac{2^{(d+1)/2}\Gamma(d/2)}{\pi^{3/2}}
 P_d^B(0)\sqrt{\frac\pi2}+O(r^{-2})
=\frac{2^{d/2}\Gamma(d/2)}\pi P_d^B(0)+O(r^{-2}).
\]
The positivity of $P_d^B(0)$ converts this additive error
into the relative error in
\eqref{eq:type-B-large-radius-absorption}.

\medskip
\emph{Type $A$.}
Here the Gaussian exponential is already explicit. The
substitution $v=y/r$ in
\eqref{eq:intro-type-A-prescribed-absorption} gives the exact
identity
\begin{equation}
\label{eq:type-A-large-radius-Gaussian-rescaling}
r^{d-1}\eee^{2r^2}p_d^A(r)
=\frac{2^{1-d/2}\Gamma(d/2)}\pi
\int_0^\infty
\frac{P_d^A(2\arsinh(y/r))}
{(1+y^2/r^2)^{(d+1)/2}}
\eee^{-2y^2}\,\dint y.
\end{equation}
The second bound in
\eqref{eq:large-radius-factor-quadratic-bounds} controls the
error after replacing the fraction by $P_d^A(0)$ by a constant
times $r^{-2}\int_0^\infty y^2\eee^{-2y^2}\,\dint y$.
Since $\int_0^\infty\eee^{-2y^2}\,\dint y
=\sqrt\pi/(2\sqrt2)$, we obtain
\[
r^{d-1}\eee^{2r^2}p_d^A(r)
=\frac{\Gamma(d/2)}{2^{(d+1)/2}\sqrt\pi}
 P_d^A(0)+O(r^{-2}).
\]
Again the leading constant is positive, proving
\eqref{eq:type-A-large-radius-absorption}.
\end{proof}

\begin{proof}[Proof of Theorem~\ref{thm:intro-brownian-short-time-absorption}]
For fixed $x\in\R^d\setminus\{0\}$, substitute
$r=\|x\|/\sqrt T$ into
\eqref{eq:type-B-large-radius-absorption}
and~\eqref{eq:type-A-large-radius-absorption}, and use
\eqref{eq:intro-brownian-absorption-scaling}.
This gives \eqref{eq:intro-type-B-short-time-absorption}
and~\eqref{eq:intro-type-A-short-time-absorption}.
The relative errors are in fact $O(T/\|x\|^2)$, and hence tend
to zero as $T\downarrow0$. The positivity of both leading
constants was proved above.
\end{proof}

\section{First absorption times for Brownian convex hulls}
\label{sec:first-absorption-times}
In this section we prove the transform formulas and the
high-dimensional results stated in
Subsections~\ref{subsec:intro-first-absorption-times-exact}
and~\ref{subsec:intro-first-absorption-times-high-dim}.
Throughout, $\tau_d(r)$ is the first absorption time of $re_1$
by the increasing hulls of a single standard Brownian motion
in $\R^d$, as defined in
\eqref{eq:intro-type-B-first-absorption-time}.
Its distribution, tail equivalents, and moment-finiteness
criterion follow from the preceding sections by the arguments
already given in the Introduction. We first evaluate its
Mellin and Laplace transforms, and then pass from the
high-dimensional Gaussian-point transition in~\cite{KabluchkoOUSpirals} to the
prescribed-point problem.

\subsection{Mellin and Laplace transforms of the first absorption time}
\label{subsec:first-absorption-time-transforms}

Recall that
$\Pp(\tau_d(r)\le T)=p_d^B(r/\sqrt T)$ and
$\tau_d(r)\overset{\mathrm d}=r^2\tau_d(1)$.
The relation with the meromorphic continuation of $M_d^B$
was established in
\eqref{eq:intro-type-B-first-absorption-radial-Mellin}.
For the explicit generating formula, we instead work directly
with Gaussian-point non-absorption and the finite sum
\eqref{eq:brownian-gaussian-nonabsorption-finite-sum}.
All beta integrals in this calculation converge in their
ordinary sense.

\begin{proof}[Proof of Theorem~\ref{thm:intro-type-B-first-absorption-transforms}]
\emph{Gaussian radial averaging and positive Mellin parameters.}
For $s\in\mathbb C$ with $0<\operatorname{Re}s<1/2$,
integration of the upper tail and the substitution
$u=T^{-1/2}$ give
\[
\Ee\bigl[(\tau_d(1))^s\bigr]
=2s\int_0^\infty u^{-2s-1}\bigl(1-p_d^B(u)\bigr)\,\dint u.
\]
The integral converges absolutely by
\eqref{eq:type-B-small-radius-nonabsorption} at zero and
$0\le1-p_d^B\le1$ at infinity.
Let $R_d$ have the $\chi_d$ distribution, as in
\eqref{eq:gaussian-radial-averaging}. Integration against its
density gives
$\Ee[R_d^{2s}]=2^s\Gamma(d/2+s)/\Gamma(d/2)$.
Since
$1-q_d^B(\sigma)=\Ee[1-p_d^B(\sigma R_d)]$,
Fubini's theorem and $u=\sigma R_d$ yield
\begin{equation}
\label{eq:type-B-absorption-time-Mellin-radial-average}
\frac{2^s\Gamma(d/2+s)}{\Gamma(d/2)}
\Ee\bigl[(\tau_d(1))^s\bigr]
=2s\int_0^\infty
\sigma^{-2s-1}\bigl(1-q_d^B(\sigma)\bigr)\,\dint\sigma.
\end{equation}
To justify this exchange, replace $s$ by
$c:=\operatorname{Re}s$ in the absolute values. Tonelli's
theorem gives the product of $\Ee[R_d^{2c}]$ and the preceding
finite integral at $c$.

We evaluate the integrals of the coefficient functions
$L_k^B$ occurring in
\eqref{eq:brownian-gaussian-nonabsorption-finite-sum}.
For the zeroth coefficient, integration by parts and the beta
integral give
\begin{equation}
\label{eq:type-B-absorption-time-zeroth-coefficient-integral}
\begin{aligned}
\int_0^\infty \sigma^{-2s-1}L_0^B(\sigma)\,\dint\sigma
&=\frac1{2\pi s}\int_0^\infty
\frac{\sigma^{-2s}}{1+\sigma^2}\,\dint\sigma
=\frac1{4s\cos(\pi s)}.
\end{aligned}
\end{equation}
The boundary terms vanish because $0<\operatorname{Re}s<1/2$.
For $k\ge1$, use
\eqref{eq:type-B-small-radius-positive-coefficients} and
integrate first in $\sigma$. For fixed $u>0$, the substitution
$t=\sigma\sinh u$ gives
\[
\int_0^{1/\sinh u}
\frac{\sigma^{-2s}}{\sqrt{1-\sigma^2\sinh^2u}}\,\dint\sigma
=\frac12\mathrm B\!\left(\frac12-s,\frac12\right)
(\sinh u)^{2s-1}.
\]
Consequently,
\begin{equation}
\label{eq:type-B-absorption-time-positive-coefficient-integral}
\int_0^\infty \sigma^{-2s-1}L_k^B(\sigma)\,\dint\sigma
=\frac{\mathrm B(1/2-s,1/2)}{2\pi(k-1)!}
\int_0^\infty u^{k-1}(\sinh u)^{2s-1}\,\dint u.
\end{equation}
These exchanges are absolutely justified by the same calculation
at $c=\operatorname{Re}s$: the last integrand has order
$u^{k+2c-2}$ at zero and decays exponentially at infinity.

\medskip
\emph{Summation in the dimension.}
Let $0<\operatorname{Re}s<1/2$ and
$|w|<1-2\operatorname{Re}s$, and put $c:=\operatorname{Re}s$
and $\eta:=|w|$. The nonnegative coefficient functions satisfy
\begin{equation}
\label{eq:type-B-absorption-time-coefficient-sum-bound}
\sum_{k=1}^\infty \eta^k
\int_0^\infty\sigma^{-2c-1}L_k^B(\sigma)\,\dint\sigma
=
\frac{\eta\,\mathrm B(1/2-c,1/2)}{2\pi}
\int_0^\infty \eee^{\eta u}(\sinh u)^{2c-1}\,\dint u
<\infty.
\end{equation}
Indeed, the last integrand has order $u^{2c-1}$ at zero and
order $\eee^{-(1-2c-\eta)u}$ at infinity.
Together with
\eqref{eq:type-B-absorption-time-zeroth-coefficient-integral},
this justifies all the following summations and integrations
absolutely. In the finite complementary sum, a fixed coefficient
$L_k^B$ occurs precisely in dimensions $d=k+1+2j$, $j\ge0$.
Thus \eqref{eq:type-B-absorption-time-Mellin-radial-average}
gives
\begin{equation}
\label{eq:type-B-absorption-time-Mellin-series-integral}
\begin{aligned}
&\sum_{d=1}^\infty
\frac{2^s\Gamma(d/2+s)}{\Gamma(d/2)}
\Ee\bigl[(\tau_d(1))^s\bigr]w^{d-1}
\\
&\qquad=
\frac{4s}{1-w^2}
\sum_{k=0}^\infty w^k
\int_0^\infty\sigma^{-2s-1}L_k^B(\sigma)\,\dint\sigma
\\
&\qquad=
\frac1{1-w^2}\left[
\frac1{\cos(\pi s)}
+\frac{2sw}{\pi}\mathrm B\!\left(\frac12-s,\frac12\right)
\int_0^\infty\eee^{wu}(\sinh u)^{2s-1}\,\dint u
\right].
\end{aligned}
\end{equation}

The convergence is locally uniform in $(s,w)$ on the stated
domain. To see this, work in a sufficiently small closed
neighborhood on which
$0<c_-\le\operatorname{Re}s\le c_+<1/2$ and
$|w|\le\eta_0<1-2c_+$. For $\sigma>0$,
\[
|\sigma^{-2s-1}|
\le\sigma^{-2c_--1}+\sigma^{-2c_+-1}.
\]
Using \eqref{eq:type-B-absorption-time-Mellin-radial-average},
the finite complementary sum, and
\eqref{eq:type-B-absorption-time-coefficient-sum-bound} at
$c_-$ and $c_+$ gives a summable majorant for the dimension
series on this neighborhood. The factor $|s|$ is bounded there.
This proves absolute and locally uniform convergence, including
for complex parameters.

To evaluate the remaining integral, set $t=\eee^{-2u}$:
\begin{equation}
\label{eq:type-B-absorption-time-hyperbolic-beta-integral}
\int_0^\infty\eee^{wu}(\sinh u)^{2s-1}\,\dint u
=2^{-2s}\mathrm B\!\left(\frac{1-2s-w}{2},2s\right).
\end{equation}
Both beta parameters have positive real parts. The reflection
and duplication formulas for the gamma function give
\[
\frac{2^{1-2s}s}{\pi}
\mathrm B\!\left(\frac12-s,\frac12\right)\Gamma(2s)
=\frac{\Gamma(1+s)}{\Gamma(1-s)\cos(\pi s)};
\]
see \cite[Equations~(5.5.3) and~(5.5.5)]{NIST:DLMF}.
Substituting this identity and
\eqref{eq:type-B-absorption-time-hyperbolic-beta-integral}
into \eqref{eq:type-B-absorption-time-Mellin-series-integral}
proves
\eqref{eq:intro-type-B-first-absorption-Mellin-generating-function}.

\medskip
\emph{The Laplace transform.}
Fix $d\ge1$ and $r,\lambda>0$, and put $a:=r\sqrt{2\lambda}$.
The distribution identity
\eqref{eq:intro-type-B-first-absorption-distribution} and
Tonelli's theorem give
\begin{equation}
\label{eq:type-B-absorption-time-Laplace-distribution-integral}
\Ee\bigl[\eee^{-\lambda\tau_d(r)}\bigr]
=\lambda\int_0^\infty
\eee^{-\lambda T}p_d^B(r/\sqrt T)\,\dint T.
\end{equation}
We insert the hyperbolic representation
\eqref{eq:type-B-prescribed-absorption-hyperbolic} and the
integral defining $G_d^B$ in
\eqref{eq:intro-type-B-analytic-kernel}.
The time integral is evaluated by
\begin{equation}
\label{eq:absorption-time-Gaussian-reciprocal-integral}
\int_0^\infty T^{-1/2}\eee^{-\lambda T-b/T}\,\dint T
=\sqrt{\frac\pi\lambda}\,\eee^{-2\sqrt{\lambda b}},
\qquad b,\lambda>0,
\end{equation}
which follows from the Bessel integral and the value of
$K_{-1/2}$ in
\cite[Equations~(10.32.10) and~(10.39.2)]{NIST:DLMF}.
After this integration, set $y=\sqrt{1+t}$ in the kernel
integral. The result is
\begin{equation}
\label{eq:type-B-absorption-time-Laplace-double-integral}
\Ee\bigl[\eee^{-\lambda\tau_d(r)}\bigr]
=\frac{2a}{\pi}\int_0^\infty P_d^B(u)
\int_1^\infty (y^2-1)^{d/2-1}
\eee^{-ay\cosh u}\,\dint y\,\dint u.
\end{equation}

For completeness, the exchanges leading to this identity do
not rely on cancellation in $P_d^B$. Replace that factor by
$|P_d^B|$ and apply Tonelli's theorem. The same time integration
and substitution reduce absolute integrability to the
finiteness of
\[
\int_0^\infty\int_1^\infty
|P_d^B(u)|(y^2-1)^{d/2-1}\eee^{-ay\cosh u}
\,\dint y\,\dint u.
\]
For $y\ge1$ and $u\ge0$, we have
$y\cosh u\ge(y+\cosh u)/2$.
Hence this integral is at most
\[
\left(\int_0^\infty |P_d^B(u)|\eee^{-(a/2)\cosh u}\,\dint u\right)
\left(\int_1^\infty (y^2-1)^{d/2-1}\eee^{-ay/2}\,\dint y\right)
<\infty.
\]
The first factor is finite by
\eqref{eq:prescribed-factors-growth}. The second is integrable
at $y=1$ because $d/2-1>-1$, and at infinity because of its
exponential factor. Fubini's theorem is therefore applicable.

The classical Bessel representation
\cite[Equation~(10.32.8)]{NIST:DLMF} gives, for $b>0$,
\begin{equation}
\label{eq:absorption-time-Bessel-integral}
\int_1^\infty (y^2-1)^{d/2-1}\eee^{-by}\,\dint y
=\frac{\Gamma(d/2)}{\sqrt\pi}
\left(\frac2b\right)^{(d-1)/2}K_{(d-1)/2}(b).
\end{equation}
Apply this with $b=a\cosh u$ in
\eqref{eq:type-B-absorption-time-Laplace-double-integral}.
It follows that
\[
\Ee\bigl[\eee^{-\lambda\tau_d(r)}\bigr]
=\frac{2^{(d+1)/2}\Gamma(d/2)}{\pi^{3/2}}a^{(3-d)/2}
\int_0^\infty
\frac{P_d^B(u)}{\cosh^{(d-1)/2}u}
K_{(d-1)/2}(a\cosh u)\,\dint u.
\]
Finally, $v=\sinh u$ gives exactly
\eqref{eq:intro-type-B-first-absorption-Laplace-transform}.
The preceding absolute-integrability argument also proves
absolute convergence of that integral, completing the proof.
\end{proof}

\subsection{First absorption of a prescribed point by high-dimensional Brownian motion}
\label{subsec:first-absorption-high-dimensions}

We prove Theorem~\ref{thm:intro-type-B-first-absorption-high-dimension}
by combining the high-dimensional Ornstein--Uhlenbeck
absorption theorem from the companion
paper~\cite[Theorems~5.8 and~5.14]{KabluchkoOUSpirals}
with concentration of the norm of an independent standard
Gaussian vector around $\sqrt d$.
Monotonicity of the prescribed-point absorption probabilities
then allows us to pass from Gaussian-point to prescribed-point
absorption.

\begin{proof}[Proof of Theorem~\ref{thm:intro-type-B-first-absorption-high-dimension}]
The ordinary Ornstein--Uhlenbeck correspondence and the
high-dimensional absorption theorem in the companion paper~\cite[Theorems~5.8 and~5.14]{KabluchkoOUSpirals}
have the following consequence. For any positive sequence
$(\sigma_d)$ and any $y\in\R$,
\begin{equation}
\label{eq:type-B-high-dimensional-Gaussian-point-input}
\frac{\log(1+\sigma_d^{-2})-2d}{2\sqrt d}\longrightarrow y
\quad\Longrightarrow\quad
q_d^B(\sigma_d)\longrightarrow\Phi(y).
\end{equation}
Indeed, the correspondence identifies $q_d^B(\sigma_d)$ with
the origin-absorption probability of the stationary
Ornstein--Uhlenbeck process with coordinate covariance
$\eee^{-|t-s|/2}$ over a time interval of length
$\log(1+\sigma_d^{-2})$. The cited transition theorem applies
to precisely this normalization.

Fix $y\in\R$, and put
\[
t_d:=\eee^{2d+2y\sqrt d},
\qquad
r_d:=t_d^{-1/2}=\eee^{-d-y\sqrt d}.
\]
By \eqref{eq:intro-type-B-first-absorption-distribution},
$\Pp(\tau_d(1)\le t_d)=p_d^B(r_d)$.
For every fixed $c>0$,
\[
\log\left(1+\left(\frac{cr_d}{\sqrt d}\right)^{-2}\right)
=2d+2y\sqrt d+\log d-2\log c+o(1).
\]
Since $\log d=o(\sqrt d)$,
\eqref{eq:type-B-high-dimensional-Gaussian-point-input} gives
\begin{equation}
\label{eq:type-B-critical-Gaussian-radius-limit}
q_d^B\!\left(\frac{cr_d}{\sqrt d}\right)
\longrightarrow\Phi(y),
\qquad c>0.
\end{equation}

Let $Z_d\sim N(0,I_d)$ be independent of the Brownian motion
in dimension $d$, and put $R_d:=\|Z_d\|$.
The identities $\Ee[R_d^2]=d$ and
$\operatorname{Var}(R_d^2)=2d$, followed by Chebyshev's
inequality, imply $R_d/\sqrt d\to1$ in probability.
Fix $0<\varepsilon<1$. Using
$q_d^B(\sigma)=\Ee[p_d^B(\sigma R_d)]$, the monotonicity of
$p_d^B$, and $0\le p_d^B\le1$, we obtain
\begin{equation}
\label{eq:type-B-high-dimensional-radial-squeeze}
\begin{aligned}
q_d^B\!\left(\frac{r_d}{(1-\varepsilon)\sqrt d}\right)
-\Pp\bigl(R_d<(1-\varepsilon)\sqrt d\bigr)
&\le p_d^B(r_d),
\\
p_d^B(r_d)
&\le q_d^B\!\left(\frac{r_d}{(1+\varepsilon)\sqrt d}\right)
+\Pp\bigl(R_d>(1+\varepsilon)\sqrt d\bigr).
\end{aligned}
\end{equation}
For the first inequality, on the event
$R_d\ge(1-\varepsilon)\sqrt d$ the argument of $p_d^B$ in the
Gaussian average is at least $r_d$; on its complement we bound
the probability by one. For the second inequality, on
$R_d\le(1+\varepsilon)\sqrt d$ that argument is at most $r_d$,
so the Gaussian average is at least
$p_d^B(r_d)\Pp(R_d\le(1+\varepsilon)\sqrt d)
\ge p_d^B(r_d)-\Pp(R_d>(1+\varepsilon)\sqrt d)$.

Both error probabilities in
\eqref{eq:type-B-high-dimensional-radial-squeeze} tend to zero.
By \eqref{eq:type-B-critical-Gaussian-radius-limit}, with
$c=(1-\varepsilon)^{-1}$ and $c=(1+\varepsilon)^{-1}$,
both Gaussian-point probabilities tend to $\Phi(y)$.
Thus $p_d^B(r_d)\to\Phi(y)$, which is
\eqref{eq:intro-type-B-first-absorption-Gaussian-window}.
Since $y\in\R$ was arbitrary and $\Phi$ is continuous, this
proves \eqref{eq:intro-type-B-first-absorption-logarithmic-CLT}.
\end{proof}

\subsection{First absorption of the unit ball by high-dimensional Brownian motion}
\label{subsec:first-absorption-unit-ball-high-dim}

We deduce the central limit theorem for the unit-ball absorption
time from the Gaussian-point transition
\eqref{eq:type-B-high-dimensional-Gaussian-point-input}
and the prescribed-point central limit theorem, using independence
of Brownian increments and a Gaussian translation estimate.
Convexity allows us to pass from $2d$ coordinate shifts to
containment of the entire unit ball.

\begin{proof}[Proof of Theorem~\ref{thm:intro-type-B-unit-ball-high-dimension}]
The convex hull of the $2d$ points $\pm\sqrt d\,e_i$,
$1\le i\le d$, contains $\mathbb B_d$, since
$\sum_{i=1}^d|x_i|\le\sqrt d\,\|x\|$.
Each of these finitely many points is absorbed almost surely
in finite time, by rotational invariance and the almost-sure
finiteness of $\tau_d(\sqrt d)$.
Consequently, $0<\tau_d(1)\le \tau_d^{\mathrm{ball}}<\infty$ almost surely.

Fix $d\ge1$ and a deterministic time $T>d^4$, and put
\[
Z_d:=W_{d^4},
\qquad
K':=\conv\{W_{d^4+t}-W_{d^4}:0\le t\le T-d^4\}.
\]
Then $Z_d\sim N(0,d^4I_d)$ is independent of $K'$, and
$Z_d+K'\subset\mathcal K_{d,T}^B$.
For every $i\in\{1,\ldots,d\}$ and
$\eta\in\{-1,1\}$, we have
\begin{equation}
\label{eq:type-B-unit-ball-Gaussian-shift-bound}
\Pp\bigl(-Z_d\in K',\,
        -Z_d+\eta\sqrt d\,e_i\notin K'\bigr)
\le\frac{\sqrt d}{\sqrt{2\pi}\,d^2}.
\end{equation}
Indeed, condition on $K'$ and all coordinates of $Z_d$
except the $i$-th. The section of $K'$ along the remaining
coordinate line is an interval, possibly empty or degenerate.
The coordinates of points that belong to this interval but leave
it after translation by $\eta\sqrt d$ form an interval of
length at most $\sqrt d$.
The remaining Gaussian coordinate has density bounded by
$1/(\sqrt{2\pi}\,d^2)$, which proves the bound.

If all the points $-Z_d\pm\sqrt d\,e_i$ belong to $K'$, then
$\mathbb B_d\subset Z_d+K'\subset\mathcal K_{d,T}^B$.
The union bound and
\eqref{eq:type-B-unit-ball-Gaussian-shift-bound} therefore give
\begin{align}
\Pp(\tau_d^{\mathrm{ball}}\le T)
\ge\Pp(-Z_d\in K')-\sqrt{\frac{2}{\pi d}}
=q_d^B\!\left(\frac{d^2}{\sqrt{T-d^4}}\right)
-\sqrt{\frac{2}{\pi d}},
\label{eq:type-B-unit-ball-Gaussian-point-comparison}
\end{align}
where the equality follows from Brownian scaling and independence.

Now fix $y\in\R$ and put $T_d:=\eee^{2d+2y\sqrt d}$.
For all sufficiently large $d$, we have $T_d>d^4$, and
\[
\frac{\log(1+(T_d-d^4)/d^4)-2d}{2\sqrt d}
=y-\frac{2\log d}{\sqrt d}\longrightarrow y.
\]
Thus \eqref{eq:type-B-high-dimensional-Gaussian-point-input}
implies
\[
q_d^B\!\left(\frac{d^2}{\sqrt{T_d-d^4}}\right)
\longrightarrow\Phi(y).
\]
Together with
\eqref{eq:type-B-unit-ball-Gaussian-point-comparison}, this yields
$\liminf_{d\to\infty}\Pp(\tau_d^{\mathrm{ball}}\le T_d)\ge\Phi(y)$.
On the other hand, $\tau_d(1)\le \tau_d^{\mathrm{ball}}$ almost surely, so
Theorem~\ref{thm:intro-type-B-first-absorption-high-dimension}
gives
\[
\limsup_{d\to\infty}\Pp(\tau_d^{\mathrm{ball}}\le T_d)
\le\lim_{d\to\infty}\Pp(\tau_d(1)\le T_d)
=\Phi(y).
\]
Since $y\in\R$ was arbitrary, this proves
\eqref{eq:intro-type-B-unit-ball-logarithmic-CLT}.

Finally, put
\[
X_d:=\frac{\log \tau_d^{\mathrm{ball}}-2d}{2\sqrt d},
\qquad
Y_d:=\frac{\log\tau_d(1)-2d}{2\sqrt d}.
\]
Both sequences converge in distribution to $N(0,1)$,
and $X_d\ge Y_d$ almost surely. Since $\arctan$ is bounded,
continuous, and increasing,
\[
0\le\Ee[\arctan X_d-\arctan Y_d]\longrightarrow0.
\]
For each $M>0$, the mean value theorem gives
\[
0\le X_d-Y_d
\le(1+M^2)(\arctan X_d-\arctan Y_d)
\quad\text{on }\{|X_d|\le M,\ |Y_d|\le M\}.
\]
Markov's inequality and tightness of $(X_d)$ and $(Y_d)$
therefore imply $X_d-Y_d\to0$ in probability, proving
\eqref{eq:intro-type-B-unit-ball-point-logarithmic-difference}.
\end{proof}

\subsection{High-dimensional absorption by Brownian bridges}
\label{subsec:high-dimensional-bridge-absorption}

We first prove the prescribed-point limit using the exact
Gaussian-point absorption formula and the Poisson central limit
theorem. For unit-ball absorption, we extract an independent
Gaussian translation from the middle part of the bridge and
apply the high-dimensional transition for the ordinary
Ornstein--Uhlenbeck process.

\begin{proof}[Proof of Theorem~\ref{thm:intro-type-A-high-dimensional-absorption}]
\emph{Prescribed-point absorption.}
Put $r_d:=T_d^{-1/2}$. Fix $c>0$ and set
\[
\theta_d:=2\arsinh\left(\frac{\sqrt{dT_d}}{2c}\right).
\]
The assumption on $T_d$ gives
$\theta_d=\log(dT_d/c^2)+o(1)$, and hence
$\theta_d/d\to1$ and
$(\theta_d-d)/\sqrt d\to y$.
Let $N_d$ have the Poisson distribution with mean $\theta_d$.
By \eqref{eq:type-A-gaussian-point-absorption}
and~\eqref{eq:intro-hyperbolic-remainders},
\[
q_d^A\left(\frac{cr_d}{\sqrt d}\right)
=
\frac{\Pp(N_d\ge d+1,\ N_d\equiv d+1\pmod 2)}
{\Pp(N_d\text{ is odd})}.
\]
The Poisson probabilities form a unimodal sequence.
Splitting the alternating sum at a mode and using Stirling's
formula therefore gives
\[
\left|
\sum_{j=d+1}^\infty(-1)^{j-d-1}\Pp(N_d=j)
\right|
\le2\max_{j\ge0}\Pp(N_d=j)
=O(d^{-1/2}).
\]
Since $\Pp(N_d\text{ is odd})=(1-\eee^{-2\theta_d})/2$,
it follows that
\[
q_d^A\left(\frac{cr_d}{\sqrt d}\right)
=\Pp(N_d\ge d+1)+o(1).
\]
The Poisson central limit theorem and
$(d+1-\theta_d)/\sqrt{\theta_d}\to-y$ now imply
\begin{equation}
\label{eq:type-A-high-dimensional-gaussian-point-limit}
q_d^A\left(\frac{cr_d}{\sqrt d}\right)
\longrightarrow\Phi(y),\qquad c>0.
\end{equation}

Let $G_d\sim N(0,I_d)$ be independent of the bridge and put
$R_d:=\|G_d\|$. Since $\Ee[R_d^2]=d$ and
$\operatorname{Var}(R_d^2)=2d$, we have
$R_d/\sqrt d\to1$ in probability.
Fix $0<\varepsilon<1$. By
\eqref{eq:gaussian-radial-averaging} and the monotonicity of
$p_d^A$, which follows from convexity and $0\in\mathcal K_d^A$,
\begin{align*}
q_d^A\left(\frac{r_d}{(1-\varepsilon)\sqrt d}\right)
-\Pp\bigl(R_d<(1-\varepsilon)\sqrt d\bigr)
\le p_d^A(r_d)
\le
q_d^A\left(\frac{r_d}{(1+\varepsilon)\sqrt d}\right)
+\Pp\bigl(R_d>(1+\varepsilon)\sqrt d\bigr).
\end{align*}
Both Gaussian-point probabilities converge to $\Phi(y)$ by
\eqref{eq:type-A-high-dimensional-gaussian-point-limit},
while the error probabilities tend to zero.
Thus $p_d^A(r_d)\to\Phi(y)$, and
\eqref{eq:intro-brownian-absorption-scaling} yields
\begin{equation}
\label{eq:type-A-high-dimensional-point-limit-proof}
\Pp(e_1\in\mathcal K_{d,T_d}^A)
=p_d^A(r_d)\longrightarrow\Phi(y).
\end{equation}

\medskip
\emph{Unit-ball absorption.}
For all sufficiently large $d$, we have $T_d>2d^4$.
Let $(\beta_t)_{0\le t\le T_d}$ be a standard Brownian bridge
in $\R^d$ from $0$ to $0$, and put
\[
Z_d:=\frac{\beta_{d^4}+\beta_{T_d-d^4}}2,
\qquad
K_d':=\conv\{\beta_t-Z_d:d^4\le t\le T_d-d^4\}.
\]
For each coordinate $i$, the covariance identity
$\Ee[\beta_s^{(i)}\beta_t^{(i)}]=\min(s,t)-st/T_d$ gives
\[
\operatorname{Var}(Z_d^{(i)})
=\operatorname{Cov}(Z_d^{(i)},\beta_t^{(i)})
=\frac{d^4}{2},
\qquad d^4\le t\le T_d-d^4.
\]
Since the coordinates are independent and the variables are
jointly Gaussian, $Z_d\sim N(0,(d^4/2)I_d)$ is independent of
$(\beta_t-Z_d)_{d^4\le t\le T_d-d^4}$, and hence of $K_d'$.
Moreover,
\[
Z_d+K_d'
=\conv\{\beta_t:d^4\le t\le T_d-d^4\}
\subset\mathcal K_{d,T_d}^A.
\]

For every $i\in\{1,\ldots,d\}$ and $\eta\in\{-1,1\}$,
\[
\Pp\bigl(-Z_d\in K_d',\,
-Z_d+\eta\sqrt d\,e_i\notin K_d'\bigr)
\le\frac{\sqrt d}{\sqrt\pi\,d^2}.
\]
Indeed, condition on $K_d'$ and all coordinates of $Z_d$
except the $i$-th. The section of $K_d'$ along the remaining
coordinate line is an interval, possibly empty or degenerate.
The points of this interval that leave it after a translation
by $\eta\sqrt d$ form a set of length at most $\sqrt d$,
and the remaining Gaussian density is bounded by
$1/(\sqrt\pi\,d^2)$.

The convex hull of the $2d$ points $\pm\sqrt d\,e_i$
contains $\mathbb B_d$, since
$\sum_{i=1}^d|x_i|\le\sqrt d\,\|x\|$.
Consequently, the union bound and the preceding estimate give
\begin{equation}
\label{eq:type-A-high-dimensional-ball-comparison}
\Pp\bigl(0\in\conv\{\beta_t:d^4\le t\le T_d-d^4\}\bigr)
-\frac{2}{\sqrt{\pi d}}
\le
\Pp(\mathbb B_d\subset\mathcal K_{d,T_d}^A)
\le p_d^A(r_d).
\end{equation}

To evaluate the probability on the left, define
\[
u(t):=\log\frac{t}{T_d-t},
\qquad
U_{u(t)}:=\frac{\beta_t}{\sqrt{t(T_d-t)/T_d}},
\qquad 0<t<T_d.
\]
This is a centered stationary Ornstein--Uhlenbeck process
with independent coordinates and covariance
\[
\Ee[U_u^{(i)}U_v^{(i)}]=\eee^{-|u-v|/2},
\qquad u,v\in\R,\quad 1\le i\le d.
\]
Indeed, for $0<s\le t<T_d$, the normalized covariance is
$\sqrt{s(T_d-t)/(t(T_d-s))}
=\eee^{-(u(t)-u(s))/2}$.
The retained interval corresponds to an OU interval of length
\[
L_d:=u(T_d-d^4)-u(d^4)
=2\log\frac{T_d-d^4}{d^4}.
\]
Multiplying the points of a set by individual positive scalars
does not change whether their convex hull contains the origin.
By stationarity, we therefore have
\[
\Pp\bigl(0\in\conv\{\beta_t:d^4\le t\le T_d-d^4\}\bigr)
=
\Pp\bigl(0\in\conv\{U_u:0\le u\le L_d\}\bigr).
\]
Furthermore,
\[
\frac{L_d-2d}{2\sqrt d}
=
\frac{\log(T_d-d^4)-4\log d-d}{\sqrt d}
\longrightarrow y.
\]
The ordinary Ornstein--Uhlenbeck transition theorem
from~\cite[Theorem~5.14]{KabluchkoOUSpirals} shows that the
last origin-absorption probability converges to $\Phi(y)$.
Combining this with
\eqref{eq:type-A-high-dimensional-ball-comparison}
and~\eqref{eq:type-A-high-dimensional-point-limit-proof}
proves the unit-ball limit and hence
\eqref{eq:intro-type-A-high-dimensional-point-ball-absorption}.
The subsequent assertion about the asymptotic equivalence
of the two absorption events follows by subtracting their
probabilities, since the unit-ball event is contained in
the prescribed-point event.
\end{proof}

\section{Intersections of two independent Brownian-motion convex hulls}
\label{sec:two-brownian-hull-intersections}

In this section we prove
Theorems~\ref{thm:intro-two-brownian-hulls-fixed-dimension}
and~\ref{thm:intro-two-brownian-first-intersection-high-dimension}.
The common ingredient is an exact identity relating Gaussian-point
absorption by the joint convex hull of two Brownian paths to
origin absorption by a stationary Ornstein--Uhlenbeck process.
We first establish this identity and compare joint-hull absorption
with intersection of the two hulls. We then consider the long-time
non-intersection probability in a fixed dimension and the first
intersection time as the dimension tends to infinity.

\smallskip
Throughout, $W^{(1)}$ and $W^{(2)}$ are independent standard Brownian
motions in $\R^d$, both starting at the origin. Their hulls
$\mathcal K_{d,T}^{B,(i)}$ and the first intersection time $H_d$
of $\mathcal K_{d,T}^{B,(1)}$ and
$e_1+\mathcal K_{d,T}^{B,(2)}$ are defined in
\eqref{eq:intro-two-brownian-hulls}
and~\eqref{eq:intro-two-brownian-first-intersection-time}.

\subsection{Gaussian-point absorption by the joint convex hull of two Broanian motions}
\label{subsec:two-brownian-gaussian-point-absorption}

For $T_1,T_2>0$, define
\begin{equation}
\label{eq:two-brownian-joint-hull}
\mathcal C_{d;T_1,T_2}
:=\conv\bigl(
\mathcal K_{d,T_1}^{B,(1)}
\cup(-\mathcal K_{d,T_2}^{B,(2)})\bigr).
\end{equation}
By symmetry, this has the same distribution as the joint convex
hull without reflection of the second path. The reflected version
will be convenient for comparing absorption with hull intersection.
Let $Z_d\sim N(0,I_d)$ be independent of both Brownian motions, and put
\begin{equation}
\label{eq:two-brownian-gaussian-point-absorption-definition}
q_{d;T_1,T_2}^{(2)}(\sigma)
:=\Pp(\sigma Z_d\in\mathcal C_{d;T_1,T_2}),
\qquad \sigma\ge0.
\end{equation}
For unit time intervals, abbreviate
$\mathcal C_d:=\mathcal C_{d;1,1}$ and
$q_d^{(2)}:=q_{d;1,1}^{(2)}$, and define
$p_d^{(2)}(r):=\Pp(re_1\in\mathcal C_d)$ for $r\ge0$.

We also recall the ordinary Ornstein--Uhlenbeck absorption
probability from the companion paper~\cite{KabluchkoOUSpirals}.
Let $(U_u)_{u\in\R}$ be a centered stationary Ornstein--Uhlenbeck
process in $\R^d$ with continuous paths, independent coordinates,
and coordinate covariance
$\Ee[U_u^{(i)}U_v^{(i)}]=\eee^{-|u-v|/2}$, and write
\begin{equation}
\label{eq:two-brownian-ou-absorption-definition}
a_d^B(L):=\Pp\bigl(0\in\conv\{U_u:0\le u\le L\}\bigr),
\qquad L>0.
\end{equation}

\begin{lemma}[Gaussian-point absorption by two Brownian paths]
\label{lem:two-brownian-gaussian-point-absorption}
For every integer $d\ge1$, every $T_1,T_2>0$, and every $\sigma>0$,
\begin{align}
q_{d;T_1,T_2}^{(2)}(\sigma)
&=a_d^B\left(
\log\left(1+\frac{T_1}{\sigma^2}\right)
+\log\left(1+\frac{T_2}{\sigma^2}\right)
\right)
\label{eq:two-brownian-gaussian-point-OU-identity}
\\
&=q_d^B\left(
\frac{\sigma^2}{\sqrt{\sigma^2(T_1+T_2)+T_1T_2}}
\right),
\label{eq:two-brownian-gaussian-point-single-hull-identity}
\end{align}
where $q_d^B$ is defined in
\eqref{eq:gaussian-point-absorption-functions}.
In particular,
\begin{equation}
\label{eq:two-brownian-unit-time-gaussian-point-identity}
q_d^{(2)}(\sigma)
=a_d^B\bigl(2\log(1+\sigma^{-2})\bigr)
=q_d^B\left(\frac{\sigma^2}{\sqrt{1+2\sigma^2}}\right).
\end{equation}
\end{lemma}

\begin{proof}
Put $L_i:=\log(1+T_i/\sigma^2)$ for $i\in\{1,2\}$.
The event $\{\sigma Z_d\in\mathcal C_{d;T_1,T_2}\}$ is
equivalent to the origin belonging to the convex hull of the two
families
\[
\{\sigma Z_d-W_t^{(1)}:0\le t\le T_1\},
\qquad
\{\sigma Z_d+W_s^{(2)}:0\le s\le T_2\}.
\]
Define a process on $[-L_2,L_1]$ by
\[
\widetilde U_u:=
\begin{cases}
\displaystyle
\frac{\sigma Z_d-W_{\sigma^2(\eee^u-1)}^{(1)}}
{\sigma\eee^{u/2}},&0\le u\le L_1,\\[7pt]
\displaystyle
\frac{\sigma Z_d+W_{\sigma^2(\eee^{-u}-1)}^{(2)}}
{\sigma\eee^{-u/2}},&-L_2\le u\le0.
\end{cases}
\]
Both definitions agree at $u=0$. This is a centered Gaussian
process with continuous paths and independent coordinates.
For $0\le u\le v\le L_1$, its coordinate covariance is
\[
\Ee[\widetilde U_u^{(i)}\widetilde U_v^{(i)}]
=\frac{\sigma^2+\sigma^2(\eee^u-1)}
{\sigma^2\eee^{(u+v)/2}}
=\eee^{-(v-u)/2}.
\]
The same calculation applies when $-L_2\le u\le v\le0$.
When $-L_2\le u\le0\le v\le L_1$, only the common Gaussian
vector $Z_d$ contributes to the covariance, giving
\[
\Ee[\widetilde U_u^{(i)}\widetilde U_v^{(i)}]
=\frac{\sigma^2}{\sigma^2\eee^{(v-u)/2}}
=\eee^{-(v-u)/2}.
\]
Thus $\widetilde U$ has the law of the stationary
Ornstein--Uhlenbeck process restricted to $[-L_2,L_1]$.

Multiplying individual vectors by positive scalars preserves
whether their convex hull contains the origin: one rescales the
coefficients in a finite convex combination and normalizes their
sum. Applying this observation to the two families above, and
then using stationarity, gives
\[
q_{d;T_1,T_2}^{(2)}(\sigma)
=\Pp\bigl(0\in\conv\{\widetilde U_u:-L_2\le u\le L_1\}\bigr)
=a_d^B(L_1+L_2).
\]
This proves \eqref{eq:two-brownian-gaussian-point-OU-identity}.

The ordinary OU--Brownian correspondence
from~\cite[Theorem~5.8]{KabluchkoOUSpirals} states that
$q_d^B(\rho)=a_d^B(\log(1+\rho^{-2}))$ for $\rho>0$.
Taking
$\rho=\sigma^2/\sqrt{\sigma^2(T_1+T_2)+T_1T_2}$ gives
\[
1+\rho^{-2}
=\left(1+\frac{T_1}{\sigma^2}\right)
 \left(1+\frac{T_2}{\sigma^2}\right).
\]
Hence \eqref{eq:two-brownian-gaussian-point-single-hull-identity}
follows, and \eqref{eq:two-brownian-unit-time-gaussian-point-identity}
is its specialization to $T_1=T_2=1$.
\end{proof}

\medskip
\noindent\emph{Comparison with intersections of the two hulls.}
The preceding lemma concerns the convex hull of a union, whereas
intersection of the two hulls is described by their Minkowski
difference. For compact convex sets $A,B\subset\R^d$ containing
$0$, write $A-B:=\{a-b:a\in A,\ b\in B\}$. Then
\begin{equation}
\label{eq:two-brownian-joint-hull-Minkowski-inclusions}
\conv(A\cup(-B))\subseteq A-B
\subseteq2\conv(A\cup(-B)).
\end{equation}
Indeed, $A$ and $-B$ are contained in the convex set $A-B$,
while $(a-b)/2\in\conv(A\cup(-B))$ for every $a\in A$ and
$b\in B$. Moreover,
$A\cap(x+B)\ne\varnothing$ if and only if $x\in A-B$.

Apply these observations to
$A=\mathcal K_{d,T}^{B,(1)}$ and
$B=\mathcal K_{d,T}^{B,(2)}$.
Brownian scaling gives
$\mathcal C_{d;T,T}\overset{\mathrm d}=\sqrt T\,\mathcal C_d$.
Consequently, for every $T>0$,
\begin{align}
p_d^{(2)}(T^{-1/2})
&\le\Pp(H_d\le T)
\le p_d^{(2)}\bigl((2\sqrt T)^{-1}\bigr),
\label{eq:two-brownian-intersection-absorption-comparison}
\\
1-p_d^{(2)}\bigl((2\sqrt T)^{-1}\bigr)
&\le\Pp(H_d>T)
\le1-p_d^{(2)}(T^{-1/2}).
\label{eq:two-brownian-nonintersection-nonabsorption-comparison}
\end{align}
These inequalities reduce both asymptotic problems to
prescribed-point absorption by $\mathcal C_d$.

Finally, $\mathcal C_d$ is rotationally invariant in distribution,
convex, and contains $0$. Hence $p_d^{(2)}$ is nonincreasing, and
conditioning on the independent Gaussian vector gives
\begin{equation}
\label{eq:two-brownian-gaussian-radial-averaging}
q_d^{(2)}(\sigma)
=\Ee\bigl[p_d^{(2)}(\sigma R_d)\bigr],
\qquad R_d:=\|Z_d\|,\quad \sigma\ge0.
\end{equation}
This is the analogue of \eqref{eq:gaussian-radial-averaging}
for the joint hull and will be used in both proofs below.

\subsection{Non-intersection probabilities in fixed dimensions}
\label{subsec:two-brownian-hull-intersections-fixed-dim}

We first determine the small-radius non-absorption probability
for the joint hull $\mathcal C_d$, using the Gaussian-point identity
from the preceding subsection and
Lemma~\ref{lem:radial-gaussian-regular-variation}.
The comparison with the Minkowski difference then yields the
bounds for the first intersection time.

\begin{proof}[Proof of Theorem~\ref{thm:intro-two-brownian-hulls-fixed-dimension}]
Fix an integer $d\ge1$, and put
$g_d(r):=1-p_d^{(2)}(r)=\Pp(re_1\notin\mathcal C_d)$ for $r\ge0$.
Since $\mathcal C_d$ is convex and contains $0$, the function
$g_d$ is nondecreasing and takes values in $[0,1]$.
By \eqref{eq:two-brownian-gaussian-radial-averaging}
and~\eqref{eq:two-brownian-unit-time-gaussian-point-identity},
\begin{equation}
\label{eq:two-brownian-joint-hull-gaussian-nonabsorption}
\Ee[g_d(\sigma R_d)]
=1-q_d^{(2)}(\sigma)
=1-a_d^B\bigl(2\log(1+\sigma^{-2})\bigr),
\qquad \sigma>0,
\end{equation}
where $R_d$ has the $\chi_d$ distribution.
The fixed-dimensional OU asymptotics
from~\cite[Theorem~5.12]{KabluchkoOUSpirals} give
\[
1-a_d^B(L)
\sim\frac{2^{2-d}}{\pi(d-1)!}\,
L^{d-1}\eee^{-L/2},
\qquad L\to\infty.
\]
Substituting $L=2\log(1+\sigma^{-2})$, and noting that
$\eee^{-L/2}=\sigma^2/(1+\sigma^2)$, we obtain
\begin{equation}
\label{eq:two-brownian-joint-hull-gaussian-small-scale}
\Ee[g_d(\sigma R_d)]
\sim\frac{2}{\pi(d-1)!}\,
\sigma^2\left(\log\frac1{\sigma^2}\right)^{d-1},
\qquad \sigma\downarrow0.
\end{equation}

The Gaussian radial average in
\eqref{eq:two-brownian-joint-hull-gaussian-small-scale}
is positive for sufficiently small $\sigma$ and regularly
varying at zero with index $2$.
We may therefore apply
Lemma~\ref{lem:radial-gaussian-regular-variation}
with $f=g_d$ and $\beta=2$.
Since $\Ee[R_d^2]=d$, the lemma yields
\begin{equation}
\label{eq:two-brownian-joint-hull-small-radius-nonabsorption}
g_d(r)\sim\frac1d\Ee[g_d(rR_d)]
\sim\frac{2}{\pi d!}\,
r^2\left(\log\frac1{r^2}\right)^{d-1},
\qquad r\downarrow0.
\end{equation}

Finally, \eqref{eq:two-brownian-nonintersection-nonabsorption-comparison}
gives
\[
g_d\bigl((2\sqrt T)^{-1}\bigr)
\le\Pp(H_d>T)\le g_d(T^{-1/2}),
\qquad T>0.
\]
By \eqref{eq:two-brownian-joint-hull-small-radius-nonabsorption},
\begin{align*}
g_d\bigl((2\sqrt T)^{-1}\bigr)
&\sim\frac{1}{2\pi d!}\,
\frac{(\log T)^{d-1}}{T},
\\
g_d(T^{-1/2})
&\sim\frac{2}{\pi d!}\,
\frac{(\log T)^{d-1}}{T},
\qquad T\to\infty.
\end{align*}
Combining these equivalents with the preceding inequalities
proves \eqref{eq:intro-two-brownian-hulls-nonintersection-bounds}.
\end{proof}

\subsection{First intersection times in high dimensions}
\label{subsec:two-brownian-hull-intersections-high-dim}

We combine the Gaussian-point identity
\eqref{eq:two-brownian-unit-time-gaussian-point-identity}
with the high-dimensional transition
\eqref{eq:type-B-high-dimensional-Gaussian-point-input}.
Concentration of the Gaussian radius and the comparison between
joint-hull absorption and hull intersection then give the desired
limit. The fixed-dimensional tail estimates are not needed here.

\begin{proof}[Proof of Theorem~\ref{thm:intro-two-brownian-first-intersection-high-dimension}]
Recall that $0<H_d<\infty$ almost surely, as noted in
Subsection~\ref{subsec:intro-two-brownian-hull-intersections}.
Fix $y\in\R$, and put
\[
T_d:=\eee^{d+y\sqrt d},
\qquad r_d:=T_d^{-1/2}.
\]
For every fixed $c>0$, set
$\sigma_d:=cr_d/\sqrt d$ and
$\rho_d:=\sigma_d^2/\sqrt{1+2\sigma_d^2}$.
By \eqref{eq:two-brownian-unit-time-gaussian-point-identity},
$q_d^{(2)}(\sigma_d)=q_d^B(\rho_d)$.
Since $1+\rho_d^{-2}=(1+\sigma_d^{-2})^2$, we have
\[
\frac{\log(1+\rho_d^{-2})-2d}{2\sqrt d}
=\frac{\log(1+dT_d/c^2)-d}{\sqrt d}\\
=y+\frac{\log d-2\log c}{\sqrt d}+o(1)
\longrightarrow y.
\]
Applying \eqref{eq:type-B-high-dimensional-Gaussian-point-input}
therefore gives
\begin{equation}
\label{eq:two-brownian-critical-gaussian-point-limit}
q_d^{(2)}\!\left(\frac{cr_d}{\sqrt d}\right)
\longrightarrow\Phi(y),
\qquad c>0.
\end{equation}

We next remove the Gaussian radius. The random variable $R_d$
in \eqref{eq:two-brownian-gaussian-radial-averaging} satisfies
$R_d/\sqrt d\to1$ in probability, as in the proof of
Theorem~\ref{thm:intro-type-B-first-absorption-high-dimension}.
Fix $c>0$ and $0<\varepsilon<1$.
By \eqref{eq:two-brownian-gaussian-radial-averaging},
monotonicity of $p_d^{(2)}$, and $0\le p_d^{(2)}\le1$,
\begin{align*}
q_d^{(2)}\!\left(\frac{cr_d}{(1-\varepsilon)\sqrt d}\right)
-\Pp\bigl(R_d<(1-\varepsilon)\sqrt d\bigr)
&\le p_d^{(2)}(cr_d),
\\
p_d^{(2)}(cr_d)
&\le q_d^{(2)}\!\left(\frac{cr_d}{(1+\varepsilon)\sqrt d}\right)
+\Pp\bigl(R_d>(1+\varepsilon)\sqrt d\bigr).
\end{align*}
Both error probabilities tend to zero.
By \eqref{eq:two-brownian-critical-gaussian-point-limit},
applied with the fixed constants $c/(1-\varepsilon)$ and
$c/(1+\varepsilon)$, both Gaussian-point probabilities
converge to $\Phi(y)$. Hence
\begin{equation}
\label{eq:two-brownian-critical-prescribed-point-limit}
p_d^{(2)}(cr_d)\longrightarrow\Phi(y),
\qquad c>0.
\end{equation}

Finally,
\eqref{eq:two-brownian-intersection-absorption-comparison}
gives
\[
p_d^{(2)}(r_d)
\le\Pp(H_d\le T_d)
\le p_d^{(2)}(r_d/2).
\]
Taking $c=1$ and $c=1/2$ in
\eqref{eq:two-brownian-critical-prescribed-point-limit},
we conclude that
\[
\Pp\!\left(H_d\le\eee^{d+y\sqrt d}\right)
\longrightarrow\Phi(y).
\]
This proves
\eqref{eq:intro-two-brownian-first-intersection-Gaussian-window}.
Since $y\in\R$ was arbitrary and $\Phi$ is continuous, it also
proves the central limit theorem
\eqref{eq:intro-two-brownian-first-intersection-logarithmic-CLT}.

In particular, tightness of $(\log H_d-d)/\sqrt d$ implies
$d^{-1}\log H_d\to1$ in probability, and hence
$H_d^{1/d}\to\eee$ in probability.
\end{proof}

\section*{Acknowledgements}
This work was supported by the DFG under Germany's Excellence Strategy
EXC 2044/2--390685587, Mathematics M\"unster:
\emph{Dynamics-Geometry-Structure}; by the DFG Priority Programme
SPP 2265, \emph{Random Geometric Systems}; and by the DFG Research
Training Group \emph{Rigorous Analysis of Complex Random Systems}
(RTG 3027, Project Number 524444762).

\section*{Declarations}

\subsection*{Statement on the use of generative AI}
This paper and our companion paper~\cite{KabluchkoOUSpirals}
grew out of an extended dialogue between the human author
and ChatGPT (OpenAI, GPT-6).
One of the starting points was the approach
of~\cite{KabluchkoZaporozhets2020}, where absorption
probabilities were studied for Gaussian polytopes, that is,
convex hulls of independent standard Gaussian points in $\R^d$.

ChatGPT played a key role throughout the mathematical development
and writing of the present paper.
Its contributions included proposing and refining mathematical
ideas, carrying out technical calculations, suggesting proof
strategies and generalizations, deriving and checking formulas,
identifying connections with related areas of probability and
convex geometry, researching relevant literature, and assisting
with the organization and writing of the manuscript.
The work developed iteratively: suggestions generated by ChatGPT
were examined, corrected, refined, or omitted from the final
manuscript, while the human author supplied mathematical direction,
additional ideas, research questions, and editorial decisions
throughout. Mathematical arguments and claims included in the final manuscript were critically examined and independently verified by the human author, who takes full responsibility for the contents of the paper.

\subsection*{Conflict of interest statement}
The author declares that he has no conflict of interest.

\subsection*{Data availability statement}
We do not analyze or generate any datasets.

\bibliographystyle{amsplain}
\bibliography{brownian_convex_hulls_d_dim_bib}

\end{document}